\documentclass[11pt]{article}
\usepackage{amssymb, amsthm, amsmath, amscd}
\newtheorem{thm}{Theorem}[section]
\newtheorem{lem}[thm]{Lemma}
\newtheorem{prop}[thm]{Proposition}
\newtheorem{cor}[thm]{Corollary}
\newtheorem{NN}[thm]{}
\theoremstyle{definition}\newtheorem{df}[thm]{Definition}
\theoremstyle{definition}\newtheorem{rem}[thm]{Remark}
\theoremstyle{definition}

\renewcommand{\phi}{\varphi}

\newcommand{\Z}{\mathbb{Z}}
\newcommand{\Q}{\mathbb{Q}}
\newcommand{\R}{\mathbb{R}}
\newcommand{\C}{\mathbb{C}}
\newcommand{\T}{\mathbb{T}}

\newcommand{\morp}{contractive completely positive linear map}

\newcommand{\hm}{homomorphism}
\newcommand{\dt}{\delta}
\newcommand{\ep}{\epsilon}
\newcommand{\andeqn}{\,\,\,{\rm and}\,\,\,}
\newcommand{\rforal}{\,\,\,{\rm for\,\,\,all}\,\,\,}
\newcommand{\CA}{$C^*$-algebra}
\newcommand{\SCA}{$C^*$-subalgebra}

\newcommand{\af}{{\alpha}}
\newcommand{\bt}{{\beta}}

\newcommand{\D}{\mathbb D}
\newcommand{\beq}{\begin{eqnarray}}
\newcommand{\eneq}{\end{eqnarray}}
\newcommand{\tforal}{\,\,\,\text{for\,\,\,all}\,\,\,}

\title{Asymptotically Unitary Equivalence and Asymptotically Inner Automorphisms
}
\author{Huaxin Lin
 }
\date{}

\begin{document}

\maketitle

\begin{abstract}
Let $C$ be a unital AH-algebra and let $A$ be a unital separable
simple \CA\, with tracial rank zero. Suppose that $\phi_1, \phi_2:
C\to A$ are two unital monomorphisms. We show that there is a
continuous path of unitaries $\{u_t: t\in [0, \infty)\}$ of $A$
such that
$$
\lim_{t\to\infty}u_t^*\phi_1(a)u_t=\phi_2(a)\tforal a\in C
$$
if and only if $[\phi_1]=[\phi_2]$ in $KK(C,A),$ $\tau\circ
\phi_1=\tau\circ \phi_2$ for all $\tau\in T(A)$ and a rotation
related map ${\tilde\eta}_{\phi_1,\phi_2}$ associated with $\phi_1$
and $\phi_2$ is zero.  In particular, an automorphism $\af$ on a
unital separable simple \CA\, $A$ in ${\cal N}$ with tracial rank
zero is asymptotically inner if and only if
$$
[\af]=[{\rm id}_A]\,\,\,\text{in}\,\,\,KK(A,A)
$$
and the rotation related map ${\tilde\eta}_{\phi_1, \phi_2}$ is
zero.

Let $A$ be a unital AH-algebra (not necessarily simple) and let
$\af\in Aut(A)$ be an automorphism. As an application, we show
that the associated crossed product $A\rtimes_{\af}\Z$ can be
embedded into a unital simple AF-algebra if and only if $A$ admits
a strictly positive $\af$-invariant tracial state.

\end{abstract}

\section{Introduction}
Given two compact metric spaces and two continuous maps from one
space to another, one of the fundamental questions is when these
two maps are equivalent in certain sense. Equivalently, one may
study \hm s from the (commutative) algebra of continuous functions
on one space to the (commutative) algebra of continuous functions
on the other. By the Gelfand transform, every separable commutative
\CA\, is isomorphic to $C(X)$ for some compact metric space. One
way to make this question (at least partially) non-commutative is
to replace the target algebra by a unital \CA.
  Let $X$ be
a compact metric space and let $B$ be a unital  \CA. Let $\phi_1,
\phi_2: C(X)\to B$ be two unital monomorphisms. A natural and
important question is when $\phi_1$ and $\phi_2$ are unitarily
equivalent, i.e., when there is a unitary $u\in B$ such that ${\rm
ad}\, u\circ \phi_1=\phi_2?$ A classical special case is the case
that $B$ is the Calkin algebra, i.e., $B=B(H)/{\cal K},$ where $H$
is a separable infinite dimensional Hilbert space and ${\cal K}$
is the compact operators on $H.$ This problem originated from the
study of essentially normal operators on the Hilbert space. It is
closely related to perturbation of normal operators and index
theory. This classical problem was  solved by Brown, Douglass and
Fillmore (\cite{BDF1} and \cite{BDF2}). Later it developed into
the Kasparov KK-theory and had a profound impact on many related
fields.

However, recently it becomes apparent that the problem is
especially
 interesting when $B$ is a unital separable simple \CA\, (note
that the Calkin algebra is a unital simple \CA). It also becomes
clear that unitary equivalence must be given the way to some more
practical and useful relation. When $B$ is assumed to be a unital
separable simple \CA\, of tracial rank zero, it was shown in
\cite{Lncd} that $\phi_1$ and $\phi_2$ are approximately unitarily
equivalent, i.e., there exists a sequence of unitaries in $B$ such
that
$$
\lim_{n\to\infty}{\rm ad}\, u_n\circ \phi_1(f)=\phi_2(f)
$$
for all $f\in C(X),$ if and only if
$$
[\phi_1]=[\phi_2]\,\,\,\text{in}\,\,\,KL(C(X),A)\andeqn \tau\circ
\phi_1=\tau\circ \phi_2
$$
for all tracial states $\tau$ of $A.$ This result was developed
based on the results and techniques in the study of classification
of amenable \CA s and has several immediate applications including
in topological  dynamical systems.

A stronger relation is the asymptotic unitary equivalence. Let $A$
and $B$ be two unital \CA s and let $\phi_1,\phi_2: A\to B$ be two
unital \hm s. These two \hm s are asymptotically unitarily
equivalent if there is a {\it continuous} path of unitaries
$\{u_t: t\in [0,\infty)\}$ of $B$ such that
$$
\lim_{t\to\infty}{\rm ad}u_t\circ \phi_1(a)=\phi_2(a)
$$
for all $a\in A.$ Note that $\phi_1$ and $\phi_2$ are ``homotopic"
by a long path. So it is rather strong equivalence. But it still
weaker than that of unitary equivalence. It turns out that, at least in
several known occasions, one needs asymptotic unitary equivalence.
It is not difficult to see that in order to have $\phi_1$ and
$\phi_2$ to be asymptotically unitarily equivalent one must have
$[\phi_1]=[\phi_2]$ in $KK(A,B)$ instead of $[\phi_1]=[\phi_2]$ in
$KL(A,B).$ However, in addition of $\tau\circ \phi_1=\tau\circ
\phi_2$ for all tracial states $\tau$ of $B,$ it requires more to
obtain asymptotic unitary equivalence.

A pioneer work in this direction was done by Kishimoto and Kumjian
(\cite{KK1} and \cite{KK2}) who studied a special case that $A=B$
and $A$ is a unital simple $A\T$-algebra of real rank zero. Among
other things, they gave a $KK$-theoretical necessary and
sufficient condition for  an automorphisms $\af$ on a unital
simple  $A\T$-algebra of real rank zero to be asymptotically
inner. An interesting application was given by Matui. Using
Kishimoto and Kumjian's result, he  showed that if $A$ is a unital
simple  $A\T$-algebra of real rank zero and $\af\in Aut(A)$ is an
automorphism of $A$ then the associated crossed product
$A\rtimes_{\af}\Z$ can be embedded into a unital simple
AF-algebra(\cite{M}).

We study \hm s from an AH-algebra $A$ (not necessarily simple nor
real rank zero) to a unital separable simple \CA\, $B$ with tracial
rank zero. It should be noted that even in the simple case, a unital
simple AH-algebra may have higher stable rank and it may have some
other unexpected properties (see \cite{jV} and \cite{aT}). Suppose
that $\phi_1, \phi_2: A\to B$ are two unital \hm s. In order for
$\phi_1$ and $\phi_2$ to be asymptotically unitarily equivalent (or
approximately unitarily equivalent), $\phi_1$ and $\phi_2$ must have
the same kernel. Since any quotient of an AH-algebra is again an
AH-algebra, to simplify the matter, we may only consider
monomorphisms. Suppose that $\phi_1$ and $\phi_2$ are asymptotically
unitarily equivalent. Then, as mentioned above, $[\phi_1]=[\phi_2]$
in $KK(A,B)$ and $\tau\circ \phi_1=\tau\circ \phi_2$ for any tracial
state $\tau$ of $B.$ As in Kishimoto and Kumjian's case, a rotation
related map ${\tilde \eta}$ from $K_1(A)$ to $K_0(B)$  associated
with the pair of monomorphisms must vanish. Details of this rotation
map will be defined below (\ref{eta}). The main result of this paper
is to establish that these conditions are also sufficient for
$\phi_1$ and $\phi_2$ being asymptotically unitarily equivalent.

The main background tools to establish  the main theorem is the
similar result for approximate unitary equivalence established in
\cite{Lncd} and the so-called Basic Homotopy Lemmas
(\cite{Lnhomp}). Suppose that there is a unitary $u\in U(B)$ such
that
$$
\|[\phi_1(b), \, u]\|<\dt
$$
for some large set of elements $b\in A$ and sufficiently small
$\dt.$ One needs to know how to find a continuous path of
unitaries $u(t)$ with $u(0)=u,$ $u(1)=1_B$ and $\|[\phi_1(a),
u(t)]\|<\ep$ for some given set of elements $a\in A.$ The Basic
Homotopy Lemma (\cite{Lnhomp}) states that it is possible when a
map $\text{Bott}(\phi_1, u)$ vanishes. This type of theorem was
first given in \cite{BEEK}. The more general version used here is
given in \cite{Lnhomp}. We also  took many ideas from Kishimoto
and Kumjian's work (\cite{KK2}) for granted. It turns out that the
proof of the main theorem is considerably shorter in the case that
$K_i(A)$ ($i=0,1$) is finitely generated (\ref{LTM}). The general case needs a
number of additional facts which require some lengthy work.

Let $B$ be a unital separable amenable simple \CA\, with tracial
rank zero which satisfies the Universal Coefficient Theorem and
let $\af: B\to B$  be an endomorphism. By the main result
mentioned above, we show that $\af$ is asymptotically inner, i.e.,
there is a continuous path of unitaries $\{u_t: t\in [0,
\infty)\}$ of $A$ such that
$$
\lim_{t\to\infty}u_t^*au_t=\af(a)
$$
for all $a\in A$ if and only if $[\af]=[{\rm id}_A]$ in $KK(A,A)$
and the associated rotation related map ${\tilde
\eta}_{\phi_1,\phi_2}$ is zero.

Now we turn to quasidiagonal \CA s.  Quasidiagonality for crossed
products were studied by Voiculescu (\cite{V1}, \cite{V2} and
\cite{V3}). Let $X$ be a compact metric space and $\af$ be a
homeomorphism on $X.$ It was proved by Pimsner (\cite{Pi}) that
$C(X)\rtimes_{\af}\Z$ is quasidiagonal if and only if $\af$ is
pseudo-non-wondering and if and only if $C(X)\rtimes_{\af}\Z$  can
be embedded into an AF-algebra. When $A$ is an AF-algebra, Nate
Brown (\cite{Bn1}) proved that $A\rtimes_{\af}\Z$ is quasidiagonal
if and only if $A\rtimes_{\af}\Z$ can be embedded into an
AF-algebra. He also gave a $K$-theoretical necessary and
sufficient condition for $A\rtimes_{\af}\Z$ being embedded into an
AF-algebra. This result was generalized by Matui as mentioned
above.  See \cite{Bn2} for further discussion on this topic.

Let $A$ be a unital AH-algebra (not necessarily simple nor real
rank zero) and let $\af\in Aut(A)$ be an automorphism. We are
interested in the problem when the crossed product
$A\rtimes_{\af}\Z$ can be embedded into a unital simple AF-algebra
(\cite{Lnemb1} and \cite{Lnemb2}).

Suppose that there is a unital embedding $h: A\rtimes_{\af}\Z\to
B,$ where $B$ is a unital simple AF-algebra. Let $\tau\in T(B)$ be
a tracial state on $B.$ Then $\tau\circ h$ gives a faithful
$\af$-invariant tracial state on $A.$ As an application of the
main result of the paper, we prove that $A\rtimes_{\af}\Z$ can be
embedded into a unital simple AF-algebra if and only if $A$ admits
a faithful $\af$-invariant tracial state. Note that if $A$ is
simple then $A$ always admits a faithful $\af$-invariant tracial
state.

There is another  question concerning the main result. If $\phi_1$
and $\phi_2$ are asymptotically unitarily equivalent, can one
actually find a continuous path of unitaries $\{u_t: t\in [0,
\infty)\}$ of $B$ such that
$$
u_0=1_B\andeqn \lim_{t\to\infty}{\rm ad}\, u_t\circ
\phi_1(a)=\phi_2(a)\rforal a\in A?
$$
In other words, $\phi_1$ and $\phi_2$ are asymptotically unitarily
equivalent with $u_t\in U_0(A)$ for all $t\in [0, \infty).$

We will show, for example, when $A=C(X)$ and $B$ is a unital
separable simple \CA\, with tracial rank zero, the answer to the
question is affirmative (see \ref{SUC}). However, in general, this
can not be done even in the case that $A=M_n(C(X)).$ There is a
$K$-theory obstacle prevent one from choosing a path starting with
$1_B.$ This phenomenon will be discussed in this paper.

The paper is organized as follows:

Section 2 provides some notations and some background materials. In
Section 3, we define the rotation related map mentioned above and
give a generalized Exel's formula for bott elements concerning a
pair of almost commuting unitaries. In Section 4, we show that in
order for two unital monomorphisms $\phi_1$ and $\phi_2$ to be
asymptotically unitarily equivalent it is necessarily to have
$[\phi_1]=[\phi_2]$ in $KK(C,A)$ and $\tau\circ \phi_1=\tau\circ
\phi_2$ for all $\tau\in T(B)$ as well as the rotation related map
associated with these two \hm s vanishes. Section 5 discusses a pair
of almost commuting unitaries and associated bott element which puts
some known results into a more general setting. In section 6, using
classification theory for unital separable amenable simple \CA s of
tracial rank zero,  we show that for $A=C(X),$ where $X$ is a finite
CW complex $X,$ the  Bott map can be constructed from prescribed
information. In section 7, we generalize the results in section 6 to
the case that $A$ is a unital AH-algebra.  In section 8, we prove a
special case of the main theorem. In section 9, we present a stable
version of the so-called Basic Homotopy Lemma and use it to
establish certain facts that are needed in the proof of the main
theorem. Section 10 is devoted to the proof of the main theorem
(\ref{TM}). In section 11, we give a couple of applications. In
particular, we solve the problem when $A\rtimes_{\af}\Z$ can be
embedded into a unital simple AF-algebra for all unital AH-algebras.
Finally, in the last section, we discuss the problem when the
continuous paths in the definition of asymptotic unitary equivalence
can be taken in the connected component of the unitary groups which
contains the identity.


\section{Preliminaries}

\begin{NN}\label{NN1}
{\rm Let $A$ be a stably finite \CA. Denote by $T(A)$ the tracial
state space of $A$ and denote by $Aff(T(A))$ the space of all real
affine continuous functions on $T(A).$ Suppose $\tau\in T(A)$ is a
tracial state. We will also use $\tau$ for the trace $\tau\otimes
Tr$ on $A\otimes M_k=M_k(A)$ (for every integer $k\ge 1$), where
$Tr$ is the standard trace on $M_k.$

Define $\rho_A: K_0(A)\to Aff(T(A))$ be the positive \hm\, defined
by $\rho_A([p])(\tau)=\tau(p)$ for each projection $p$ in
$M_k(A).$

}
\end{NN}

\begin{NN}
{\rm  Let $A$ be a \CA. Denote by $SA=C_0((0,1), A)$ the
suspension of $A.$ }
\end{NN}

\begin{NN}

{\rm A \CA\, $A$ is an AH-algebra if $A=\lim_{n\to\infty}(A_n,
\psi_n),$ where each $A_n$ has the form $P_nM_{k(n)}(C(X_n))P_n,$
where $X_n$ is a finite CW complex (not necessarily connected) and
$P_n\in M_{k(n)}(C(X_n))$ is a projection.

We use $\psi_{n, \infty}: A_n\to A$ for the induced \hm. Note that
every separable commutative \CA\, is an AH-algebra. Every
AF-algebra and every $A\T$-algebra are AH-algebras.

}
\end{NN}

\begin{NN}
{\rm Denote by ${\cal N}$ the class of separable amenable \CA\,
which satisfies the Universal Coefficient Theorem. }

\end{NN}

\begin{NN}

{\rm

A unital separable simple \CA\, $A$ with tracial rank zero is
quasi-diagonal, of real rank zero, stable rank one and has weakly
unperforated $K_0(A).$ Every simple AH-algebra $A$ with real rank
zero, stable rank one and weakly unperforated $K_0(A)$ has tracial
rank zero (\cite{Lnah}).
 We refer to \cite{Lnplms}  and \cite{Lntr0} for the definition of tracial rank zero
and its properties. Unital separable simple \CA s in ${\cal N}$
with tracial rank zero can be classified by their $K$-theory (see
(\cite{Lnann} and \cite{Lnduke}).

}

\end{NN}

\begin{NN}
{\rm  Let $A$ and $B$ be two unital \CA s and let $\phi: A\to B$
be a \hm. One can extend $\phi$ to a \hm\, from $M_k(A)$ to
$M_k(B)$ by $\phi\otimes {\rm id}_{M_k}.$ In what follows we will
use $\phi$ again for this extension without further notice.}

\end{NN}

\begin{NN}

{\rm Let $A$ and $B$ be two \CA s and let $L_1, L_2: A\to B$ be a
map.  Suppose that ${\cal F}\subset A$ is a subset and  $\ep>0.$
We write
$$
L_1\approx_{\ep} L_2\,\,\,\text{on}\,\,\, {\cal F},
$$
if
$$
\|L_1(a)-L_2(a)\|<\ep\tforal a\in {\cal F}.
$$

}

\end{NN}

\begin{NN}
{\rm Let $A$ be a unital \CA. Denote by $U(A)$ the group of
unitaries in $A$ and denote by $U_0(A)$ the path connected
component of $U(A)$ containing the identity.

By $Aut(A)$ we mean the group of automorphisms on $A.$ Let $u\in
U(A).$ We write ${\rm ad}\, u$ the inner automorphism defined by
${\rm ad}\, u(a)=u^*au$ for all $a\in A.$

}
\end{NN}

\begin{df}\label{Kund}
{\rm Let $A$ be a \CA.  Following Dadarlat and Loring (\cite{DL}),
denote
$$
\underline{K}(A)=\oplus_{i=0,1}
K_i(A)\bigoplus_{i=0,1}\bigoplus_{k\ge 2}K_i(A,\Z/k\Z).
$$
Let $B$ be a unital \CA. If furthermore, $A$ is assume to be
separable and satisfy the Universal Coefficient Theorem
(\cite{RS}), by \cite{DL},
$$
Hom_{\Lambda}(\underline{K}(A), \underline{K}(B))=KL(A,B).
$$

Here $KL(A,B)=KK(A,B)/Pext(K_*(A), K_*(B)).$ (see \cite{DL} for
details).

 Let $k\ge 1$ be an integer. Denote
$$
F_k\underline{K}(A)=\oplus_{i=0,1}K_i(A)\bigoplus_{n|k}K_i(A,\Z/k\Z).
$$
Suppose that $K_i(A)$ is finitely generated ($i=0,2$). It follows
from \cite{DL} that there is an integer $k\ge 1$ such that
\beq\label{dkl1}
Hom_{\Lambda}(F_k\underline{K}(A), F_k\underline{K}(B))
=Hom_{\Lambda}(\underline{K}(A), \underline{K}(B)).
\eneq
}
\end{df}

\begin{df}\label{Dbot2}
{\rm Let
$A$ and $B$ be  two unital \CA s.  Let $h: A\to B$ be a \hm\, and
$v\in U(B)$ such that
$$
h(g)v=vh(g)\,\rforal\, g\in A.
$$
 Thus we
obtain a \hm\, ${\bar h}: A\otimes C(S^1)\to B$ by ${\bar
h}(f\otimes g)=h(f)g(v)$ for $f\in A$ and $g\in C(S^1).$ The
tensor product induces two injective \hm s:
\beq\label{dbot01}
\bt^{(0)}&:& K_0(A)\to K_1(A\otimes C(S^1))\\
 \bt^{(1)}&:&
K_1(A)\to K_0(A\otimes C(S^1)).
\eneq
The second one is the usual Bott map. Note, in this way, one
writes
$$K_i(A\otimes C(S^1))=K_i(A)\oplus \bt^{(i-1)}(K_{i-1}(A)).$$
We use $\widehat{\bt^{(i)}}: K_i(A\otimes C(S^1))\to
\bt^{(i-1)}(K_{i-1}(A))$ for the projection to
$\bt^{(i-1)}(K_{i-1}(A)).$

For each integer $k\ge 2,$ one also obtains the following
injective \hm s:
\beq\label{dbot02}
\bt^{(i)}_k: K_i(A, \Z/k\Z))\to K_{i-1}(A\otimes C(S^1), \Z/k\Z),
i=0,1.
\eneq
Thus we write
\beq\label{dbot02-1}
K_{i-1}(A\otimes C(S^1), \Z/k\Z)=K_{i-1}(A,\Z/k\Z)\oplus
\bt^{(i)}_k(K_i(A, \Z/k\Z)),\,\,i=0,1.
\eneq
Denote by $\widehat{\bt^{(i)}_k}: K_{i}(A\otimes C(S^1),
\Z/k\Z)\to \bt^{(i-1)}_k(K_{i-1}(A,\Z/k\Z))$ similarly to that of
$\widehat{\bt^{(i)}}.,$ $i=1,2.$ If $x\in \underline{K}(A),$ we
use ${\boldsymbol{\beta}}(x)$ for $\bt^{(i)}(x)$ if $x\in K_i(A)$
and for $\bt^{(i)}_k(x)$ if $x\in K_i(A, \Z/k\Z).$ Thus we have a
map ${\boldsymbol{ \bt}}: \underline{K}(A)\to
\underline{K}(A\otimes C(S^1))$ as well as
$\widehat{\boldsymbol{\bt}}: \underline{K}(A\otimes C(S^1))\to
 {\boldsymbol{ \bt}}(\underline{K}(A)).$ Therefore one may write
 $\underline{K}(A\otimes C(S^1))=\underline{K}(A)\oplus {\boldsymbol{ \bt}}( \underline{K}(A)).$

On the other hand ${\bar h}$ induces \hm s ${\bar h}_{*i,k}:
K_i(A\otimes C(S^1)), \Z/k\Z)\to K_i(B,\Z/k\Z),$ $k=0,2,...,$ and
$i=0,1.$

We use $\text{Bott}(h,v)$ for all  \hm s ${\bar h}_{*i,k}\circ
\bt^{(i)}_k.$ We write
$$
\text{Bott}(h,v)=0,
$$
if ${\bar h}_{*i,k}\circ \bt^{(i)}_k=0$ for all $k\ge 1$ and
$i=0,1.$

We will use $\text{bott}_1(h,v)$ for the \hm\, ${\bar
h}_{1,0}\circ \bt^{(1)}: K_1(A)\to K_0(B),$ and
$\text{bott}_0(h,u)$ for the \hm\, ${\bar h}_{0,0}\circ \bt^{(0)}:
K_0(A)\to K_1(B).$

Since $A$ is unital, if $\text{bott}_0(h,v)=0,$ then $[v]=0$ in
$K_1(B).$

In what follows, we will use $z$ for the standard generator of
$C(S^1)$ and we will often identify $S^1$ with the unit circle
without further explanation. With this identification $z$ is the
identity map from the circle to the circle.

}
\end{df}

\begin{NN}\label{ddbot}
{\rm Given a finite subset ${\cal P}\subset \underline{K}(A),$
there exists a finite subset ${\cal F}\subset A$ and $\dt_0>0$
such that
$$
\text{Bott}(h, v)|_{\cal P}
$$
is well defined, if
$$
\|[h(a),\, v]\|=\|h(a)v-vh(a)\|<\dt_0\tforal a\in {\cal F}
$$
(see 2.10 of \cite{Lnhomp}).

There is $\dt_1>0$ (\cite{Lo}) such that $\text{bott}_1(u,v)$ is
well defined for any pair of unitaries $u$ and $v$ such that
$\|[u,\, v]\|<\dt_1.$ As in 2.2 of \cite{ER}, if $v_1,v_2,...,v_n$
are unitaries such that
$$
\|[u, \, v_j]\|<\dt_1/n,\,\,\,j=1,2,...,n,
$$
then
$$
\text{bott}_1(u,\,v_1v_2\cdots v_n)=\sum_{j=1}^n\text{bott}_1
(u,\, v_j).
$$

By considering  unitaries $z\in  {\widetilde{A\otimes C}}$
($C=C_n$ for some commutative \CA\, with torsion $K_0$ and
$C=SC_n$), from the above, for a given unital separable
\CA\, $A$ and a given finite subset ${\cal P}\subset
\underline{K}(A),$ one obtains a universal constant $\dt>0$ and a
finite subset ${\cal F}\subset A$ satisfying the following:
\beq\label{ddbot-1}
\text{Bott}(h,\, v_j)|_{\cal P}\,\,\, \text{is well
defined}\andeqn \text{Bott}(h,\, v_1v_2\cdots v_n)=\sum_{j=1}^n
\text{Bott}(h,\, v_j),
\eneq
for any unital \hm\, $h$ and unitaries $v_1, v_2,...,v_n$ for
which
\beq\label{ddbot-2}
\|[h(a),\, v_j]\|<\dt/n,\,\,\,j=1,2,...,n
\eneq
for all $a\in {\cal F}.$

If furthermore, $K_i(A)$ is finitely generated, then (\ref{dkl1})
holds. Therefore, there is a finite subset ${\cal Q}\subset
\underline{K}(A),$ such that
$$
\text{Bott}(h,v)
$$
is well defined if $\text{Bott}(h, v)|_{\cal Q}$ is well defined
(see also 2.3 of \cite{Lnhomp}).

See section 2 of \cite{Lnhomp} for the further information. }

\end{NN}

We will use  the following the theorems frequently.

\begin{thm}{\rm (Corollary 17.9 of \cite{Lnhomp})}\label{LNHOMP}
Let $A$ be a unital AH-algebra and let $\ep>0$ and ${\mathcal
F}\subset A$ be a finite subset. Suppose that $B$ is a unital
separable simple \CA\, with tracial rank zero and $h: A\to B$ is a
unital monomorphism. Then there exists $\dt>0,$ a finite subset
${\mathcal G}\subset A$ and a finite subset ${\mathcal P}\subset
\underline{K}(A)$ satisfying the following: Suppose that there is
a unitary $u\in B$ such that
\beq\label{NCT1C-1}
\|[h(a), u]\|<\dt\tforal f\in {\mathcal G}\andeqn
\rm{Bott}(h,u)|_{{\mathcal P}}=0.
\eneq
 Then there exists a
continuous path of unitaries $\{u_t:t\in [0,1]\}$ such that
$$
u_0=u,\,\,\, u_1=u, \|[h(a), v_t]\|<\ep\rforal f\in {\mathcal
F}\andeqn t\in [0,1],
$$
$$
\|u_t-u_{t'}\|\le (2\pi+\ep)|t-t'|\tforal t,t'\in [0,1] \andeqn
$$
$$
\rm{Length}(\{u_t\})\le 2\pi+\ep.
$$

\end{thm}

We will also use the following

\begin{thm}{\rm (Theorem 3.4 of \cite{Lncd}; also Theorem 3.6 of \cite{Lnemb2})}\label{CDT}
Let $C$ be a unital AH-algebra and let $B$ be a unital separable
simple \CA\, with tracial rank zero. Suppose that $\phi_1, \phi_2:
C\to B$ are two unital monomorphisms. Then there exists a sequence
of unitaries $\{u_n\}$ of $A$ such that
$$
\lim_{n\to\infty}{\rm ad}\, u_n\circ \phi_1(a)=\phi_2(a)\tforal
a\in C
$$
if and only if
$$
[\phi_1]=[\phi_2]\,\,\,\text{in}\,\,\, KL(C,B)\andeqn \tau\circ
\phi_1=\tau\circ \phi_2
$$
for all $\tau\in T(B).$
\end{thm}

\section{Rotation maps and  Exel's trace formula}

\begin{NN}\label{Mtorus}
{\rm Let $A$ and $B$ be two unital \CA s. Suppose that $\phi, \psi:
A\to B$ are two monomorphisms. Define
\beq\label{Dmt-1}
M_{\phi, \psi}=\{x\in C([0,1], B): x(0)=\phi(a)\andeqn
x(1)=\psi(a) \,\,\,\text{for some}\,\,\, a\in A\}.
\eneq
When $A=B$ and $\phi={\rm id}_A,$ $M_{\phi, \psi}$ is the usual
mapping torus.  We may call $M_{\phi, \psi}$ the (generalized)
mapping torus of $\phi$ and $\psi.$

 This notation
will be used throughout of this article.  Thus one obtains an exact
sequence:
\beq\label{Dmt-2}
0\to SB \stackrel{\imath}{\to} M_{\phi, \psi}\stackrel{\pi_0}{\to}
A\to 0.
\eneq

Suppose that $A$ is a separable amenable \CA.  From (\ref{Dmt-2}),
one obtains an element in $Ext(A,SB).$ In this case we identify
$Ext(A,SB)$ with $KK^1(A,SB)$ and $KK(A,B).$

Suppose that $[\phi]=[\psi]$ in $KL(A,B).$ The mapping torus
$M_{\phi, \psi}$ corresponds a trivial element in $KL(A, B).$ It
follows that there are two exact sequences:
\beq\label{Dmt-3}
&&0\to K_1(B)\stackrel{\imath_*}{\to} K_0(M_{\phi, \psi})\stackrel{(\pi_0)_*}{\to} K_0(A)\to
0\andeqn\\\label{Dmt-4}
&&0\to K_0(B)\stackrel{\imath_*}{\to}K_1(M_{\phi, \psi}) \stackrel{(\pi_0)_*}{\to}
K_1(A)\to 0.
\eneq
which are pure extensions of abelian groups.

 }

\end{NN}

\begin{df}\label{Dr}
{\rm Suppose that $T(B)\not=\emptyset.$  Let $u\in M_l(M_{\phi,
\psi})$ be a unitary which is a piecewise smooth function on
$[0,1].$ For each $\tau\in T(B),$ we denote by $\tau$  the trace
$\tau\otimes Tr$ on $M_l(B),$ where $Tr$ is the standard trace on
$M_l$ as in \ref{NN1}. Define
\beq\label{Dr-1}
R_{\phi,\psi}(u)(\tau)={1\over{2\pi i}}\int_0^1
\tau({du(t)\over{dt}}u(t)^*)dt.
\eneq

Since
$$
\tau({du(t)\over{dt}}u(t)^*)=-\tau(u(t){du(t)^*\over{dt}}),
$$
$R_{\phi, \psi}(u)(\tau)$ is real.

We now assume that $[\phi]=[\psi]$ in $KL(A,B).$ We also assume
that
\beq\label{Dr-2}
\tau(\phi(a))=\tau(\psi(a))\rforal a\in A\andeqn \tau\in T(B).
\eneq

Exactly as in section 2 of \cite{KK2}, one has the following
statement:}

{\it When $[\phi]=[\psi]$ in $KL(A,B)$ and {\rm (\ref{Dr-2})}
holds, there exists a \hm\,
$$
R_{\phi, \psi}: K_1(M_{\phi, \psi})\to Aff(T(B))
$$
defined by
$$
R_{\phi, \psi}([u])(\tau)={1\over{2\pi i}}\int_0^1
\tau({du(t)\over{dt}}u(t)^*)dt.
$$
}


\end{df}

If $p$ is a projection in $M_l(B)$ for some integer $l\ge 1,$ one
has $\imath_*([p])=[u],$ where $u\in M_{\phi, \psi}$ is a unitary
defined by
$$
u(t)=e^{2\pi it}p+(1-p)\,\,\,\text{for}\,\,\, t\in [0,1].
$$
It follows that
$$
R_{\phi, \psi}(\imath_*([p]))(\tau)=\tau(p)\tforal \tau\in T(B).
$$
In other words,
$$
R_{\phi, \psi}(\imath_*([p]))=\rho_B([p]).
$$
Thus one has, exactly as in 2.2 of \cite{KK2},  the following:

\begin{lem}\label{DrL} When $[\phi]=[\psi]$ in $KL(A,B)$ and
{\rm (\ref{Dr-2})} holds, the following diagram commutes:
$$
\begin{array}{ccccc}
K_0(B) && \stackrel{\imath_*}{\longrightarrow} && K_1(M_{\phi, \psi})\\
& \rho_B\searrow && \swarrow R_{\phi,\psi} \\
& & Aff(T(B)) \\
\end{array}
$$

\end{lem}

\vspace{0.2in}

\begin{df}\label{eta}
If furthermore, $[\phi]=[\psi]$ in $KK(A,B)$ and $A$ satisfies the
Universal Coefficient Theorem, using Dadarlat-Loring's notation,
one has the following splitting exact sequence:
\beq\label{eta-1-}
0\to \underline{K}(SB)\,{\stackrel{[\imath]}{\to}}\,
\underline{K}(M_{\phi,\psi})\,{\stackrel{[\pi_0]}{\rightleftarrows}}_{\theta}
\,\,\underline{K}(A)\to 0.
\eneq
In other words there is  $\theta\in Hom_{\Lambda}(\underline{K}(A),
\underline{K}(M_{\phi,\psi}))$ such that $[\pi_0]\circ \theta=[\rm
id_A].$ In particular, one has a monomorphism $\theta|_{K_1(A)}:
K_1(A)\to K_1(M_{\phi, \psi})$ such that $[\pi_0]\circ
\theta|_{K_1(A)}=({\rm id}_A)_{*1}.$ Thus, one may write
\beq\label{eta-1}
K_1(M_{\phi, \psi})=K_0(B)\oplus K_1(A).
\eneq
Suppose also that $\tau\circ \phi=\tau\circ \psi$ for all $\tau\in
T(A).$  Then  one obtains the \hm\,
\beq\label{eta-2}
R_{\phi,\psi}\circ \theta|_{K_1(A)}: K_1(A)\to Aff(T(B)).
\eneq

We say a  rotation related map  vanishes,  if there exists a such
splitting map $\theta$ such that
$$R_{\phi, \psi}\circ \theta=0.$$

Note that there is nothing standard about the splitting map
$\theta.$ To keep the same notation as in \cite{KK2}, when
$R_{\phi,\psi}\circ \theta=0$ for some such $\theta,$ we write
$${\tilde \eta}_{\phi,\psi}=0.$$

In what follows  whenever we write ${\tilde \eta}_{\phi, \psi}=0,$
we mean  $\theta(K_1(A))\in {\rm ker}R_{\phi,\psi}$ for some
$\theta$ so that (\ref{eta-1-}) holds. Here we do not define a map
${\tilde \eta}_{\phi,\psi}$ but simply use
$\tilde{\eta}_{\phi,\psi}=0$ as a notation.

Thus, $\theta$ also gives the following:
$$
{\rm ker}R_{\phi,\psi}={\rm ker}\rho_B\oplus K_1(A).
$$

Under the assumption that $[\phi]=[\psi]$ in $KK(A,B)$ and
$\tau\circ \phi=\tau\circ \psi $ for all $\tau\in T(B),$ and if, in
addition, $K_i(A)$ is torsion free, such $\theta$ exists whenever
$\rho_B(K_0(B))=R_{\psi,\phi}(K_1(M_{\phi,\psi}))$ and the following
splits:
$$
0\to {\rm ker}\rho_A\to {\rm ker}R_{\phi, \psi}\to K_1(A)\to 0.
$$

\end{df}

\vspace{0.2in}

\begin{lem}\label{Rotrange}
 Let $C$ be a unital separable amenable \CA\, and let $B$ be a
unital \CA. Suppose that $\phi, \psi: C\to B$ are two unital
monomorphisms such that $[\phi]=[\psi]$ in $KL(C,B)$ and
\beq\label{rot-1}
\tau\circ \phi=\tau\circ \psi\rforal \tau\in T(B).
\eneq

Suppose that $u\in M_l(C)$ is a unitary and $w\in U(M_l(B))$ such
that
\beq\label{rot-2}
\|\phi(u){\rm ad}\, w(\psi(u^*))-1\|<2.
\eneq

Then
\beq\label{rot-3}
{1\over{2\pi i}}\tau(\log (\phi(u){\rm ad}\, w(\psi(u^*))))\in
R_{\phi, \psi}(K_1(M_{\phi, \psi})).
\eneq

\end{lem}

\begin{proof}
By replacing $\phi$ and $\psi$ by $\phi\otimes{\rm id}_{M_l}$ and
by $\psi\otimes {\rm id}_{M_l},$ respectively, we may reduce the
general case to the case that $u\in U(C)$ and $w\in U(B).$ Note
that $\phi({\rm diag}(u,1))={\rm diag}(\phi(u),1)$ and $\psi({\rm
diag}(u,1))={\rm daig}(\psi(u),1).$ Also that
$$
\log(\phi({\rm diag}(u,1)){\rm ad}({\rm diag}(w,w^*))(\psi({\rm
diag}(u^*,1))))\\
={\rm diag}(\log(\phi(u){\rm ad}w(\psi(u^*))),0).
$$
Thus, by replacing $u$ by ${\rm diag}(u,1)$ and $w$ by ${\rm
diag}(w,w^*),$ we may assume that $w\in U_0(B).$

It follows that there are $b_1,b_2,...,b_m\in B_{s.a}$ with
$\|b_j\|\le 1$ such that
\beq\label{rot-4}
w=\prod_{j=1}^m e^{2\pi i b_j}.
\eneq
There is a smooth branch of $\log$ defined on the ${\rm
sp}(\phi(u){\rm ad}\, w(\psi(u^*)).$ Let
\beq\label{rot-5}
h={1\over{2\pi i}}\log(\phi(u){\rm ad}\, w(\psi(u^*))
\eneq
with $\|h\|\le 1.$
 Define a unitary $U\in M_{\phi, \psi}$ as
follows
\beq\label{rot-5+}
U(t)=\begin{cases}\phi(u^*)e^{2(2m+1)\pi i h
t}\,\,\,&\text{if}\,\,\, t\in
[0,{1\over{2m+1}})\\
e^{ -2\pi i(2-(2m+1)t) b_1}\prod_{j=2}^me^{-2\pi
ib_j}\psi(u^*)w\,\,\,
&\text{if}\,\,\, t\in [{1\over{2m+1}}, {2\over{2m+1}})\\
e^{-2\pi i (k+1-(2m+1)t)b_k}\prod_{j=k+1}^{m}e^{-2\pi i
b_j}\psi(u^*)w\,\,\,&\text{if}\,\,\, t\in [{k\over{2m+1}},
{k+1\over{2m+1}})\\
e^{-2\pi i (m+1-(2m+1)t)b_m}\psi(u^*)w \,\,\,&\text{if}\,\,\, t\in
[{m\over{2m+1}}, {m+1\over{2m+1}})\\
\psi(u^*)e^{2\pi i (m+1+k-(2m+1)t)b_k}\prod_{j=+1}^{m}e^{-2\pi i
b_j}\,\,\,&\text{if}\,\,\, t\in [{m+k\over{2m+1}},
{m+k+1\over{2m+1}})\\
\psi(u^*)e^{2\pi i (2m+1-(2m+1)t)b_m}\,\,\,&\text{if}\,\,\, t\in
[{2m\over{2m+1}}, {2m+1\over{2m+1}}]\,.\\
\end{cases}
\eneq
Note that
$$
U(0)=\phi(u^*)\andeqn U(1)=\psi(u^*).
$$
So $U(t)$ is indeed in $M_{\phi, \psi}.$ Moreover, it is piecewise
smooth function on $[0,1].$ It follows that
\beq\label{rot-6}
{1\over{2\pi i}}\int_0^1 \tau({dU(t)\over{dt}}U(t)^*)dt\in
R_{\phi, \psi}(K_1(M_{\phi, \psi})).
\eneq
We also compute that
\beq\label{rot-7}
{1\over{2\pi i}}\int_0^1\tau({dU(t)\over{dt}}U(t)^*)dt &=&
{1\over{2\pi
i}}\sum_{k=1}^{2m+1}\int_{k-1\over{2m+1}}^{k\over{2m+1}}
\tau({dU(t)\over{dt}}U(t)^*)dt\\
&=&{1\over{2\pi i}}(2\pi i\tau(h)+\sum_{k=1}^m 2\pi i\tau(b_k)
-\sum_{k=1}^m 2\pi i\tau(b_k))\\
&=& \tau(h).
\eneq
Thus $\tau(h)\in R_{\phi, \psi}(K_1(M_{\phi, \psi})).$

\end{proof}

The following is a generalization of the Exel trace formula for
bott element.

\begin{thm}\label{Exel}
There is $\dt>0$ satisfying the following:
 Let $A$ be a unital
separable simple \CA\, of tracial rank zero and let $u, v\in U(A)$
be two unitaries such that
\beq\label{Exel-1}
\|uv-vu\|<\dt.
\eneq
Then $\rm{bott}_1(u,v)$ is well defined and
\beq\label{Exel-2}
\rho_A(\rm{bott}_1(u,v))(\tau)={1\over{2\pi
i}}(\tau(\log(vuv^*u^*)))\tforal \tau\in T(A).
\eneq

\end{thm}

\begin{proof}
There is $\dt_1>0$ (one may choose $\dt_1=2$) such that for any
pair of unitaries $\text{bott}_1(u,v)$ is well defined (see
\ref{ddbot}). There is also $\dt_2>0$ satisfying the following: if
two pair of unitaries $u_1, v_1, u_2, v_2$ such that
\beq\label{exl-1}
\|u_1-u_2\|<\dt_2,\,\,\,\|v_1-v_2\|<\dt_2
\eneq
as well as
$$
\|[u_1,\, v_1]\|<\dt_1/2\andeqn \|[u_2,\,v_2]\|<\dt_1/2,
$$
then
\beq\label{exl-2}
\text{bott}_1(u_1,v_1)=\text{bott}_1(u_2, v_2).
\eneq
We may also assume that
$$
\|v_1u_1v_1^*u_1-1\|<1\,\,\,\text{whenever}\,\,\,
\|[u_1,\,v_1]\|<\dt_2.
$$

 Let
$$
F=\{z\in S^1: |z-1|<1+1/2\}.
$$
Let $\log: F\to (-\pi, \pi)$ be a smooth branch of logarithm.

 We choose
$\dt=\min\{\dt_1/2, \dt_2/2\}.$ Now fix a pair of unitaries $u,
v\in A$ with
\beq\label{exl-3}
\|[u,\,v]\|<\dt.
\eneq

For $\ep>0,$ there is $\dt_3>0$ such that
\beq\label{exl-3+1}
|\log(t)-\log (t')|<\ep
\eneq
provided that $|t-t'|<\dt_3$ and $t,t'\in F.$

Choose $\ep_1=\min\{\ep, \dt_3/4, \dt/2\}.$

Since $TR(A)=0,$ There is a finite dimensional \SCA\, $B$ of $A,$
a projection $p\in A$ with $1_B=p,$  unitaries $u', v'\in B,$ and
$u'', v''\in (1-p)A(1-p)$ such that
\beq\label{exl-4}
\|u'+u''-u\|<\ep_1,\,\,\, \|v'+v''-v\|<\ep_1\\\label{exl-5}
\andeqn \tau(1-p)<\ep_1\tforal \tau\in T(A).
\eneq
In particular,
\beq\label{exl-6}
\text{bott}_1(u'+u'',v'+v'')=\text{bott}_1(u,v).
\eneq
Therefore
\beq\label{exl-6+1}
\rho_A(\text{bott}_1(u'+u'',v'+v''))(\tau)=\rho_A(\text{bott}_1(u,v))(\tau)
\eneq
for all $\tau\in T(A).$

 Note that
\beq\label{exl-7}
\|[u',\, v']\|<\dt\andeqn \|[u'',\,v'']\|<\dt.
\eneq
We also have
\beq\label{exl-7+1}
\|v'u'(v')^*(u')^*-p\|<2\andeqn \|v''u''(v'')^*(u'')^*-(1-p)\|<2.
\eneq

Write $B=\oplus_{j=1}^r M_{l(j)}$ and
\beq\label{exl-8}
u'=\oplus_{j=1}^r u'(j)\andeqn v'=\oplus_{j=1}^rv'(j),
\eneq
where $u'(j), v'(j)\in M_{l(j)}$ are unitaries. By Exel's trace
formula (\cite{Ex}), one has
\beq\label{exl-9}
Tr_j(\text{bott}_1(u'(j), v'(j)))={1\over{2\pi i}}
Tr_j(\log(v'(j)u'(j)v(j)^*u(j)^*)),\,\,\,j=1,2,...,r,
\eneq
where $Tr_j$ is the standard trace on $M_{l(j)}.$ Suppose that
$\tau\in T(A)$ then there are $\lambda_j\ge 0$ such that
\beq\label{exl-10}
\sum_{j=1}^r {\lambda_j\over{l(j)}}=1\andeqn
\tau|_{B}=\sum_{j=1}^r{\lambda_j\over{l(j)}}Tr_j.
\eneq
It follows that
\beq\label{exl-11}
\tau(\text{bott}_1(u',v'))={1\over{2\pi
i}}\tau(\log(v'u'(v')^*(u')^*)).
\eneq

We also have
\beq\label{exl-12}
&&\hspace{-0.3in}\tau(\log((v'+v'')(u'+u'')(v'+v'')^*(u'+u'')^*)))\\
&&=
\tau(\log(v'u'(v')^*(u')^*))+\tau(\log(v''u''(v'')^*(u'')^*)).
\eneq
Note that
\beq\label{exl-13}
&&|{1\over{2\pi
i}}\tau(\log(v''u''(v'')^*(u'')^*))|<\tau(1-p)<\ep\\
&&\hspace{-0.2in}\andeqn
\tau(\text{bott}_1(u'', v''))<\ep
\eneq
for all $\tau\in T(A).$ It follows that
\beq\label{exl-14}
|\rho_A(\text{bott}(u,v))(\tau)-{1\over{2\pi
i}}\tau(\log((v'+v'')(u'+u'')(v'+v'')^*(u'+u'')^*)))| <2\ep
\eneq
for all $\tau\in T(A).$

By(\ref{exl-4}), we have
\beq\label{exl-14+1}
\|vuv^*u^*-(v'+v'')(u'+u'')(v'+v'')^*(u'+u'')^*\|<4\ep_1.
\eneq
Thus, by the choice of $\ep_1$ and  by (\ref{exl-4}),
\beq\label{exl-15}
|{1\over{2\pi
i}}[\tau(\log((v'+v'')(u'+u'')(v'+v'')^*(u'+u'')^*))-
\tau(\log(vuv^*u))]|<\ep
\eneq
for all $\tau\in T(A).$ Thus, by (\ref{exl-15}) and
(\ref{exl-14}),
\beq\label{exl-16}
|\rho_A(\text{bott}_1(u,v))(\tau)-{1\over{2\pi
i}}\tau(\log(vuv^*u))|<3\ep
\eneq
for all $\tau\in T(A)$ and for all $\ep.$ Let $\ep\to 0,$ we
obtain
\beq\label{exl-17}
\rho_A(\text{bott}_1(u,v)={1\over{2\pi i}}\tau(\log(vuv^*u))
\eneq
for all $\tau\in T(A).$

\end{proof}

\section{Asymptotic unitary equivalence}

\begin{lem}\label{NecL}
Let $A$ be a unital AH-algebra, and let  $B$ be  a  unital
separable simple \CA\, of tracial rank zero. Suppose that $\phi_1,
\phi_2: A\to B$ are two unital \hm s such that there is a
continuous path of unitaries $\{u(t): t\in [0,\infty)\}\subset B$
such that
\beq\label{NecL1}
\lim_{t\to\infty}{\rm ad}\, u(t)\circ \phi_1(a)=\phi_2(a)\rforal
a\in A.
\eneq
Then there is a continuous piecewise smooth path of unitaries
$\{w(t): t\in [0,\infty)\}\subset B$ such that
\beq\label{NecL2}
\lim_{t\to\infty}{\rm ad}\, v(t)\circ \phi_1(a)=\phi_2(a)\rforal
a\in A.
\eneq
\end{lem}

\begin{proof}
Define $w(t)=u(0)^*u(t).$ Then $w(0)=1$ and $w(t)\in U_0(B)$ for
each $t\in [0,\infty).$ Let $\{{\cal F}_n\}$ be an increasing
sequence of finite subsets of $A$ whose union is dense in $A.$
Without loss of generality, we may assume that, if $t\ge n,$
\beq\label{NecL3}
{\rm ad}\, u(0)w(t)\circ \phi_1\approx_{1/2^n}
\phi_2\,\,\,\text{on}\,\,\, {\cal F}_n.
\eneq
Let $\dt_n>0,$ ${\cal G}_n\subset A$ and ${\cal P}_n\subset
\underline{K}(A)$  be finite subsets  be required by \ref{LNHOMP}
(17.9 of \cite{Lnhomp}) corresponding to $\phi_2,$ $1/2^n$ and
${\cal F}_n,$ $n=1,2,....$ We may also assume that $
\text{Bott}(h, u')|_{{\cal P}_n} $ is  well-defined for any unital
\hm\, $h: A\to B$ and a unitary $u'\in B$ whenever $\|h(a),\,
u]\|<\dt_n$ for all $a\in {\cal G}_n.$

Let $\eta_n>0$ (in place of $\dt$) and let ${\cal G}_n'\subset A$
be a finite subset required by Lemma 9.3 of \cite{Lnhomp} for
$L=2\pi+1,$ ${\cal G}_n$ and $\dt_n/2$ (in place of $\ep$). We may
assume that $\eta_n<\dt_n/2$ and ${\cal G}_n'\supset {\cal G}_n.$

There are $\{s_n\}$ with $s_{n+1}>s_n\ge n$ such that
\beq\label{NecL4}
{\rm ad}\, u(0)w(t)\circ \phi_1\approx_{\eta_n/2}
\phi_2\,\,\,\text{on}\,\,\,{\cal G}_n'.
\eneq
It follows that
\beq\label{NecL5}
\|\phi_2(a),\, w(s(n))^*w(t)]\|<\dt_n\rforal a\in {\cal
G}_n'\andeqn t\ge s(n).
\eneq
Moreover, since $w(s(n))^*w(s(n))=1,$ we conclude that
\beq\label{NecL6}
\text{Bott}(\phi_2, w(s(n))^*w(s(n+1)))|_{{\cal P}_n}=0
\eneq
Combining \ref{LNHOMP} and the proof of Lemma 9.3 of
\cite{Lnhomp}, we obtain  a continuous and piecewise smooth path
$\{z_n(t): t\in [s(n), s(n+1)]\}$ such that
\beq\label{NecL7}
z_n(s(n))=1,\,\,\,z_n(s(n+1))=w(s(n))^*w(s(n+1))
\eneq
\beq\label{NecL7+}
\|[\phi_2(a),\, z_n(t)]\|<1/2^n\rforal a\in {\cal F}_n.
\eneq

It is standard to find a piecewise smooth path of unitaries
$\{v(t): t\in [0, s(1)]\}$ such that
\beq\label{NecL8}
v(0)=u(0)\andeqn v(s(1))=u(0)w(s(n)).
\eneq

Define $v(t)=u(0)w(s(n))z_n(t)$ if $t\in [s(n), s(n+1)),$
$n=1,2,....$ It is easy to check, by (\ref{NecL4}) and
(\ref{NecL7+}), that
$$
\lim_{t\to\infty}{\rm ad}\, v(t)\circ \phi_1(a)=\phi_2(a)
$$
for all $a\in A.$

\end{proof}

\begin{lem}\label{zeroV}
Let $A$  be a unital stably finite \CA\, and $u(t)\in C([0,1)),
A)$ be a unitary which is piecewise smooth. Suppose that
$w(t)=u(t)^*zu(t)$ for some unitary $z\in A.$
Then
\beq\label{zeroV1}
 \tau({dw(t)\over{dt}}w(t)^*)=0\tforal \tau\in T(A)\andeqn \tforal t\in [0,1).
\eneq
\end{lem}

\begin{proof}
Let $y=u(0),$ $z'=yzy^*$ and $v(t)=y^*u(t).$ Then $v(0)=1.$
Therefore by considering $v(t)^*z'v(t),$ without loss of
generality, we may assume that $u(0)=1.$

Let $0<d<1.$ We obtain $h_1(t), h_2(t),..., h_k(t)\in C([0,d], A)$
such that
$$
u(t)=\exp(ih_1(t))\cdot \exp(ih_2(t))\cdots \exp(ih_k(t))
$$
for $t\in [0,d].$ We claim that
\beq\label{zeroV-2}
\tau({dw(t)\over{dt}}w^*(t))=0\rforal t\in [0,d]
\eneq
and for all $\tau\in T(A).$ We prove this by induction on $k.$ For
$k=1,$ we may write $u(t)=\exp(ih(t))$ for $t\in [0,d].$ Thus
\beq\label{zeroV-3}
\tau({dw(t)\over{dt}}w(t)^*)&=&
\tau([-{dh(t)\over{dt}}u(t)^*zu(t)+u^*(t)z{dh(t)\over{dt}}u(t)]u(t)^*z^*u(t))\\
&=& \tau(-{dh(t)\over{dt}})+\tau(u^*(t)z{dh(t)\over{dt}}z^*u(t))\\
&=&-\tau({dh(t)\over{dt}})+\tau({dh(t)\over{dt}})=0
\eneq
for all $\tau\in T(A).$ Suppose that the claim holds for $k.$
Write $u(t)=\exp(ih(t))v(t),$ where $v(t)$ can be written as
product of $k$ exponentials. In particular,
\beq\label{zeroV-4}
\tau({d(v(t)^*zv(t))\over{dt}}(v(t)^*z^*v(t))dt)=0\tforal \tau\in
T(A)\andeqn t\in [0,d].
\eneq
Then
\beq\label{zeroV-4+}
\tau({dw(t)\over{dt}}w(t)^*)&=&\tau(-{dh(t)\over{dt}})+\tau(e^{-ih(t)}{d(v(t)^*zv(t))\over{dt}}(v(t)^*z^*v(t))e^{ih(t)})\\
&&+\tau(e^{ih(t)}v(t)^*zv(t){dh(t)\over{dt}}(v(t)^*z^*v(t)))\\
&=&-\tau({dh(t)\over{dt}})+\tau({d(v(t)^*zv(t))\over{dt}}(v(t)^*z^*v(t))+\tau({dh(t)\over{dt}})=0.
\eneq
Therefore
$$
\tau({dw(t)\over{dt}}w(t)^*)=0\rforal \tau\in T(A)\andeqn \rforal
t\in [0,d].
$$
It follows that
$$
\tau({dw(t)\over{dt}}w(t)^*)=0\tforal \tau\in T(A)\andeqn t\in
[0,1).
$$

\end{proof}

The following is a modification of (i) $\Rightarrow$ (ii) of 3.1
of \cite{KK2}.

\begin{thm}\label{NecT}
Let $A$ be a unital separable \CA\, satisfying the Universal
Coefficient Theorem and let $B$ be a unital separable  \CA.
Suppose that $\phi_1, \phi_2: A\to B$ are unital monomorphisms
such that
\beq\label{NecM1}
\lim_{to\infty}{\rm ad}\, u(t)\circ \phi_1(a)=\phi_2(a)\tforal
a\in A
\eneq
for some continuous and piecewise smooth path of unitaries
$\{u(t): t\in [0, \infty)\}\subset B.$  Then
\beq\label{NecM2}
&&[\phi_1]=[\phi_2],\,\,\, \tau\circ \phi_1=\tau\circ \phi_2\tforal
\tau\in T(A)\\\label{NecM2+}
&&\hspace{-0.2in}\andeqn {\tilde \eta}_{\phi_1, \phi_2}=0.
\eneq
\end{thm}

\begin{proof}

By changing parameter, we may assume that
$$
\lim_{t\to 1}{\rm ad}\, u(t)\circ \phi_1(a)=\phi_2(a)\rforal a\in
A
$$
and for some continuous path of unitaries $\{u(t): t\in
[0,1)\}\subset B.$

 Let
$$
M_{\phi_1, \phi_2}=\{f\in C([0,1], B): f(0)=\phi_1(a)\andeqn
f(1)=\phi_2(a)\,\,\,\text{for\,\,\,some}\,\,\, a\in A\}.
$$
Define $U(t)\in M_{2}(C([0,1],B)$ as follows.
\beq\label{NecM-3-}
U(t)=T_t\begin{pmatrix} 1 & 0\\
                                0 &
                                u(0)\end{pmatrix}T_t^{-1},
                                \eneq
for $t\in [0,1/2),$ where
$$
T_t=\begin{pmatrix}\cos ({\pi t}) & -\sin ({\pi t})\\
                     \sin({\pi t}) & \cos ({\pi t})
                     \end{pmatrix}
                     $$
   and define
\beq\label{NecM-4}
U(t)=\begin{pmatrix} u(2(t-1/2)) & 0\\
                          0 &1\end{pmatrix}.
\eneq

Define $\Phi: A\to M_2(M_{\phi_1, \phi_2})$ by
\beq\label{NecM-1}
\Phi(a)(t)&= &U(t)^*\begin{pmatrix}\phi_1(a) & 0\\
                                    0  & 0\end{pmatrix}U(t)\tforal
                                    t\in [0,1)
                                    \andeqn\\
\Phi(a)(1)&=&\begin{pmatrix}\phi_2(a) & 0\\
                              0 & 0\end{pmatrix}.
\eneq
$\Phi$ is a monomorphism. Moreover, $\pi_0\circ \Phi={\rm id}_A.$
This implies that the extension given by the mapping torus
$M_{\phi_1,\phi_2}$ is stably trivial. The map $\Phi$ also gives
the following splitting exact sequence:
$$
0\to \underline{K}(SB)\to
\underline{K}(M_{\phi_1,\phi_2})\hspace{0.1in}{\stackrel{\pi_0}{\rightleftarrows}}_{[\Phi]}
\hspace{0.05in}\underline{K}(A)\to 0.
$$

It follows that
$$
[\phi_1]=[\phi_2]\,\,\,\text{in}\,\,\, KK(A,B).
$$
It is also clear that
$$
\tau\circ\phi_1=\tau\circ\phi_2\tforal \tau\in T(B).
$$
To show that ${\tilde \eta}_{\phi_1, \phi_2}=0,$ we note that
$U(t)$ is piecewise smooth and continuous.

Let $z\in U(M_{2k}(A)).$ We write ${\tilde \Phi}$ for $\Phi\otimes
{\rm id}_{M_k}$
and
$$
{\tilde U(t)}={\rm diag}(\overbrace{U(t), U(t),...,U(t)}^k)\rforal
t\in [0,1].
$$
Then
\beq\label{NecM-3}
{\tilde \Phi}(z)(t)={\tilde U}(t)^*z{\tilde U}(t).
\eneq
It follows from \ref{zeroV} that
\beq\label{NecM-4+}
\int_0^1 \tau({d{\tilde \Phi}(z)\over{dt}} {\tilde
\Phi}(z)(t)^*)dt=0\tforal \tau\in T(B).
\eneq
This implies that
$$
R_{\phi_1,\phi_2}\circ\Phi_{*1}=0.
$$
Thus ${\tilde \eta}_{\phi_1,\phi_2}=0.$

\end{proof}

\begin{cor}\label{NecC}
Let $A$ be a unital AH-algebra and let $B$ be a unital separable
simple \CA\, with tracial rank zero. Suppose that $\phi_1, \phi_2:
A\to B$ are two unital monomorphisms which are asymptotically
unitarily equivalent. Then
\beq\label{Necc}
[\phi_1]=[\phi_2]\,\,\,\text{in}\,\,\,KK(A,B), {\tilde
\eta}_{\phi_1,\phi_2}=0\andeqn\\
\tau\circ \phi_1=\tau\circ \phi_2\tforal \tau\in T(B).
\eneq
\end{cor}

\begin{proof}
This follows from \ref{NecT} and \ref{NecL} immediately.

\end{proof}

\section{Almost commuting unitaries and bott element}

\begin{lem}\label{small}
Let $A$ be a unital separable simple \CA\, with real rank zero and
stable rank one and weakly unperforated $K_0(A).$ Suppose that
$x\in K_0(A)$ and there is a projection $p\in A$ such that
$$
2|\rho_A(x)(\tau)|<\tau(p)\tforal \tau\in T(A).
$$
Then there exist two projections $p_1,\, p_2\in pAp$ such that
$$
[p_1]-[p_2]=x.
$$
\end{lem}

\begin{proof}
It follows from \cite{BH} that $\rho_A(K_0(A))$ is dense in
$Aff(T(A)).$ Therefore there is a projection $q\in A$ such that
$$
|\rho_A(x)(\tau)|<\tau(q)<{\tau(p)\over{2}}\tforal \tau\in T(A).
$$
Let $y=[q]+x\in K_0(A).$ Then
$$
\rho_A(y)(\tau)>0\tforal \tau\in T(A).
$$
It follows that there is a projection $e\in M_K(A)$ for some
integer $K\ge 1$ such that $[e]=y.$  Then
$$
\tau(p)>\tau(e)\tforal \tau\in T(A).
$$
Then there is a projection $e'\in pAp$ such that $[e']=[e]$ in
$K_0(A).$ We may also assume that $q\le p.$ However,
$$
[e']-[q]=x.
$$

\end{proof}

\begin{lem}\label{hitK}
Let  $1>\ep>0.$ For any unital separable simple \CA\, $A$ with
real rank zero, stable rank one and weakly unperforated $K_0(A),$
any unitary $u\in U(A)$ with ${\rm sp}(u)=S^1,$ there exists
$\dt>0$ satisfying the following: for  any $x\in K_0(A)$ with
\beq\label{hit1}
|\rho_A(x)(\tau)|<\dt\rforal \tau\in T(A),
\eneq
there exists a unitary $v\in U_0(A)$ such that
\beq\label{hit2}
\|[u,\,v]\|<\ep\andeqn \rm{bott}_1(u,v)=x.
\eneq

\end{lem}

\begin{proof}
Fix $\ep>0.$ There exists an integer $n_0\ge 1$ such that
\beq\label{hit3}
|\omega_n-1|<\ep/2\rforal n\ge n_0,
\eneq
where $\omega_n=e^{2\pi i/n}.$

By 3.3 of \cite{Lncd}, for the given unitary $u,$ there is
$\dt_1>0$ and a finite subset ${\cal F}\subset C(S^1)$ satisfying
the following: for any other unitary $u_1\in U(A)$ with
\beq\label{hit4}
[u_1]=[u]\andeqn |\tau(f(u_1))-\tau(f(u))|<\dt_1
\eneq
for all $f\in {\cal F}$ and $\tau\in T(A),$ then there exists a
unitary $W\in U(A)$ such that
\beq\label{hit5}
\|W^*u_1W-u\|<\ep/4.
\eneq

Choose $\dt={\dt_1\over{8n_0}}.$ There exists a unital AF-algebra
$C$ such that
\beq\label{hit6}
(K_0(C), K_0(C)_+, [1_C])=(\rho_A(K_0(A)), \rho_A(K_0(A))_+, 1).
\eneq
Suppose that $x\in K_0(A)$ such that
\beq\label{hit7}
|\rho_A(x)(\tau)|<\dt\tforal \tau\in T(A).
\eneq
There are  projection $p_0^{(1)}\in A$ such that
\beq\label{hit8}
2\dt>\tau(p_0^{(1)})>2|\rho_A(x)(\tau)|\tforal \tau\in T(A).
\eneq
It follows from \ref{small} that there are two projections $p_1',
p_2'\in p_0^{(1)}Ap_0^{(1)}$ such that
\beq\label{hit9}
[p_1']-[p_2']=x.
\eneq
Put $y=\rho_A(x).$ There is a projection ${\bar p_0^{(1)}}\in C$
such that
\beq\label{hit9+}
[\bar p_0^{(1)}]=\rho_A([p_0^{(1)}]).
\eneq
In $C,$ there is a finite dimensional unital \SCA\, $B\subset
{\bar p_0^{(1)}}C{\bar p_0^{(1)}}$ and two projections $q_1,
q_2\in B$ such that
\beq\label{hit10}
[j(q_1)]-[j(q_2)]=y
\eneq
where $j: B\to C$ is the (unital) embedding. Thus we obtain a
unital monomorphism $\psi: B\to p_0^{(1)}Ap_0^{(1)}$ such that
\beq\label{hit11}
\rho_A([\psi(q_1)]-[\psi(q_2)])=y.
\eneq
Write $B=M_{n_1}\oplus M_{n_2}\oplus\cdots \oplus M_{n_K}.$ Since
$C$ is a simple AF-algebra, we may choose (larger) $B$ so that
$n_k\ge n_0,$ $k=1,2,...,K.$ We may write
\beq\label{hit12}
\rho_A(\psi_*(k_1, k_2,...,k_K))=y,
\eneq
where $|k_i|\le n_i$ are integers, $i=1,2,...,K.$

Let $E_k=\psi(1_{M_{n_k}}),$ where we view $1_{M_{n_k}}$ as a
projection in $B,$ $k=1,2,...,K,$ Let $\{e_{i,j,k}\}$ be a system
of matrix unit for $M_{n_k}.$

Define
\beq\label{hit13}
S_{n_k}=\begin{pmatrix} 0 & 0 & \ldots  & 0 & 1\cr
                 1 & 0 & \ldots & 0 & 0\cr
                 0 & 1 & \ldots & 0 & 0\cr
                 \vdots &\vdots& \ddots & \vdots& \vdots\cr
                 0 & 0 & \ldots & 1 & 0
                 \end{pmatrix}\andeqn
                 W_{n_k}=\begin{pmatrix} \omega_{n_k} &0 &\ldots &0 \cr
                                               0 & \omega_{n_k}^2 &\ldots &0\cr
                                               \vdots & \vdots &\ddots & \vdots  \cr
                                                  0 & 0 &\ldots  &
                                                  \omega_{n_k}^{n_k}
                                                 \end{pmatrix} . \eneq
We compute that
\beq\label{hit13+}
{1\over{2\pi i}}{\rm Tr}(\log(W_{n_k}S_{n_k}W_{n_k}^*S_{n_k}^*))=1.
\eneq
Let $n_0'=\max\{|k_i|:1\le i\le K\}.$  By (\ref{hit8}), there are $n_0'-1$ mutually
orthogonal and mutually equivalent projections
$p_0^{(2)},p_0^{(3)},...,p_0^{(n_0')}\in
(1-p_0^{(1)})A(1-p_0^{(1)})$ for which $p_0^{(i)}$ is equivalent
to $p_0^{(1)}.$
Thus there is a unitary $X_i\in A$
such that
$$
X_i^*p_0^{(1)}X_i=p_0^{(i)},\,\,\,i=1,2,...,n_0'.
$$
Note we assume that $X_1=1_A.$

We write $S_{n_k}=\sum_{i,j}^{n_k} a_{i,j}e_{i,j,k},$ where
$a_{i,j}$ are zero's or $1.$
Define $T_s=\sum_{m=1}^{k_s}X_m^*S_{n_k}X_m$ and
$e_s=\sum_{m=1}^{k_s}X_m^*p_0^{(1)}X_m,$ $s=1,2,...,K.$ Define
$p=\sum_{s=1}^Ke_s.$
Define $v_{1,s}=\sum_{m=1}^{k_s}X_i^*(\sum_{i=1}^{n_s}
\omega_{n_s}^ie_{i,i,s})X_i,$ $s=1,2,...,K.$ Define
$u_1=\sum_{s=1}^K T_s$ and  define $v_1$ to
the sum of $v_{1,s}$ or $v_{1,s}^*$ depending $k_j$ is
non-negative or negative. Then one computes easily
\beq\label{hit14}
\|u_1v_1-v_1u_1\|\le \max_k\{|1-\omega_{n_k} |\}\le
|1-\omega_{n_0}|.
\eneq
Moreover, by (\ref{hit12}), the Exel trace formula and
(\ref{hit13+}),
\beq\label{hit15}
\rho_A(\text{bott}_1(u_1, v_1))=y.
\eneq
Let $p_1\in (1-p)A(1-p)$ such that $1-p-p_1\not=0$ and
\beq\label{hit16}
\tau(1-p-p_1)<\dt\tforal \tau\in T(A).
\eneq
Let $p_2,\,p_2\in (1-p-p_1)A(1-p-p_1)$ be two mutually orthogonal
non-zero projections such that $p_2+p_3=1-p-p_1.$

It follows from 4.3 of \cite{Lnhomp} that there is a unitary
$u_2\in p_1Ap_1$ with finitely many points in the spectrum of
$u_2$ such that
\beq\label{hit17}
|\tau_1(f(u_2))-\tau(f(u))|<\dt/2\rforal \tau\in T(A)
\eneq
for all $f\in {\cal F},$ where
$\tau_1(b)={\tau(b)\over{\tau(p_1)}}$ for all $b\in p_1Ap_1.$ Let
$x_1=x-\text{bott}_1(u_1, v_1).$ Then
\beq\label{hit18}
x_1\in {\rm ker}\rho_A.
\eneq
It follows from \cite{LnKT} that there exists a unital \hm\,
(necessarily injective) $\phi_1: C(S^1\times S^1)\to p_2Ap_2$ such
that
\beq\label{hit19}
\text{bott}_1(\phi_1(z\otimes 1),\phi_1(1\otimes z))=x_1\andeqn
(\phi_1)_{*1}=0,
\eneq
where $z$ is the identity function on the unit circle.

There is a unitary $u_3\in p_3Ap_3$ such that $[u_3]=[u].$

Now define $U=u_1\oplus u_2\oplus\phi_1(z\otimes 1)\oplus u_3$ and
$v'=v_1\oplus p_21\oplus \phi_1(1\otimes z)\oplus  p_3.$ Since
$u_1\in \psi(B),$ $u_2$ has finite spectrum and $(\phi_1)_{*1}=0,$
$$
[U]=[u].
$$
One also easily verify that
\beq\label{hit20}
\|[U,\,v']\|\le |1-\omega_n|<\ep/2\andeqn \text{bott}_1(U, v')=x.
\eneq
Moreover,
\beq\label{hit21}
|\tau(f(U))-\tau(f(u))|<\dt\tforal \tau\in T(A)
\eneq
for all $f\in {\cal F}.$ By 3.3 of \cite{Lncd}, there is a unitary
$W\in U(A)$ such that
\beq\label{hit22}
\|W^*UW-u\|<\ep/4.
\eneq
Now define $v=W^*v'W.$ It follows that
\beq\label{hit23}
\|uv-vu\| &\le &\|[U,v']\|+\|W^*UW-u\|\\
&<& |1-\omega_{n_0}|+\ep/4<\ep
\eneq
Moreover,
\beq\label{hit24}
\text{bott}_1(u,v)=x.
\eneq

\end{proof}

\begin{rem}\label{rehitk}
{\rm It is the application of 3.3 of \cite{Lncd} that requires
$\dt$ depending on $u.$ However, if we apply Theorem 4.6 of
\cite{Lncd}, then we can choose $\dt$ independent of $u$ as well
as $A$ but $\dt$ is required to be dependent on a given measure
distribution on the unit circle.

}

\end{rem}

The following is certainly known and follows from 2.6.11 of
\cite{Lnbk} immediately.

\begin{lem}\label{spnotfull}
Let $X$ be a compact proper subset of the unit circle. Then there
exists $\dt>0$ satisfying the following. If $A$ is a unital \CA,
$u$ and $v$ are two unitaries with ${\rm sp}(u)\subset X$ and
\beq\label{spn-1}
\|[u,\, v]\|<\dt,
\eneq
then \beq
\rm{bott}_1(u,v)=0.
\eneq
\end{lem}

\begin{proof}
Let $\ep=1/2.$ Let $d>0$ so that $S^1\setminus X$ contains an arc
of length $d.$ Let $\dt>0$ be as in 2.6.11 of \cite{Lnbk}. Then,
by 2.6.11 of \cite{Lnbk}, there is a selfadjoint element $h\in A$
such that $\exp(ih)=u$ and
\beq\label{spn-2}
\|[\exp(ith),\, v]\|<\ep\rforal t\in [0,1].
\eneq
Let $u(t)=\exp(ith).$ Then $u(0)=1$ and $u(1)=u.$ It follows from
(\ref{spn-2}) that
$$
\text{bott}_1(u,v)=0.
$$

\end{proof}

\section{Finite CW complex}

\begin{df}\label{Dkind-0}
{\rm  Let $f: S^1\to S^1$ be a degree $k$ map ($k>1$), i.e., a
continuous map with the winding number $k.$ Following 4.2 of
\cite{EG1}, denote by $T_{II,k}=D^2\cup_f S^1,$ the connected
finite CW complex obtained by attaching a $2$-cell $D^2$ to $S^1$
via the map $f.$ Then
\beq\label{dfkind1}
K_0(C(T_{II,k}))=\Z\oplus k\Z\andeqn K_1(T_{II,k})=\{0\}.
\eneq

Let $g: S^2\to S^2$ be a degree $k$ map ($k>1$). Let
$T_{III,k}=D^3\cup_g S^2$ be the connected finite CW complex
obtained by attaching a $3$-cell $D^3$ to $S^2$ via the map $g.$
Then
\beq\label{dfkind2}
K_0(C(T_{III,k}))=\Z\andeqn K_1(C(T_{III,k}))=\Z/k\Z
\eneq
(see 4.2 of \cite{EG1}). }

\end{df}

\begin{df}\label{Dfkind1}
{\rm A unital AH-algebra is said to be of first kind if it is
inductive limit of \CA s with the form $M_l(C(S^1))$ and
$PM_l(C(T_{II,k})P.$

$A$ is said to be of second kind if  it is inductive limit of
finite direct sums of \CA s with the form $PM_l(C(T_{II,k})P$ and
$PM_l(C(T_{III,k}))P.$ }

\end{df}

The following statements follows immediately from \cite{EG1} and
\cite{Lnduke}

\begin{prop}\label{eg1}
Let $A\in {\cal N}$ be a unital separable simple \CA\, of tracial
rank zero.

{\rm (1)} Suppose that $K_1(A)$ is torsion free ($i=0,1$). Then
$A$ is a unital  AH-algebra of the first kind.

{\rm (2)} Suppose that $K_1(A)$ is torsion. Then $A$ is a unital
AH-algebra of the second kind.
\end{prop}

\begin{lem}\label{ATemb1}
 Let $A$ be a unital separable simple \CA\, with tracial
rank zero.

{\rm (1)} Suppose that there is a \hm\,  $\gamma: \Z^k\to K_1(A).$
Then there exists a unital simple AH-algebra $A_0$ of the first
kind with tracial rank zero and with $K_1(A_0)=\Z^k,$ and there
exists a unital monomorphism $\phi_0: A_0\to A$ such that
$\phi_{*0}$ gives an order isomorphism and $\phi_{*1}=\gamma.$

{\rm (2)} There is a unital simple AH-algebra $A_1$ of the second
kind with tracial rank zero with $K_1(A_1)={\rm Tor}(K_1(A))$ and
a unital monomorphism $\phi_1: A_1\to A$ such that $(\phi_1)_{*0}$
is an order isomorphism and $(\phi_1)_{*1}={\rm
id}|_{Tor(K_1(A))}.$

\end{lem}

\begin{proof}
By \cite{EG1}, there is a unital simple AH-algebra $A_0$ of real
rank zero and slow dimension growth such that
\beq\label{atemb1-1}
K_1(A_0)&=&\Z^k\andeqn\\ (K_0(A_0),
 K_0(A_0)_+, [1_{A_0}]) &=&(K_0(A), K_0(A)_+, [1_A]).
 \eneq
 Note that $A_0$ has tracial rank zero (see \cite{LnTAF}).
It also follows from 4.6 of \cite{LnKT} that there is a unital
monomorphism $\phi_0: A_0\to A$ such that $(\phi_0)_{*0}$ is the
order isomorphism and $(\phi_0)_{*1}=\gamma.$ This proves (1).

For (2), by \cite{EG1}, there is a unital simple AH-algebra $A_1$
with real rank zero and slow dimension growth such that
\beq\label{atemb1-2}
K_1(A_1)&=&{\rm Tor}(K_1(A))\andeqn\\
K_0(A_1),
 K_0(A_1)_+, [1_{A_1}]) &=&(K_0(A), K_0(A)_+, [1_A]).
 \eneq
As the above, one obtains a unital monomorphism $\phi_1: A_1\to A$
such that $(\phi_1)_{*0}$ gives the order isomorphism and
$(\phi_1)_{1*}={\rm id}_{{\rm Tor}(K_1(A))}.$

\end{proof}

The following theorem follows immediately from a version of the
statement in \cite{hlx}

\begin{lem}\label{emb1}
Let $X$ be a compact metric space  with $K_1(C(X))={\mathbb
Z}^k\oplus G_0,$ where $G_0$ is a finite group.

Let $A$ be unital separable simple \CA\,  with tracial rank zero
and let $\gamma: K_1(C(X))\to K_1(A).$ Suppose that
$\lambda:C(X)_{s.a.}\to Aff(T(A))$ is a positive linear map. Then
there exists a unital simple AH-algebra $A_0$ of tracial rank zero
with $K_1(A_0)=\Z^k$ and a unital monomorphism $\phi_0: A_0\to A$
such that $\phi_{*0}$ gives an order isomorphism and
$\phi_{*1}=\gamma|_{\Z^k}.$  Moreover, for any $\dt>0$ and any
finite subset ${\cal G}\subset C(X)_{s.a},$ there exists a unital
monomorphism $h: C(X)\to A_0$ such that
\beq\label{emb1-1}
\phi_0\circ h_{*1}|_{\Z^k}=\gamma|_{\Z^k},
\,\,\,h_{*1}|_{G_0}=0\,\,\,and\,\,\, |\tau\circ
h(a)-\lambda(a)(\tau)|<\dt\tforal a\in {\cal G}.
\eneq

\end{lem}

\begin{proof}
Choose two non-zero mutually orthogonal projections $p_1, p_2\in
A_0$ such that
\beq\label{emb1-2}
\tau(p_j)<\dt/4\tforal \tau\in T(A_0),\,\,\,j=1,2.
\eneq
Note that $K_i(pA_0p)=K_i(A_0),$ $i=0,1.$ It follows from 4.7 of
\cite{LnKT} that there is a unital \hm\, $h_0: C(X)\to p_1A_0p_1$
such that
\beq\label{emb1-3}
(h_0)_{1*}|_{\Z^k}={\rm id}_{\Z^k}.
\eneq
There is a non-zero projection $p_1\in A_0.$ It follows from 3.6
of \cite{hlx} that there is a unital \hm\, $h_{00}: C(X)\to
(1-p_1-p_2)A_0(1-p_1-p_2)$ with finite dimensional range such that
\beq\label{emb1-4}
|\tau\circ h_{00}(f)-\lambda(f)(\tau)|<\dt/4\tforal a\in {\cal G}.
\eneq
It is well known that there exists a unital monomorphism $h_1:
C(X)\to p_2A_0p_2$ which factors through $C([0,1]).$ Now define
$h: C(X)\to A_0$ by $h(f)=h_0(f)+h_{00}(f)+h_1(f)$ for all $f\in
C(X).$

\end{proof}

\begin{lem}\label{app}
Let $X$ be a connected finite CW complex with
$K_1(C(X))=\Z^k\oplus G_0,$ where $G_0$ is a finite group, and let
$A$ be a unital separable simple \CA\, of tracial rank zero. Let
$C=PM_r(C(X))P,$ where $l>0$ is an integer, $P\in M_l(C(X))$ is a
projection.

 Suppose that $h: C\to A$ is a unital monomorphism.
Then, for any $\ep>0$ and any finite subset ${\cal F}\subset
C(X),$ there is a projection $p\in A,$ a unital monomorphism
$\phi_0: A_0\to pAp,$ where $A_0$ is a unital simple AH-algebra of
the first kind with tracial rank zero and with $K_1(A_0)=\Z^k$
such that $(\phi_0)_{*1}=h_{*1}|_{\Z^k},$ and there is a unital
simple AH-algebra $A_1$ of second kind with tracial rank zero and
with $K_1(A_1)={\rm Tor}(K_1(A)),$ a unital monomorphism $\phi_1:
A_1\to (1-p)A(1-p)$ with $(\phi_1)_{*1}= {\rm id}_{{\rm
Tor}(K_1(A))}$ and $(\phi_1)_{*0}$ is an order isomorphism,
 and a unital
monomorphism $h_1: C\to A_0$ with $(h_1)_{*1}|_{{\mathbb
Z}^k}={\rm id}_{\Z^k}$ and a unital \hm\, $h_2: C\to A_0$ with
$(h_2)_{*1}|_{\Z^k}=0$ and $(h_2)_{*1}|_{G_0}=h_{*1}|_{G_0}$ such
that
\beq\label{app1}
\|\phi_0\circ h_1(f)\oplus \phi_1\circ h_2(f)-h(f)\|<\ep\tforal
f\in {\cal F}.
\eneq
\end{lem}

\begin{proof}
We first assume that $C=C(X).$ The general case will be dealt with
at the end of this proof.

Write $K_0(C)=\Z\oplus G_{00},$ where $G_{00}={\rm ker}\rho_C.$

 Fix $\ep>0$ and a finite subset
${\cal F}\subset C(X).$ Let $\dt>0,$ $\gamma>0,$ a finite subset
${\cal G}\subset C(X)$ and a finite subset ${\cal P}\subset
\underline{K}(C)$ required by 3.3 of\cite{Lncd} for $h,$ $\ep$ and
${\cal F}.$ Without loss of generality (by taking even smaller
$\dt$ for example), we may assume that ${\cal G}\subset
C(X)_{s.a.}.$

Let $\lambda: C(X)_{s.a}\to Aff(T(A))$ be induced by $h.$  Choose
a non-zero projection $e\in A$ such that $\tau(e)<\gamma/2$ for
all $\tau\in T(A).$ Let $A_0$ and $\phi_0$ be as in \ref{emb1} but
we replace $A$ by $(1-e)A(1-e),$ and $\dt$ by $\gamma/2.$ We also
use $h_1$ for $h$ given in \ref{emb1}. Define
\beq\label{app-2-1}
\kappa=[h]-[\phi_0\circ h_1]\in KL(C,
A)=Hom_{\Lambda}(\underline{K}(C),\underline{K}(A)) .
\eneq
 In particular, by identify $\Z^k$ with the free part of
$K_1(C),$
\beq\label{app-2k}
\kappa|_{\Z^k}=0\andeqn \kappa|_{G_0}=h_{*1}|_{G_0}.
\eneq
Moreover, $\kappa([1_c])=[e]$ and $\rho_A(\kappa(G_{00}))=0.$

By (2) of \ref{ATemb1}, there is a unital simple AH-algebra $A_1$
of the second kind with tracial rank zero and with $K_1(A_1)={\rm
Tor}(K_1(A)),$ and a unital monomorphism $\phi_1: A_1\to eAe$ such
that $(\phi_1)_{*0}$ is an order isomorphism and
\beq\label{app-2}
(\phi_1)_{*1}={\rm id}|_{{\rm Tor}(K_1(A))}.
\eneq
From
$$
{\small \put(-160,0){$K_0(A_1)$} \put(0,0){$K_0(A_1,{\Z }/k\Z )$}
\put(180,0){${\rm Tor}(K_1(A))$} \put(-85,-40){$K_0(A)$}
\put(0,-40){$K_0(A,{\Z}/k{\Z})$} \put(105,-40){$K_1(A)$} \put(-85,
-70){$K_0(A)$} \put(0,-70){$K_1(A, {\Z}/k{\Z})$}
\put(105,-70){$K_1(A)$} \put(-160,-110){$K_0(A_1)$}
\put(0,-110){$K_1(A_1,{\Z}/k{\Z})$} \put(180,-110){${\rm
Tor}(K_1(A))$\,,} \put(-120, 2){\vector(1,0){95}}
\put(70,1){\vector(1,0){95}} \put(-123,-3){\vector(1,-1){30}}
\put(30,-3){\vector(0,-1){25}} \put(180,-2){\vector(-1,-1){30}}
\put(-45,-38){\vector(1,0){35}} \put(70,-38){\vector(1,0){25}}
\put(-147, -90){\vector(0,1){85}} \put(-75,-60){\vector(0,1){15}}
\put(115, -45){\vector(0,-1){15}} \put(190,-7){\vector(0,-1){85}}
\put(-7,-68){\vector(-1,0){35}} \put(95,-68){\vector(-1,0){25}}
\put(-123,-102){\vector(1,1){30}} \put(175,
-105){\vector(-1,1){30}} \put(30,-104){\vector(0,1){30}} \put(-5,
-108){\vector(-1,0){100}} \put(170,-108){\vector(-1,0){95}}
\put(-112,-12){$[\phi_1]$} \put(13, -15){$[\phi_1] $} \put(145,
-12){$[\phi_1]$} \put(-135,-92){$[\phi_1]$} \put(13,
-88){$[\phi_1]$} \put(160, -88){$[\phi_1]$} }
$$
one calculates that all maps $[\phi_1]$ are injective and
$[\phi_1]$ from $K_0(A_1, \Z/k\Z)$ to $K_0(A, \Z/k\Z)$ is an
isomorphism.  Since $(h_1)_{*1}$ maps $G_0$ into ${\rm
Tor}(K_1(A)),$ by (\ref{app-2k}), then there is $x\in
Hom_{\Lambda}(\underline{K}(C), \underline{K}(A_1))$ such that
\beq\label{app-3}
x\times [\phi_1]=\kappa.
\eneq
Note that $x([1_C])=[1_{A_1}]$ and $x(G_{00})\in {\rm
ker}\rho_{A_1}.$ It follows from 4.7 of \cite{LnKT} that there
exists a unital \hm\, $h_2: C\to A_1$ such that
\beq\label{app-4}
[h_2]=x.
\eneq

Define $h_3: C(X)\to A$ by $h_3=\phi_0\circ h_1\oplus \phi_1\circ
h_2.$ Then,
\beq\label{app-5}
[h_3]=[h]\,\,\,\text{in}\,\,\,KL(C, A)
\eneq
We also have
\beq\label{app-6}
|\tau\circ h_3(f)-\tau\circ h(f)|<\gamma\rforal f\in {\cal G}.
\eneq
It follows from 3.3 of \cite{Lncd} that there exists a unitary
$u\in U(A)$ such that
\beq\label{app-7}
{\rm ad}\, u\circ h_3\approx_{\ep} h_2\,\,\,\text{on}\,\,\, {\cal
F}.
\eneq
This proves the lemma in the case that $C=C(X).$

Now consider the case that $C=M_l(C(X)).$ Let $e_{1,1}\in
M_l(C(X))$ be a constant rank one projection. Let
$e_0=h(e_{1,1}).$ By replacing $A$ by $e_0Ae_0$ and $h$ by
$h|_{e_{1,1}Ce_{1,1}},$ one easily sees that this case reduces to
the case that $C=C(X).$

Now consider the general case that $C=PM_l(C(X))P.$ There is $r\ge
1$ such that there is projection $q\in M_r(C)$ such that
$(q+P)M_r(C)(q+P)\cong M_{l_1}(C(X)).$ We then first extend $h$ to
$h\otimes {\rm id}_{M_r}: M_r(C)\to M_r(A).$ We then obtain the
desired \hm s $h_1', h_2': M_r(C)\to M_r(A)$ as well as $\phi_0'$
and $\phi_1'.$ We then restrict $h_1'$ and $h_2'$ on $PM_r(C)P=C.$
It is easy to see that the general case can be reduced to the case
that $C=M_l(C(X)).$

\end{proof}

\begin{lem}\label{KTemb}
Let $X$ be a connected finite CW complex, $C=PM_l(C(X))P,$ where
$P\in M_l(C(X))$ is a projection, and let $A$ be a unital
separable simple \CA\, with tracial rank zero. Suppose that $h: C
\to A$ is a unital monomorphism. Suppose that $\kappa\in
Hom_{\Lambda}(\underline{K}(C\otimes C(S^1)), \underline{K}(A))$
such that
\beq\label{KTemb-1}
\kappa|_{\underline{K}(C)}=[h]\andeqn
\rho_A(\kappa({{\boldsymbol{\bt}}(K_1(C))}))=0.
\eneq
Then, for any $\ep>0$ and any finite subset ${\cal F}\subset C,$
there exists a unitary $u\in U(A)$ such that
\beq\label{KTemb-2}
\|[h(a), u]\|<\ep\tforal a\in {\cal F}\andeqn \rm{Bott}(h,
u)=\kappa|_{{\boldsymbol{\bt}}(\underline{K}(C))}.
\eneq

\end{lem}

\begin{proof}
Suppose that $d={\rm dim} X.$  We will prove this for the case
that $C=M_r(C(X)).$ The general case follows from the fact that
there is a projection $Q\in M_{d+1}(C)$ such that
$(Q+P)M_{d+1}(C)(Q+P) \cong M_{l'}(C(X)).$

Define $H: C\otimes C(S^1)\to A$ by $H(f\otimes g)=h(f) g(1)$ for
all $f\in C$ and $g\in C(S^1).$

Since $K_i(C)$ is finitely generated, from \ref{ddbot}, as in 2.10
of \cite{Lnhomp}, there is $\dt_0>0,$ a finite subset ${\cal
G}_0\subset C$ and a finite subset ${\cal P}_0\subset
\underline{K}(C)$ satisfy the following: for any unital \hm\, $h':
C\to A$ and any unitary $u'\in U(A)$ with
$$
\|[h'(a), u']\|<\dt_0,
$$
$\text{Bott}(h',u')|_{{\cal P}_0}$ is well defined and
$\text{Bott}(h',u')|_{{\cal P}_0}$ defines $\text{Bott}(h',u').$

 Let $\ep>0$ and let ${\cal F}\subset C$ be a finite
subset of the unit ball. Denote $S=\{1_{C(S^1)}, z\},$ where $z\in
C(S^1)$ is the identity function on the unit circle.

Let $\dt_1>0$  and ${\cal F}_1\subset C$ be finite subset such
that
$$
\text{Bott}(h',u')|_{{\cal P}_0}=\text{Bott}(h'',u'')|_{{\cal
P}_0}
$$
for any two pairs of $h',$ $u'$ and $h'',$ $u''$ for which
$$
\|[h'(a), \,u']\|<\dt_1\andeqn \|[h''(a), \,u'']\|<\dt_1
$$
for all $a\in {\cal F}_1,$ provided that
\beq\label{KTemb-2-1}
h'\approx_{\dt_1} h''\,\,\,\text{on}\,\,\, {\cal F}_1\andeqn
\|u'-u''\|<\dt_1
\eneq

Let $\ep_1=\min\{\ep/2, \dt_0/2, \dt_1/2\}$ and ${\cal F}_3={\cal
F}\cup{\cal G}_0\cup {\cal F}_2.$

 Let $\dt>0, \gamma>0$ and ${\cal G}\subset C$ be a finite
subset which are required by 3.3 of \cite{Lncd} for $\ep_1$ and
${\cal F}_3$ and $h$ above. We may assume that $\dt<\ep_1$ and
${\cal G}\supset {\cal F}_3.$ Moreover, without loss of
generality, we may assume that ${\cal G}$ is a subset of the unit
ball of $C.$

Define $B=C\otimes C(S^1).$ We may write $K_0(B)=\Z\oplus G_0,$
where $G_0={\rm ker}\rho_B.$ Suppose that $k\ge 1$ is an integer
such that
\beq\label{KTemb-3-1}
kx=0\tforal x\in {\rm Tor}(K_i(B)),\,\,\, i=0,1.
\eneq

Put $m=k!.$ Then, by 2.11 of \cite{DL},
\beq\label{KTemb-3}
Hom_{\Lambda}(F_m\underline{K}(B),F_m\underline{K}(B))\cong
Hom_{\Lambda}(\underline{K}(B), \underline{K}(B)).
\eneq

 Since $A$ is a simple \CA\, with tracial rank zero,
it is easy to find two  mutually orthogonal projections $e, e'\in
A$ such that
\beq\label{KTemb-4}
\tau(e)<\gamma/3\tforal \tau\in T(A)\andeqn
[e]+k![e']=[1_A]\,\,\,\text{in}\,\,\, K_0(A).
\eneq

Therefore
\beq\label{KTemb-5}
{\overline{[e]}}=\overline{[1_A]}\,\,\, \text{in}\,\,\,
K_0(A,\Z/n\Z),
\eneq
where ${\overline{[e]}}$ and ${\overline{[1_A]}}$ are images of
$[e]$ and $[1_A]$ in $K_0(A, \Z/n\Z),$ respectively,
$n=2,3,...,m.$

 Let $\kappa'\in Hom_{\Lambda}(F_m\underline{K}(B), F_m\underline{K}(A))$ be defined by
\beq\label{KTemb-7}
\kappa'([1_B])&=&[e],\,\,\,\kappa'|_{G_0}=\kappa|_{G_0}\\\label{KTemb-7++}
\kappa'|_{K_1(B)}&=&\kappa|_{K_1(B)}\andeqn
\kappa'|_{K_i(B,\Z/n\Z)},\,\,\,n=2,3,...,m.
\eneq
Note that, by (\ref{KTemb-5}), so defined $\kappa'$ is indeed in
$Hom_{\Lambda}(F_m\underline{K}(B), F_m\underline{K}(A)).$
Furthermore,
\beq\label{KTemb-7+1}
\rho_A(\kappa(x))=0\tforal x\in G_0.
\eneq
 By (\ref{KTemb-3}),
\beq\label{KTemb-8}
\kappa'\in KK(B,A)_+.
\eneq
(see \cite{LnKT}) It follows from 4.7 of \cite{LnKT} that there
exists a unital \hm\, $\phi_1: B\to eAe$ such that
\beq\label{KTemb-9}
[\phi_1]=\kappa'.
\eneq
Let $e_1\in (1-e)A(1-e)$ be a projection such that
$\tau(e_1)<\gamma/3$ for all $\tau\in T(A).$

 It follows from 2.25 (1) of \cite{Lnemb2} that there exists a unit
 monomorphism $\phi_2: B\to e_1Ae_1$ which factors through
 $M_r(C([0,1]).$
 Let $\xi\in X$ be a point and let $Y=X\setminus \{\xi\}.$
 Then
 \beq\label{KTemb-9+}
[\phi_2|_{M_r(C_0(Y))}]=0\,\,\,\text{in}\,\,\,KL(C_0(Y), A).
\eneq

It follows from 3.6 of \cite{hlx} that there exists a unital \hm\,\\
$\phi_3: C\to (1-e-e_1)A(1-e-e_1)$ with finite dimensional range
such that
\beq\label{KTemb-10}
|\tau\circ \phi_3(a)-\tau\circ \phi_3(a)|<\gamma/3
\eneq
for all  $a\in {\cal G}\otimes S$ and for all $\tau\in T(A).$
Since $\phi_3$ has finite dimensional range,
\beq\label{KTemb-10+}
[\phi_3|_{M_r(C_0(Y))}]=0\,\,\,\text{in}\,\,\,KL(C_0(Y), A).
\eneq

Now define $\phi: B\to A$ by
\beq\label{KTemb-11}
\phi(f)=\phi_1(f)+\phi_2(f)+\phi_3(f)\tforal f\in B.
\eneq
By (\ref{KTemb-7}),(\ref{KTemb-7+1}), (\ref{KTemb-9+}) and
(\ref{KTemb-10+}), one computes that
\beq\label{KTemb-12}
[\phi]=\kappa.
\eneq
Define $\psi: C\to A$ by $\psi(f)=\phi(f\otimes 1_{C(S^1)})$ for
$f\in C.$ Then
\beq\label{KTemb-13}
[h]=[\psi]\andeqn |\tau(h(f))-\tau(\psi(f))|<\gamma
\eneq
for all $f\in {\cal G}$ and for all $\tau\in T(A).$ It follows
from 3.3 of \cite{Lncd} that there exits a unitary $w\in A$ such
that
\beq\label{KTemb-14}
{\rm }\, w\circ\psi(f)\approx_{\ep_1} h(f)\rforal f\in {\cal F}_3.
\eneq
Put
$$
u={\rm ad}\, w(\phi(1\otimes z)).
$$
Then,
\beq\label{KTemb-15}
\|[u,\, h(a)]\|<\ep_1 \rforal a\in{\cal F}_3.
\eneq
By (\ref{KTemb-2-1}),
\beq\label{KTemb-16}
\text{Bott}(h,u)=\text{Bott}(\psi,
u)=\kappa|_{{\boldsymbol{\bt}}(\underline{K}(C))}.
\eneq

\end{proof}

\begin{lem}\label{ebot1}
Let $C=\oplus_{i=1}^rC_i,$ where $C_i=P_iM_{l_i}(C(X_i))P_i,$
$P_i\in M_{l_i}(C(X_i))$ is a projection and
 $X_i=S^1,$  $X_i=T_{II,m_i}$ or
$X_i=T_{III,m_i}.$ Let $A$ be a unital separable simple \CA\, of
tracial rank zero and $\psi: C\to A$ is a unital monomorphism.

Let $x_i$ be a generator of $K_1(C_i)$ (if
$K_1(C(X_i))\not=\{0\}$).
  For any $\ep>0$ and any
finite subset ${\cal F}\subset C,$ there is $\dt>0$ satisfying the
following: For any $y_i\in K_0(A),$ if there is a \hm\, $\kappa_i:
K_1(C_i)\to K_0(A)$ with $\kappa_i(x_i)=y_i$ and
\beq\label{ebot1-1}
|\rho_A(y_i)(\tau)|<\dt\tforal \tau\in T(A),
\eneq
there exists a unitary $u\in U(A)$ such that
\beq\label{ebot1-2}
\|[\psi(c),\,u]\|<\ep\tforal c\in {\cal F}\andeqn \rho_A\circ
\rm{bott}_1(\psi, u)(x_i)=\rho_A(y_i).
\eneq

\end{lem}

\begin{proof}
Let $f_i\in C$ be the identity of the $i$-th summand and let
$p_i=\psi(f_i).$ By considering each summand and $p_iAp_i,$ one
sees that one can reduce the general case to the case that
$C=PM_l(C(X))P,$ where $P\in M_l(C(X))$ is a projection and
$X=S^1,$ $X=T_{II,k}$ or $X=T_{III,k}$ for some integer $k>1.$

The case that $X=S^1$ follows from \ref{hitK}  immediately.

Since $K_1(C(T_{II,k}))=\{0\},$ it remains to consider the case
that $X=T_{III,k}.$ Let $v\in M_R(C)$ be a unitary such that
$K_1(C)=\Z/k \Z$ is generated by $[v].$ Suppose that $y\in K_0(A)$
such that there is a \hm\, $\kappa: K_1(C)\to K_0(A)$ with
\beq\label{ebot-1-2}
\kappa([v])=y.
\eneq
Since $k[v]=0,$ it follows that $ky=0.$ Therefore one must have
\beq\label{ebot-1-2+}
\rho_A(y)=0.
\eneq

The lemma then follows when one applies \ref{KTemb} to this
situation.

\end{proof}

\begin{lem}\label{2ebot}
Let $C\in {\cal N}$ be a  unital separable simple \CA\,  with
tracial rank zero, $A$ be  a  unital separable simple \CA\,  with
tracial rank zero and let $h: C\to A$ be a unital \hm.

Let  $G_0$ be a subgroup of $K_1(C)$ generated by
$x_1,x_2,...,x_m.$ Then, for any $\ep>0$ and any finite subset
${\cal F}\subset C,$ there exists $\dt>0$ satisfying the
following:

If $\af: G_0\to K_0(A)$ is a \hm\, such that
\beq\label{2ebot-1}
|\rho_A(\af(x_i))(\tau)|<\dt\tforal \tau\in
T(A),\,\,\,i=1,2,...,m,
\eneq
then there exists a unitary $u\in U(A)$ such that
\beq\label{2ebot-2}
\|[h(a),\, u]\|<\ep\tforal a\in {\cal F}\andeqn
\rho_A\circ\rm{bott}_1(h,u)(x_i)=\rho_A(\af(x_i))(\tau)
\eneq
for all $\tau\in T(A),$ $i=1,2,...,m.$
\end{lem}

\begin{proof}
It follows from \cite{Lnduke} that $C$ is a unital simple
AH-algebra with real rank zero and with no dimension growth. It
follows \cite{EG1} that $C=\lim_{n\to\infty}(A_n,\phi_n)$ where
each $A_n$ must be the form: $M_l(C(S^1))$ and $PM_l(C(X))P,$
where $P\in M_l(C(X))$ is a projection and $X$ has the form
$T_{II,k}$ or $T_{III,k}.$ Moreover, $\phi_n$ may be assumed to be
injective (see also \cite{EGL}).

There is a finitely generated subgroup $G_0'\subset K_1(A_n)$ for
some large $n$ such that $(\phi_{n,\infty})_{*1}(G_0')\supset
G_0.$ Thus, by considering $h\circ \phi_{n,\infty}$ and $A_n,$ we
effectively reduce the problem to the case that $C=A_n.$ The lemma
then follows from \ref{ebot1}.

\end{proof}

\begin{lem}\label{KK}
Let $C\in {\cal N}$ be a unital amenable separable \CA\, with
finitely generated $K_i(C)$ ($i=0,1$) and let ${\cal P}\subset
\underline{K}(C)$ be a finite subset.  Then, there is $\dt>0$ and
a finite subset ${\cal G}\subset C$ satisfying the following: If
$A$ is a unital \CA\, and if $h: C\to A$ is a unital \hm\, and
$u\in U(A)$ is a unitary such that
\beq\label{KK1}
\|[h(a), \, u]\|<\dt\tforal a\in {\cal G},
\eneq
there is an element $\kappa\in
Hom_{\Lambda}(\underline{K}(C\otimes C(S^1)),\underline{K}(A))$
such that
\beq\label{KK2}
[h]=\kappa|_{\underline{K}(C)}\andeqn
\rm{Bott}(h,u)=\kappa|_{\boldsymbol{\bt}(\underline{K}(C))}.
\eneq

\end{lem}

\begin{proof}
This follows from Proposition 2.4 and Lemma 2.8 of \cite{Lnhomp}.

\end{proof}

\begin{lem}\label{Memb}
Let $X$ be a finite CW complex and let $C=PM_l(C(X))P,$ where
$P\in M_l(C(X))$ is a projection. Let $x_1, x_2,...,x_m\in K_1(C)$
be generators of $K_1(C).$  Let $A$ be a unital simple separable
\CA\, with tracial rank zero and let $h: C\to A$ be a unital
monomorphism. For any $\ep>0$ and any finite subset ${\cal
F}\subset C,$ there is $\dt>0$ satisfying:

For any \hm\, $\af: K_1(C)\to K_0(A)$ with
\beq\label{Memb-1}
|\rho_A(\af(x_i))(\tau))|<\dt\tforal \tau\in
T(A),\,\,\,i=1,2,...,m,
\eneq
there exists a unitary $u\in U_0(A)$ such that
\beq\label{Memb-2}
\rm{bott}_1(h,u)(x_i)=\af(x_i),\,\,\,i=1,2,...,m
\eneq
and
\beq\label{Memb-2+}
\|[h(a),\, u]\|<\ep\tforal a\in {\cal F}.
\eneq

\end{lem}

\begin{proof}
By considering each summand, it is clear that we may assume that
$X$ is connected. We now write
\beq\label{Memb-3}
K_1(C)=\Z^k\oplus G_0,
\eneq
where $G_0$ is a finite subgroup and $k$ is a positive integer.

Let $\ep>0$ and a finite subset ${\cal F}\subset C$ be given.

 By
applying \ref{app}, without loss of generality, we may assume that
\beq\label{Memb-4}
h=\phi_0\circ h_1+\phi_1\circ h_2,
\eneq
where $\phi_0,$ $\phi_1,$ $h_1$ and $h_2$ are as described in
\ref{app}. Let $A_0,$ $A_1$ and $p\in A$ be as described in
\ref{app}.

Let ${\cal P}_0=\{x_1,x_2,...,x_m\},$ let $F_{01}$ be the subgroup
of $K_1(A_0)$ generated by $\{(h_1)_{*1}(x_i): i=1,2,...,m\}$ and
let $F_{02}$ be the subgroup generated by $\{(h_2)_{*1}(x_i):
i=1,2,...,m\}.$

Let ${\cal P}_1\subset K_0(C)$ be a finite subset which includes
$[1_C].$


There is $\af_i: K_1(A_0)\to K_0(A),$ $i=1,2,$ such that
\beq\label{Memb-4-1}
\af_1\circ (h_1)_{*1}|_{\Z^k}=\af|_{\Z^k}\andeqn \af_2\circ
(h_2)_{*2}|_{G_0}=\af|_{G_0}.
\eneq
Let $\ep_1$ be $\dt$  and ${\cal F}_0$ (in place of ${\cal F}$) in
\ref{ddbot} (for $A=C$).

Let $\dt_1>0$ be required by \ref{2ebot} for $A_0,$ $\phi_0,$ the
subgroup generated by $\{(h_1)_{*1}(x_i),
(h_1)_{*1}(x_2),...,(h_1)_{*1}(x_m)\},$ $\ep_1/3$ and ${\cal
F}_0\cup h_1({\cal F}).$ Let $\dt_2>0$ be required by \ref{2ebot}
for $A_1,$ $\phi_1,$ the subgroup generated by $\{(h_2)_{*1}(x_i),
(h_2)_{*1}(x_2),...,(h_2)_{*1}(x_m)\},$ $\ep_1/3$ and ${\cal
F}_0\cup h_2({\cal F}).$

Choose $\dt=\min\{\dt_1, \dt_2\}.$ Suppose (\ref{Memb-1}) holds.

By applying \ref{2ebot}, there are unitaries $u_1\in pAp$ and
$u_2\in (1-p)A(1-p)$ such that
\beq\label{Memb-4+1}
\|\phi_0(a), \,u_1]\|<\ep_1/2\tforal a\in {\cal
F}_1',\,\,\,\rho_A\circ \text{bott}_1(\phi_0,
u_1)|_{F_{01}}=\rho_A\circ \af_1|_{F_{01}},
\eneq
where ${\cal F}_1'=h_1({\cal F}),$ and
\beq\label{Memb-4+2}
 \|\phi_1(b),\,u_2]\|<\ep_1/2\tforal b\in {\cal F}_1''\andeqn
\rho_A\circ \text{bott}_1(\phi_1,u_2)|_{F_{02}}= \rho_A\circ
\af_2|_{F_{02}}.
 \eneq
Put $u'=u_1+u_2.$ It is then easy to see (with the assumption
(\ref{Memb-4}))that
\beq\label{Memb-4+3}
\|[h(a),\, u']\|<\ep_1/2\andeqn \rho_A\circ
\text{bott}_1(h,u')=\rho_A\circ \af.
\eneq
We also assume that $\text{Bott}(h,u')$ is well defined.

There is $\kappa_1\in Hom_{\Lambda}(\underline{K}(C),
\underline{K}(A))$ such that
\beq\label{Memb-5-1}
\kappa_1|_{K_1(C)}=-\af+\text{bott}_1(h,u')\andeqn
\kappa_1|_{K_0(C)}=\text{bott}_0(h,u').
\eneq
Define
 $\kappa\in
Hom_{\Lambda}(\underline{K}(C\otimes C(S^1)), \underline{A})$ such
that
\beq\label{Memb-5-2}
\kappa|_{\underline{K}(C)}=[h]\andeqn
\kappa({\boldsymbol{\bt}({K}(C)}))=\kappa_1.
\eneq

In particular,
\beq\label{Memb-5-3}
\rho_A(\kappa(\boldsymbol{\bt}(K_1(C))))=\{0\}.
\eneq

It follows from \ref{KTemb} that there is $u''\in U(A)$ such that
\beq\label{Memb-5}
\text{Bott}_1(h, u'')=\kappa|_{\boldsymbol{\bt}({K}(C))}
\eneq
Put $u=u'(u'')^*.$ Then, $u\in U_0(A).$ Moreover, one computes
that
\beq\label{Memb-6}
\text{bott}_1(h, u)&=&\text{bott}_1(h,u')-\text{bott}_1(h,u'')\\
&=&\text{bott}_1(h,u')+\af-\text{bott}_1(h,u')=\af.
\eneq

\end{proof}

\section{AH-algebras}

We will generalize the results in the previous section to the case
that $C$ is a unital AH-algebra. Let $A=\lim_{n\to\infty}(A_n,
\psi_n),$ where $A_n=P_nM_{l(n)}(C(X_n))P_n,$ where $X_n$ is a
finite CW complex, $P_n\in M_{l(n)}(C(X_n))$ is a projection. The
main issue is that the maps $\psi_n$ may not be injective so
$\phi_{n, \infty}(A_n)$ has the form $\phi_{n,
\infty}(P_n)M_{l(n)}(C(Y_n))\phi_{n, \infty}(P_n),$ where $Y_n$ is
a compact subset of $X_n$ which can no longer be assumed to be a
finite CW complex.

\begin{lem}\label{induc}
Let $X$ be a  finite CW complex, let $C=PM_l(C(X))P$  and let $A$
be a unital separable simple \CA\, with tracial rank zero. Suppose
that $h: C\to A$ is a unital \hm\, and $\phi: D\to A$ is the
monomorphism induced by $h,$ where $D=C/{\rm ker} h$ which has the
form $P_0M_l(C(F))P_0$ and $F$ is a compact subset of $X.$

Then, there exists a decreasing sequence of finite CW complex
$X_i\subset X$ which contains $F$ and a sequence of unital
monomorphism $h_n: C_n\to A,$ where $C_n=P_nM_l(C(X_i))P_n,$
$P_n=\pi^{(n)}(P)$ and $\pi^{(n)}: PM_l(C(X))P\to C_n$ is the
quotient map, such that
\beq\label{induc-1}
[\pi_n]\times [\phi]= [h_n]\,\,\,\text{in}\,\,\, KK(C_n, A)\andeqn
[h_{n}]=[s_n]\times [h_{n+1}]\,\,\,\text{in}\,\,\, KK(C_{n}, A),
\eneq
$n=1,2,...,$ where $\pi_n: C_n\to D$ and $s_n: C_n\to C_{n+1}$ are
the quotient maps. Moreover
\beq\label{induc-2}
&&\lim_{n\to\infty}\|h_n\circ \pi^{(n)}(f)-h(f)\|=0
\andeqn\\\label{induc-2+1}
&& \lim_{n\to\infty}\|h_n\circ
\pi^{(n)}(f)-h_{n+1}\circ s_{n}\circ \pi^{(n)}(f)\|=0
\eneq
for all $f\in C.$
\end{lem}

\begin{proof}
It is easy to check that the lemma for the general case follows
from the case that $C=M_l(C(X)).$ It is then easy to see that it
further reduces to the case that $C=C(X).$ It is also clear that
we may reduce the general case to the case that $X$ is connected.
Let $F$ be the spectrum of $h.$ So $F$ is a compact subset of $X.$
Put $D=C(F).$

Let $\{{\cal F}_n\}$ be an increasing sequence of finite subsets
of $C(X)$ whose union is dense in $C(X).$

Let $\sigma_n>0$ such that
\beq\label{induc-3}
|f(x)-f(x')|<{1\over{2^n}}\rforal f\in {\cal F}_n,
\eneq
provided that ${\rm dist}(x,x')<\sigma_n.$ One easily finds a
finite CW complex $X_n\subset X$ such that
\beq\label{induc-4}
X_n\subset \{x\in X: {\rm dist}(x, F)<\sigma_n\},\,\,\,n=1,2,....
\eneq
We may assume that $X_n\supset X_{n+1}.$ Note that
$\cap_{n=1}^{\infty} X_n=F.$

Let $\pi^{(n)}: C\to C_n,$ $\pi_n: C(X_n)\to C(F)$ and $s_n:
C(X_n)\to C(X_{n+1})$ be the quotient maps. Define $\phi_n:
C(X_n)\to A$  so that $\phi_n\circ\pi^{(n)}(f)=h(f)$ for all $f\in
C(X).$ One also has $\phi\circ \pi_n=\phi_n.$

Lemma 13.7 of \cite{Lnhomp} provides a monomorphism\, $h_n:
C(X_n)\to A$ such that
\beq\label{induc-5}
h_n\circ \pi^{(n)}\approx_{1/2^n} h\,\,\,\text{on}\,\,\,{\cal
F}_n.
\eneq
The limit formula (\ref{induc-2+1}) also follows.

Moreover, in the proof of Lemma 13.7 of \cite{Lnhomp} (see
(e13.29) and the notation there),
$$
[h_{00}^{(i)}]=[\phi_n|_{C_0(\Omega_i)}]\,\,\,\text{in}\,\,\,
KK(C_0(\Omega_i)),A).
$$

It follows that
\beq\label{induc-6}
[h_n]=[\phi_n]=[\pi_n]\times [h]\,\,\,\text{in}\,\,\, KL(C_n,A).
\eneq

\end{proof}

\begin{cor}\label{IND}
Let $C$ be a unital AH-algebra and let $A$ be a unital separable
simple \CA\, with tracial rank zero. Suppose that $h: C\to A$ is a
unital monomorphism.

Then, there exists a sequence of \CA s $C_n$ with the form
$C_n=P_nM_{l_n}(C(X_n))P_n,$ where $P_n\in M_{l_n}(C(X_n))$ is a
projection and $X_n$ is a finite CW complex satisfying the
following:
\beq\label{IND-1}
C=\lim_{n\to\infty}(C_n, \psi_n),
\eneq
where each $\psi_n$ is a unital \hm, and  there is a sequence of
unital monomorphisms $h_{n}: C_n\to A$ such that
\beq\label{IND-2}
[h_{n}]=[\psi_{n, \infty}]\times
[h|_{\psi_{n,\infty}(C_n)}]\,\,\,\text{in}\,\,\, KK(C_n, A)\andeqn
\eneq
\beq\label{IND-3}
[h_{n}]=[\psi_n]\times [h_{n+1}]\,\,\,\text{in}\,\,\, KK(C_{n},A),
\eneq
$n=1,2,....$ Moreover
\beq\label{IND-4}
&&\lim_{n\to\infty}\|h\circ \psi_{k,
\infty}(a)-h_{n}\circ \psi_{k,n}(a)\|=0\andeqn\\\label{IND-4+}
&&\lim_{n\to\infty}\|h_{n+1}\circ \psi_{k,n+1}(a)-h_{n}\circ \psi_{k,n}(a)\|=0
\eneq
for all $a\in C_k,$ $k=1,2,...$

\end{cor}

\begin{proof}
We may write $C=\lim_{n\to\infty}(B_n, \phi_n),$ where $C_n$ has
the form $B_n=P_nM_{k(n)}(C(X_n))P_n,$ where $X_n$ is a finite CW
complex, $k(n)\ge 1$ is an integer and $P_n\in M_{k(n)}(C(X_n))$
is a projection.  There exists a finite subset ${\cal F}_n\subset
B_n$ such that $\phi_n({\cal F}_n)\subset {\cal F}_{n+1}$ and
$\cup_{n=1}^{\infty}\phi_{n, \infty} ({\cal F}_n)$ is dense in
$C.$

Let $D_n=\phi_{n, \infty}(B_n).$ Then
$D_n=P_n'M_{k(n)}(C(F_n))P_n',$ where $F_n\subset X_n$ is a
compact subset of $X_n$ and $P_n'=\phi_{n,\infty}(P_n).$  By
\ref{induc}, there exists a finite CW complex $Y_n\subset X_n$
such that $Y_n\supset F_n$ and $C_n=Q_nM_{k(n)}(Y_n)Q_n$ and a
surjective \hm\, $s_n: B_n\to C_n$ and a unital monomorphism $h_n:
C_n\to A$ such that
\beq\label{IND-5}
[h_n]=[\gamma_n]\times [h|_{\gamma_n(C_n)}]\,\,\,\text{in}\,\,\,
KK(C_n, A),
\eneq
where $\gamma_n: C_n\to D_n$ is the quotient map.

Let $j_n: D_n\to D_{n+1}$ be the injective \hm. From the proof of
\ref{induc}, we may also assume that ${\rm ker} s_n\subset {\rm
ker}s_{n+1}\circ \phi_n.$ Thus we may assume that there is a \hm\,
$\psi_n: C_n\to C_{n+1}$ such that
\beq\label{IND-6}
s_{n+1}\circ \phi_n=\psi_n\circ s_n\andeqn j_n\circ
\gamma_n=\gamma_{n+1}\circ \psi_n \,\,\,n=1,2,....
\eneq
 Moreover,
\beq\label{IND-7}
h_n\approx_{1/2^n} h\circ \gamma_n\,\,\,\text{on}\,\,\,{\cal F}_n.
\eneq
From (\ref{IND-6}), one shows that $C=\lim_{n\to\infty}(C_n,
\psi_n).$ Moreover, $\psi_{n,\infty}=\gamma_n,$ $n=1,2,....$ It is
then easy to see from (\ref{IND-7}) that (\ref{IND-4}) and
(\ref{IND-4+}) hold. Moreover, (\ref{IND-2}) and (\ref{IND-3})
follow from (\ref{IND-5}).

\end{proof}

\begin{lem}\label{CMemb}
Let $C$ be a unital AH-algebra, let $A$ be a unital separable
simple \CA\, with tracial rank zero and let $h: C\to A$ be a
unital monomorphism.  Write $C=\lim_{n\to\infty} (C_n, \psi_n)$ so
that $C_n$ and $\psi_n$ satisfy the conditions in \ref{IND}.
 For any $\ep>0,$ any finite subset ${\cal F}\subset C$ and $x_1,x_2,...,x_k\in K_1(C),$ there
 exist $\dt>0,$ an integer $n\ge 1$ with
$ x_1,x_2,...,x_k\in (\psi_{n, \infty})_{*1}(K_1(C_n))$
 satisfying the following:

 For any $\alpha:
 (\psi_{n, \infty})_{*1}(K_1(C_n))\to K_0(A)$ with
 \beq\label{CMemb1}
 |\rho_A(\af\circ (\psi_{n, \infty})_{*1}(y_j))|<\dt\rforal j=1,2,...,k,
 \eneq
 where $\{y_1,y_2,...,y_m\}$ forms a set of generators of $K_1(C_n),$
 there is a unitary $u\in U_0(A)$ such that
 \beq\label{CMemb2}
 {\rm{bott}_1}(h, u)(x_i)=\af(x_i), \,\,\,i=1,2,...,k
\eneq
and
\beq\label{CMemb3}
\|[h(a), \, u]\|<\ep\rforal a\in {\cal F}.
\eneq

\end{lem}

\begin{proof}
Let $\{C_n\}$ and $\psi_n$ be as in \ref{IND}. We may assume that
$\{x_1,x_2,...,x_k\}\subset (\psi_{{n_0},
\infty})_{*1}(K_1(C_{n_0})).$ Without loss of generality, we may
assume that ${\cal F}_1\subset C_{n_0}$ is a finite subset so that
$\psi_{n_0, \infty}({\cal F}_1)\supset {\cal F}.$

There is $\eta>0$ and a finite subset ${\cal F}_2\subset C_{n_0}$
such that
$$
\text{bott}_1(h', v')=\text{bott}_1(h'',v''),
$$
provided that $\|[h'(a), \, v']\|<\eta,$ $\|[h''(a), v'']\|<\eta$
and
$$
\|h'(a)-h''(a)\|<\eta\tforal a\in  {\cal F}_2
$$
for any unitaries $v'$ and $v''$ and \hm s $h', h''$  from
$C_{n_0}$ above.

We now assume that ${\cal F}_2\supset {\cal F}_1.$

There is $n>n_0$ such that
\beq\label{CMemb-1}
h_n\approx_{\eta} h\circ \psi_{n,\infty}\,\,\,\text{on}\,\,\,{\cal
F}_2,
\eneq
where $h_n: C_n\to A$ is a unital monomorphism given by \ref{IND}.

Suppose that $\{y_1,y_2,...,y_m\}\subset K_1(C_n)$ is a set of
generators. Let $\ep/2>\eta>0.$ Let $\dt>0$ be as in \ref{Memb}
associated with $C_n,$ $y_1,y_2,...,y_m$ (in place of
$x_1,x_2,...,x_m$) and with $\eta$ (in place of $\ep$) and ${\cal
F}_2$ (in place of ${\cal F}$). By \ref{Memb}, there is a unitary
$u\in U_0(A)$ such that
\beq\label{CMemb-1+}
&&\|[h_n(a), u]\|<\ep/2\rforal a\in {\cal F}_1\andeqn\\\label{CMemb-2}
&&\text{bott}_1(h_n, u)(y_j)=\af\circ (\psi_{n, \infty})_{*1}(y_j),\,\,\,j=1,2,...,m.
\eneq
Let $x_j=(\psi_{n_0,\infty})_{*1}(z_j)$ for some $z_j\in
K_1(C_{n_0}),$ $j=1,2,...,m.$ By (\ref{CMemb-1+}) and the choice
of $\eta,$ we conclude  that
\beq\label{CMemb-4}
\text{bott}_1(h\circ \psi_{n_0,\infty},u)(z_j)=\af\circ
(\psi_{n_0, \infty})_{*1}(z_j), \,\,\,j=1,2,...,k.
\eneq
It follows that
\beq\label{CMemb-5}
\text{bott}_1(h, u)(x_i)=\af(x_i),\,\,\,i=1,2,...,k.
\eneq
Moreover, we have
\beq\label{CMemb-6}
\|[h(a),\, u]\|<\ep\tforal a\in {\cal F}.
\eneq

\end{proof}

\begin{lem}\label{CCM}
Let $C$ be a unital AH-algebra, let $A$ be a unital separable
simple \CA\, with tracial rank zero and let $h: C\to A$ be a
unital monomorphism.  Write $C=\lim_{n\to\infty} (C_n, \psi_n)$ so
that $C_n$ and $\psi_n$ satisfy the conditions in \ref{IND}.
 For any $\ep>0,$ any finite subset ${\cal F}\subset C$ and
 ${\cal P}\subset \underline{K}(C),$ there is
  an integer $n\ge 1$ with
$ {\cal P}\subset  [\psi_{n, \infty}](\underline{K}(C_n))$ and an
integer $k(n)\ge n$
satisfying the following:

 For  $\kappa\in Hom_{\Lambda}(\underline{K}(C_{k(n)}\otimes C(S^1)),\underline{K}(A))$
with
 \beq\label{CCM1}
 \kappa|_{\underline{K}(C_{k(n)})}=[h\circ \psi_{k(n),\infty}]\andeqn
 \rho_A\circ\kappa|_{{\boldsymbol{\bt}}(K_1(C_{k(n)}))}=0,
 \eneq
there is a unitary $u\in U(A)$ such that
 \beq\label{CCM2}
 {\rm{Bott}}(h\circ \psi_{n, \infty}, u)=
 \kappa\circ [\psi_{n, k(n)}]|_{\boldsymbol{\bt}(\underline{K}(C_n))},
\eneq
and
\beq\label{CCM3}
\|[h(a), \, u]\|<\ep\tforal a\in {\cal F}.
\eneq
\end{lem}

\begin{proof}
Let $\{C_n\}$ and $\psi_n$ be as in \ref{IND}. We may assume that
${\cal P}\subset [\psi_{n, \infty}](\underline{K}(C_n)).$ Without
loss of generality, we may assume that ${\cal F}_1\subset C_n$ is
a finite subset so that $\psi_{n, \infty}({\cal F}_1)\supset {\cal
F}.$ Let ${\cal G}\subset C_n$ be a finite subset and $\dt>0$ such
that
\beq\label{CCM-3-1}
\text{Bott}(h', u')=\text{Bott}(h'',u')
\eneq
for any pair of unital \hm s $h', h'': C_n\to A$ and any unitary
$u'\in A,$ provided that
$$
h'\approx_{\dt} h''\,\,\,\text{on}\,\,\,{\cal G},
$$
$$
\|[h'(a),\, u']\|<\dt\andeqn \|[h''(a),\, u']\|<\dt
$$
for all $a\in {\cal G}.$

We may assume that ${\cal G}\supset {\cal F}_1$ and $\dt<\ep/2.$
Choose a sufficiently large $k(n)\ge n$ such that $h_{k(n)}:
C_{k(n)}\to A$ is the monomorphism given by \ref{IND} so that
\beq\label{CCM-3-2}
h_{k(n)}\circ \psi_{n, k(n)}\approx_{\dt} h\circ
\psi_{n,\infty}\,\,\,\text{on}\,\,\,{\cal G}.
\eneq

 Then, by \ref{KTemb}, there exists a unitary
$u\in U(A)$ such that
\beq\label{CCM4}
\text{Bott}(h_{k(n)},
u)=\kappa|_{{\boldsymbol{\bt}}(\underline{K}(C_{k(n)}))}\andeqn
\|[h_{k(n)}(a),\, u]\|<\ep/2\rforal a\in {\cal F}_1.
\eneq
It follows from (\ref{CCM-3-2}) that we may assume that
\beq\label{CCM5}
\text{Bott}(h\circ \psi_{n,\infty},u)=\kappa\circ
[\psi_{n,k(n)}]|_{{\boldsymbol{\bt}}(\underline{K}(C_n))}\andeqn
\|[h(a),\, u]\|<\ep\rforal a\in {\cal F}.
\eneq
\end{proof}

\begin{lem}\label{CMhit}
Let $C$ be a unital AH-algebra, let $A$ be a unital separable
simple \CA\, with tracial rank zero and let $h: C\to A$ be a
unital monomorphism.  Write $C=\lim_{n\to\infty} (C_n, \psi_n)$ so
that $C_n$ and $\psi_n$ satisfy the conditions in \ref{IND}.
 For any $\ep>0,$ any finite subset ${\cal F}\subset C$ and ${\cal P}\subset \underline{K}(C),$ there is
 $\dt>0$ and
  an integer $n\ge 1$ with
$ {\cal P}\subset  [\psi_{n, \infty}](\underline{K}(C_n))$ and an
integer $k(n)\ge n$  satisfying the following:

 For  $\kappa\in Hom_{\Lambda}(\underline{K}(C_{k(n)}\otimes C(S^1)),\underline{K}(A))$
with
\beq\label{CMhit-1}
\kappa|_{\underline{K}(C_{k(n)})}=[h\circ \psi_{k(n),
\infty}]\andeqn
|\rho_A(\kappa({\boldsymbol{\bt}}(x_i))(\tau)|<\dt,
\,\,\,i=1,2,...,k
\eneq
for all $\tau\in T(A),$ where $\{x_1,x_2,...,x_k\}$ forms a set of
generators for $K_1(C_{k(n)}),$ then there exists a unitary $u\in
U(A)$ such that
\beq\label{CMhit-2}
\|[h(a),\, u]\|<\ep\tforal a\in {\cal F} \andeqn
{\rm{Bott}}(h\circ\phi_{n,
\infty},u)=\kappa\circ[\psi_{n,k(n)}]\circ{\boldsymbol{\bt}}.
\eneq

\end{lem}

\begin{proof}
This is a combination of \ref{CMemb} and \ref{CCM}. Fix $\ep>0,$
${\cal F}$ and ${\cal P}$  as in the lemma. Let $n$ and $k(n)\ge
n$ be in \ref{CCM} associated with $\ep/2$ and ${\cal F}.$

Let $\eta>0$ and ${\cal G}\subset C_{k(n)}$ be a finite subset
such that both $\text{Bott}(h', u')$ and $\text{Bott}(h'',u')$ are
well defined and
\beq\label{CMhit-3+}
\text{Bott}(h', u')=\text{Bott}(h'',u')
\eneq
for any pair of unital \hm s $h',\, h'': C_{k(n)}\to A$ and a
unitary $u'\in U(A)$ provided that
\beq\label{CMhit-4}
h'\approx_{\dt}h''\,\,\,\text{on}\,\,\,{\cal G}
\eneq
$$
\|[h'(a),\, u']\|<\eta \andeqn \|[h''(a), \, u']\|<\eta
$$
for all $a\in {\cal G}.$ We may assume that
$\psi_{k(n),\infty}({\cal G})\supset {\cal F}.$

Let $k$ be the largest order of finite order elements in
$K_i(C_{k(n)})$ ($i=0,1$). Let ${\cal Q}\in
F_{k!}\underline{K}(C_{k(n)})$ be a finite generator set (see
\ref{Kund}, \ref{Dbot2})). Let $\dt_0$ and ${\cal G}_1\subset
C_{k(n)}$ be a finite subset required by \ref{ddbot} associated
with ${\cal Q}.$ We may also assume that $\eta<\min\{\dt_0/4,
\ep/2\}$  and ${\cal G}\supset {\cal G}_1.$

 Choose $n_1\ge k(n)$ such that
\beq\label{CMhit-3}
h_{K}\approx_{\eta} h\circ \psi_{K,
\infty}\,\,\,\text{on}\,\,\,\psi_{k(n), K}({\cal G})
\eneq
for all $K\ge n_1.$ Let $\{z_1,z_2,...,z_m\}$ form a set of
generators of $K_1(C_{k(n)}).$ Let
$z_j'=(\psi_{k(n),\infty})_{*1}(z_j),$ $j=1,2,...,m.$ Let $\dt>0$
and $k(n)'\ge n_1$ (in place of $n$) associated with $\eta$ (in
place of $\ep$),
 $\psi_{k(n),\infty}({\cal G})$ (in place of ${\cal F}$) and $\{z_1',z_2',...,z_m'\}$ (in place
 of $x_1,x_2,...,x_k$) as in \ref{CMemb}.

Suppose that  $\kappa\in
Hom_{\Lambda}(\underline{K}(C_{k(n)'}\otimes
C(S^1)),\underline{K}(A))$ with
\beq\label{CMhit-4+}
\kappa|_{{\underline{K}}(C_{k(n)'})}=[h\circ
\psi_{k(n)',\infty}]\andeqn \vert
\rho_{A}(\kappa({\boldsymbol{\bt}}(x_j))(\tau)\vert<\dt\rforal
\tau\in T(A),
\eneq
$j=1,2,...,k,$ where $\{x_1,x_2,...,x_k\}$ forms a set of
generators of $K_1(C_{k(n)'}).$ Let $\af:
(\psi_{k(n)',\infty})_{*1}(K_1(C_{k(n)'}))\to K_1(A)$ be defined
by $\af\circ (\psi_{k(n)',\infty})_{*1}=\kappa|_{K_1(C_{k(n)'})}.$

It follows from \ref{CMemb} that there is a unitary $v\in U_0(A)$
such that
\beq\label{CMhit-5}
\|[h(a),\, v]\|&<&\eta\rforal a\in \psi_{k(n),\infty}({\cal G})\andeqn\\
\text{bott}_1(h,v)((\psi_{k(n),\infty})_{*1}(z_j))&=&\af((\psi_{k(n),\infty})_{*1}(z_j)))\\\label{CMhit-5+}
&=&\kappa\circ( (\psi_{k(n),k(n)'})_{*1}(z_j)),
\eneq
$j=1,2,...,m.$ By the choice of $\eta$ and ${\cal G},$ there
exists $\kappa'\in Hom_{\Lambda}(\underline{K}(C_{k(n)}\otimes
C(S^1)),\underline{K}(A))$ such that
\beq\label{CMhit-5++}
\text{Bott}(h\circ \psi_{k(n),\infty},
u)=\kappa'|_{{\boldsymbol{\bt}}(\underline{K}(C_{k(n)}))}.
\eneq
Let
\beq\label{CMhit-6}
\kappa_1=\kappa\circ[\psi_{k(n), k(n)'}\otimes {\rm
id}_{C(S^1)}]-\kappa'\in
Hom_{\Lambda}(\underline{K}(C_{k(n)}\otimes
C(S^1)),\underline{K}(A)).
\eneq
Since $\{z_1,z_2,...,z_m\}$ generates $K_1(C_{k(n)}),$ by
(\ref{CMhit-5+}),
\beq\label{CMhit-7}
\rho_A(\kappa_1\circ {\boldsymbol{\bt }}(K_1(C_{k(n)})))=0.
\eneq
Now by applying \ref{CCM}, we obtain a unitary $u_1\in U(A)$ such
that
\beq\label{CMhit-8}
\|[h(a),\, u_1]\|<\ep/2\rforal a\in {\cal F}\andeqn
\text{Bott}(h\circ \psi_{n,\infty},
u_1)=\kappa_1|_{{\boldsymbol{\bt}}(\underline{K}(C_n))}.
\eneq
Finally, put $u=vu_1.$ Then, since $\eta<\dt_0/4,$ by
(\ref{CMhit-8}),  (\ref{CMhit-6}), (\ref{CMhit-5++}) and
(\ref{CMhit-5+}),
\beq\label{CMhit-9}
\|[h(a), \, u]\|&<&\ep \rforal a\in {\cal F}\andeqn\\
 \text{Bott}(h\circ \psi_{n,\infty}, u)&=& \text{Bott}(h\circ \psi_{k(n),\infty}\circ \psi_{n, k(n)},v)+
 \text{Bott}(h\circ \psi_{n, \infty}, u_1)\\
 &=& \kappa'|_{{\boldsymbol{\bt}}([\psi_{n, k(n)}](\underline{K}(C_n))}
 +(\kappa_1)|_{{\boldsymbol{\bt}}([\psi_{n,k(n)}](\underline{K}(C_n))}\\
 &=&\kappa|_{{\boldsymbol{\bt}}([\psi_{n,k(n)'}](\underline{K}(C_n)))}.
 \eneq
We use $k(n)'$ for $k(n)$ in the statement.

\end{proof}

\begin{cor}\label{LKT}
Let $C$ be a unital AH-algebra with finitely generated $K_i(C)$
($i=0,1$) and let $A$ be a unital separable simple \CA\, with
tracial rank zero. Suppose that $h: C \to A$ is a unital
monomorphism. Suppose that $\kappa\in
Hom_{\Lambda}(\underline{K}(C\otimes C(S^1)), \underline{K}(A))$
such that
\beq\label{LKT-1}
\kappa|_{\underline{K}(C)}=[h]\andeqn
\rho_A(\kappa({{\boldsymbol{\bt}}(K_1(C))}))=0.
\eneq
Then, for any $\ep>0$ and any finite subset ${\cal F}\subset C,$
there exists a unitary $u\in U(A)$ such that
\beq\label{LKT-2}
\|[h(a), u]\|<\ep\tforal a\in {\cal F}\andeqn \rm{Bott}(h,
u)=\kappa|_{{\boldsymbol{\bt}}(\underline{K}(C))}.
\eneq

\end{cor}

\begin{proof}
Let $\ep>0,$ let ${\cal F}\subset C$ be a finite subset and let
$\kappa\in Hom_{\Lambda}(\underline{K}(C\otimes C(S^1)),
\underline{K}(A))$ satisfying (\ref{LKT-1}) be given. Write
$C=\lim_{n\to\infty}(C_n, \psi_n)$ as in \ref{IND}. Moreover, let
$h_n: C_n\to A$ be the monomorphism described in \ref{IND}. Since
$K_i(C)$ is finitely generated ($i=0,1$), there is $k\ge 1$ such
that
$$
Hom_{\Lambda}(F_k\underline{K}(C), F_k\underline{A}))=
Hom_{\Lambda}(\underline{K}(C), \underline{K}(A)).
$$
Let  ${\cal P}\subset \underline{K}(C)$  be a finite subset which
generates $F_k\underline{K}(C).$ We may well assume that ${\cal
P}\subset [\psi_{n, \infty}](\underline{K}(C_n)).$ Let $k(n)\ge n$
be as in \ref{CCM} (for $\ep,$ ${\cal F}$ and ${\cal P}$).

Define $\kappa_{k(n)}\in
Hom_{\Lambda}(\underline{K}(C_{k(n)}\otimes C(S^1)),
\underline{K}(A))$ as follows:
\beq\label{LTK-3}
(\kappa_{k(n)})|_{\underline{K}(C_{k(n)})}=[h\circ\psi_{k(n),
\infty}]\andeqn
(\kappa_{k(n)})|_{{\boldsymbol{\bt}}(\underline{K}(C_{k(n)}))}=
\kappa\circ {\boldsymbol{\bt}}([\psi_{k(n),
\infty}](\underline{K}(C_{k(n)})).
\eneq

Note that
$$
\rho_A(\circ \kappa_{k(n)}|_{{\boldsymbol{\bt}}(K_1(C_{k(n)}))}=0.
$$
We then apply \ref{CCM}.

*

\end{proof}

\section{AH-algebras with finitely generated $K$-theory}

\begin{lem}\label{LTRBOT}
Let $C$ be a unital AH-algebra
 and let $A$ be a unital separable  simple \CA\,
with tracial rank zero. Suppose that $\phi_1, \phi_2: C\to A$ are
two unital monomorphisms. Suppose that
\beq\label{ltrbot-1}
[\phi_1]&=&[\phi_2]\,\,\,\text{in}\,\,\,KL(C,A),\\
\tau\circ \phi_1&=&\tau\circ \phi_2\tforal \tau\in T(A)\andeqn\\
R_{\phi_1,\phi_2}(K_1(M_{\phi_1,\phi_2}))&=&\rho_A(K_0(A)).
\eneq
Then, for any increasing sequence of finite subsets $\{{\cal
F}_n\}$ of $C$ whose union  is dense in $C,$ any increasing
sequence of finite subsets $\{{\cal P}_n\}$ of $K_1(C)$ with
$\cup_{n=1}^{\infty}{\cal P}_n= K_1(C)$ and any decreasing
sequence of positive number $\{\dt_n\}$ with
$\sum_{n=1}^{\infty}\dt_n<\infty,$ there exists a sequence of
unitaries $\{u_n\}$ in $U(A)$ such that
\beq\label{ltrbot-2}
{\rm ad}\, u_n\circ \phi_1\approx_{\dt_n}
\phi_2\,\,\,\text{on}\,\,\, {\cal F}_n
\eneq
and
\beq\label{ltrbot-3}
\rho_A({\rm bott}_1(\phi_2, u_n^*u_{n+1})(x))=0\tforal x\in {\cal
P}_n
\eneq
for all sufficiently large $n.$
\end{lem}

\begin{proof}
By (\ref{ltrbot-1}),  it follows  from  \ref{CDT} (see 3.4 of
\cite{Lncd}) that there exists a sequence of unitaries
$\{u_n\}\subset A$ such that
\beq\label{LTM-3-}
\lim_{n\to\infty}{\rm ad}\, v_n\circ \phi_1(a)=\phi_2(a)\tforal
a\in C.
\eneq

Without loss of generality, we may assume that ${\cal
P}_n=\{z_1,z_2,...,z_n\}$ and we may also assume that ${\cal F}_n$
are in the unit ball of $C$ and $\cup_{n=1}^{\infty}{\cal F}_n$ is
dense in the unit ball of $C$ and we write that
$C=\lim_{n\to\infty}(C_n, \psi_n)$ (associated $h=\phi_1$)
satisfying the conclusion of \ref{IND}.

Let $1/2>\ep_n'>0$ so that $\text{Bott}(h', u')|_{{\cal P}_n} $ is
well defined for any unital monomorphism $h': C\to A$ and any
unitary $u'\in A$ for which
\beq\label{ltr-+}
\|[h'(a), u']\|<\ep_n'\tforal a\in {\cal F}_{n}.
\eneq
Moreover, we may assume that $\text{Bott}(h'\circ
\psi_{n,\infty},u')$ is well defined whenever (\ref{ltr-+}) holds.

Put $\ep_n''=\min\{\ep_n'/2, {1/2^{n+1}}, \dt_n/2\}.$

Let $\dt'_n>0$ be as $\dt$ in \ref{CMemb} and integer $k(n)$ (in
place of $n$) corresponding to $\ep_n''/2$ (in place of $\ep$),
${\cal F}_n$ (in place of ${\cal F}$) and ${\cal
P}_n=\{z_1,z_2,..., z_n\}$ (in place of $\{x_1,x_2,...,x_k\}).$
Note here we assume that $\{z_1,z_2,...,z_n\}\subset
(\psi_{k(n),\infty})_{*1}(K_1(C_{k(n)})).$ So there are $y_1,
y_2,...,y_n\in K_1(C_{k(n)})$ such that
$(\psi_{k(n),\infty})_{*1}(y_i)=z_i,$ $i=1,2,...,n.$

Suppose that $\{y_1,y_2,...,y_{m(n)}\}$ ($m(n)\ge n$) forms  a set of generators
of $K_1(C_{k(n)}).$  Put $\ep_n=\min\{\ep_n''/2, \dt_n'/2\}.$

We may assume that
\beq\label{ltrbot-4}
{\rm ad}\, v_n\circ \phi_1\approx_{\ep_n}
\phi_2\,\,\,\text{on}\,\,\, {\cal F}_{k(n)}
\eneq
Moreover we may assume that
\beq\label{ltrbot-5}
\|[\phi_1(a),\, v_nv_{n+1}^*]\|<\ep_n\tforal a\in {\cal F}_{k(n)}
\eneq
Thus $\text{bott}_1(\phi_1\circ \psi_{k(n),\infty}, v_nv_{n+1}^*)$
is well defined.






Let $w_1,w_2,...,w_n,...,w_{m(n)}\in M_l(C_{k(n)})$ for some $l\ge 1$ be
unitaries such that $[w_i]=(\psi_{k(n),\infty})_{*1}(y_i),$
$i=1,2,...,m(n).$  
 Put
$$
{\bar v}_n={\rm diag}(\overbrace{v_n,v_n,\dots, v_n}^{l}).
$$
We will continue to use $\phi_i$ for $\phi_i\otimes{\rm
id}_{M_l},$ $i=1,2.$

By choosing larger finite subset ${\cal F}_{k(n)}$  in (\ref{ltrbot-4}) and (\ref{ltrbot-5}),
we may also assume that, for $1\le j\le m(n),$
\beq\label{ltrbot-6}
\|\phi_2(w_j){\rm ad}\, {\bar
v}_n(\phi_1((w_j)^*))-1\|<(1/4)\sin(2\pi \ep_n),\,\,\,n=1,2,....
\eneq
Put
$$
h_{j,n}=\log({1\over{2\pi i}}\phi_2(w_j){\rm ad}\, {\bar
v}_n(\phi_1((w_j)^*))),\,\,\,j=1,2,...,n,...,m(n), n=1,2,....
$$
Then for any $\tau\in T(A),$
\beq\label{ltrbot-7}
\tau(h_{j,n})<\ep_n<\dt_n',\,\,\,j=1,2,...,m(n).
\eneq
By the assumption that
$R_{\phi_1,\phi_2}(M_{\phi_1,\phi_2})=\rho_A(K_0(A))$ and by
applying \ref{Rotrange}, we conclude that
\beq\label{ltrbot-8}
\widehat{h_{j,n}}(\tau)=\tau(h_{j,n})\in \rho_A(K_0(A)).
\eneq

By applying 6.1, 6.2 and 6.3  of \cite{Lnemb2}, and by passing to
a subsequence if necessary, we obtain a \hm\, $\af_n': K_1(C_n)\to
\rho_A(K_0(A))$ such that
\beq\label{ltrbot-9}
\af_n'(y_j)(\tau)=\widehat{h_{j,n}}(\tau)=\tau(h_{j,n}),
j=1,2,...,m(n).
\eneq

Since $\af_n'(K_1(C_n))$ is free, it follows from (\ref{ltrbot-8})
that there is a \hm\, $\af_n: K_1(C_n)\to K_0(A)$ such that
\beq\label{ltrbot-10}
\rho_A\circ \af_n(y_j)(\tau)=\tau(h_{j,n}),\,\,\,
j=1,2,...,m(n).
\eneq
It follows from \ref{CMemb} that there exists a unitary $U_n\in
U_0(A)$ such that
\beq\label{ltrbot-11}
&&\|[\phi_2(a),\, U_n]\|<\ep_n''/2\rforal a\in
{\cal F}_n\\\label{ltrbot-11+}
&&\hspace{-0.3in}\andeqn\,\,
\rho_A(\text{bott}_1 (\phi_2, U_n))(z_j)=-\rho_A\circ \af_n(z_j),
\eneq
$j=1,2,...,n.$

 By the Exel trace formula \ref{Exel} and (\ref{ltrbot-9}), we have
\beq\label{ltrbot-12}
\tau(h_{j,n})&=&-\rho_A(\text{bott}_1(\phi_2,
U_n))(z_j)(\tau)\\\label{ltrbot-12+}
&=&-\tau({1\over{2\pi i}}\log( {\bar U_n}\phi_2(w_j){\bar U_n}^*
\phi_2(w_j^*)))
\eneq
for all $\tau\in T(A),$ $j=1,2,...,m,$ where
\beq\label{ltrbot-13}
{\bar U}_n={\rm diag}(\overbrace{U_n,U_n\cdots, U_n}^l)
\eneq
 Define $u_n=v_nU_n,$ $n=1,2,....$
Put
$$
{\bar u_n}={\rm diag}(\overbrace{u_n,u_n,\dots, u_n}^l).
$$

By 6.1 of \cite{Lnemb2}), and by (\ref{ltrbot-11+}) and
(\ref{ltrbot-12+}), we compute that
\beq\label{ltrbot-14}
&&\hspace{-0.3in}\tau(\log({1\over{2\pi i}}(\phi_2(w_j){\rm ad}\, {\bar
u_n}(\phi_1(w_j^*)))))\\
&=&\tau(\log({1\over{2\pi i}}({\bar U_n}\phi_2(w_j^*){\bar
U_n}^*{\bar v}_n^*
(\phi_1(w_j)){\bar v_n})))\\
&=&\tau(\log({1\over{2\pi i}}({\bar U_n}\phi_2(w_j){\bar
U_n}^*\phi_2(w_j^*)\phi_2(w_j){\bar v}_n^*
(\phi_1(w_j^*)){\bar v_n})))\\
&=&\tau(\log({1\over{2\pi i}}({\bar U_n}\phi_2(w_j){\bar
U_n}^*\phi_2(w_j^*))))\\
&&\hspace{0.6in}+\tau({1\over{2\pi i}}\log(\phi_2(w_j){\bar v}_n^*
(\phi_1(w_j^*)){\bar v_n})))\\\label{ltrbot-15}
&=&\rho_A(\text{bott}_1(\phi_2, U_n))(z_j)(\tau)+\tau(h_{j,n})=0
\eneq
for all $\tau\in T(A).$


Let
\beq\label{ltrbot-16}
b_{j,n}&=&\log({1\over{2\pi i}}{\bar u_n}\phi_2(w_j){\bar
u_n}^*\phi_1(w_j^*))\andeqn\\
 b_{j,n}'&=&\log({1\over{2\pi
i}}\phi_2(w_j){\bar u_n}^*{\bar u_{n+1}}\phi_2(w_j^*){\bar
u_{n+1}}^*{\bar u_n})),
\eneq
$j=1,2,...,n$ and $n=1,2,....$
 We have, by (\ref{ltrbot-15}),
\beq\label{ltrbot-17}
\tau(b_{j,n})&=&\tau(\log({1\over{2\pi i}}{\bar
u_n}\phi_2(w_j){\bar
u_n}^*\phi_1(w_j^*)))\\
&=& \tau(\log{1\over{2\pi i}}{\bar u_n}^*{\bar
u_n}\phi_2(w_j){\bar u_n}^*\phi_1(w_j^*){\bar u_n})\\
&=&\tau(\log{1\over{2\pi i}}\phi_2(w_j){\bar u_n}^*\phi_1(w_j^*)
{\bar u_n})=0
\eneq
for all $\tau\in T(A),$ $j=1,2,...,n$ and $n=1,2,....$  Note also
$\tau(b_{j,n+1})=0$ for all $\tau\in T(A)$ and for $j=1,2,...,n.$
Note that
\beq\label{ltrbot-18}
{\bar u_n}^*e^{2\pi i b_{j,n}'}{\bar u_n}=e^{2\pi i
b_{j,n}}e^{-2\pi i b_{j,n+1}}.
\eneq
Thus, by 6.1 of \cite{Lnemb2} and by (\ref{ltrbot-18}), we compute
that
\beq\label{ltrbot-19}
\tau(b_{j,n}')=\tau(b_{j,n})-\tau(b_{j,n+1})=0\rforal \tau\in
T(A).
\eneq
It follows the Exel trace formula (\ref{Exel}) and
(\ref{ltrbot-19}) that
\beq\label{ltrbot-20}
&&\hspace{-0.3in}\rho_A(\text{bott}_1(\phi_2,u_n^*u_{n+1}))(z_j^*)(\tau)\\
&=&\tau(\log({\bar u_n}^*{\bar
u_{n+1}}\phi_2(w_j^*){\bar u_{n+1}}^*{\bar u_n}\phi_2(w_j))\\
&=&\tau(\log(\phi_2(w_j){\bar u_n}^*{\bar
u_{n+1}}\phi_2(w_j^*){\bar u_{n+1}}^*{\bar u_n}))=0
\eneq
for all $\tau\in T(A).$ It follows that (\ref{ltrbot-3}) holds. By
(\ref{ltrbot-4}) and (\ref{ltrbot-11}), one concludes that
(\ref{ltrbot-2}) also holds.

\end{proof}

\begin{thm}\label{LTM}
Let $C$ be a unital AH-algebra with finitely generated $K_i(C)$
($i=0,1$) and let $A$ be a unital separable  simple \CA\,
with tracial rank zero. Suppose that $\phi_1, \phi_2: C\to A$ are
two unital monomorphisms. Suppose that
\beq\label{LTM-1}
[\phi_1]&=&[\phi_2]\,\,\,\text{in}\,\,\,KL(C,A)\\
\tau\circ \phi_1&=&\tau\circ \phi_2\tforal \tau\in T(A)\andeqn\\
R_{\phi_1,\phi_2}(K_1(M_{\phi_1,\phi_2}))&=&\rho_A(K_0(A)).
\eneq
Then there exists a continuous path of unitaries $\{U(t):t\in
[1,\infty)\}$ of $A$ such that
\beq\label{LTM-2}
\lim_{t\to\infty}{\rm ad}\, U(t)\circ \phi_1(a)=\phi_2(a)\tforal
a\in C.
\eneq

\end{thm}

\begin{proof}
Let $\{{\cal F}_n\}$ be an increasing sequence of finite subsets
of $C$ which is dense in $C$ and $\{\ep_n\}$ be a decreasing
sequence of positive numbers such that
$\sum_{n=1}^{\infty}\ep_n<\infty.$

Let $\dt_n>0,$ ${\cal G}_n\subset C$ be a finite subset and ${\cal
P}_n\subset \underline{K}(C)$ be a finite subset required by
\ref{LNHOMP} associated with $A,$ $\phi_2,$ $\ep_n/2$ and ${\cal
F}_n.$

We may assume that ${\cal F}_n\subset {\cal G}_n$ and
$\dt_n>\dt_{n+1}$ and $\dt_n<\ep_n/2,$ $n=1,2,....$ Since $K_1(C)$
is finitely generated, let $\{z_1,z_2,...,z_k\}$ be a set of
generators of $K_1(C).$

Since $K_i(C)$ is finitely generated, there is $k\ge 1$ such that
$$
Hom_{\Lambda}(F_k\underline{K}(C),
F_k\underline{K}(A))=Hom_{\Lambda}(\underline{K}(C),
\underline{K}(A)).
$$
Let ${\cal Q}$ be a finite subset of generators of
$F_k\underline{K}(C).$ We may assume that ${\cal Q}\subset {\cal
P}_n,$ $n=1,2,...$

Let $\dt>0$ and ${\cal G}_0$ (in place of ${\cal F}$) be as in
\ref{ddbot} associated with $C$ and and $\dt_n<\dt/3.$

By \ref{LTRBOT}, there is a sequence of unitaries $u_n\in U(A)$
such that
\beq\label{LTM-3}
{\rm ad}\, u_n\circ
\phi_1\approx_{\dt_n/4}\phi_2\,\,\,\text{on}\,\,\,{\cal
G}_n\andeqn
\eneq
\beq\label{LTM-4}
\rho_A(\text{bott}_1(\phi_2,
u_n^*u_{n+1})(z_j))=0,\,\,\,j=1,2,...,k.
\eneq
Since $K_i(C)$ is finitely generated ($i=0,1$), we may also
assume, without loss of generality, that $\text{Bott}(\phi_2,
u_n^*u_{n+1})$ is well defined. Moreover, we may also assume, by
\ref{KK}, that there is $\kappa_n\in
Hom_{\Lambda}(\underline{K}(C\otimes C(S^1)), \underline{K}(A))$
such that
\beq\label{LTM-5}
(\kappa_n)|_{\underline{K}(C)}=[\phi_2]\andeqn
(\kappa_n)|_{{\boldsymbol{\bt}}(\underline{K}(C))}=
\text{Bott}(\phi_2, u_n^*u_{n+1}).
\eneq
By (\ref{LTM-4}),
\beq\label{LTM-6}
\rho_A(\kappa_n(K_1(C)))=0.
\eneq

We now construct $V_n$ as follows. Define $V_1=u_1.$ It follows
from \ref{LKT} that there exists a unitary $w_1\in U(A)$ such that
\beq\label{LTM-7}
\|[\phi_2(a),\, w_2]\|<\dt_2/4\tforal a\in {\cal G}_{2}\andeqn
\text{Bott}(\phi_2,w_2)=\kappa_2\circ \boldsymbol{\bt}.
\eneq
Define $V_2=u_2w_2^*.$ Then
\beq\label{LTM-18}
\|[\phi_2(a), \,V_1^*V_2]\|<\dt_1/2\tforal a\in {\cal F}_1.
\eneq
Moreover,
\beq\label{LTM-18+}
\text{Bott}(\phi_2, V_1^*V_2)&=&\text{Bott}(\phi_2, u_1^*u_2)
+\text{Bott}(\phi_2,w_2^*)\\
&=&\text{Bott}(\phi_2,u_1^*u_2)
-\text{Bott}(\phi_2,w_2)=0.
\eneq

Now by applying \ref{LNHOMP}, there is a continuous path of
unitaries $\{Z_1(t): t\in [0,1]\}$ such that
\beq\label{LTM-19}
Z_1(0)=1,\,\, Z_1(1)=V_1^*V_2\andeqn \|[\phi_2(a), \,
Z_1(t)]\|<\ep_1/2
\eneq
for all $a\in {\cal F}_1$ and for all $t\in [0,1].$

Suppose $W_n$ and $V_n$ have been constructed which satisfying the
following:
\beq\label{LTM-19+}
&&V_n=u_nW_n^*,\,\,\,\|[\phi_2(a),\, W_n]\|<\dt_n/2\tforal a\in {\cal G}_{n}\\
&&\andeqn \text{Bott}(\phi_2,W_n)=
\text{Bott}(\phi_2, V_{n-1}^*u_n).
\eneq
Moreover,
\beq\label{LTM-19++}
\rho_A(\text{bott}_1(\phi_2, W_n)))=0\andeqn \|[\phi_2(a),\,
V_{n-1}^*V_n]\|<\dt_{n-1}/2
\eneq
for all $a\in {\cal G}_{n-1}.$

 It follows from \ref{LKT} that there exists $W_{n+1}\in U(A)$
such that
\beq\label{LTM-20}
&&\|[\phi_2(a),\, W_{n+1}]\|<\dt_{n+1}/4\rforal a\in {\cal
G}_{n+1}\andeqn\\
&& \text{Bott}(\phi_2,
W_{n+1})=\kappa_{n+1}\circ
{\boldsymbol{\bt}}+\text{Bott}(\phi_2,W_n) .
\eneq
Note also that
$$
\rho_A(\text{bott}_1(\phi_2, W_{n+1})))=0.
$$
Define $V_{n+1}=u_{n+1}W_{n+1}^*.$ Then ($\dt/3<\dt_n$)
\beq\label{LTM-21}
\hspace{-0.5in}\text{Bott}(\phi_2,
V_n^*V_{n+1})&=&\text{Bott}(\phi_2, V_n^*u_{n+1})
+\text{Bott}(\phi_2,W_{n+1}^*)\\\label{LTM-21+}
\hspace{-0.5in}&=&-\text{Bott}(\phi_2,
W_n)+\text{Bott}(\phi_2,u_n^*u_{n+1}) -\text{Bott}(\phi_2,
W_{n+1})=0
\eneq
Moreover,
\beq\label{LTM-22}
\|[\phi_2(a), \,V_n^*V_{n+1}]\|<\dt_n/2\tforal a\in {\cal G}_{n}.
\eneq
By induction, we obtain a sequence of $\{W_n\}$ and $\{V_n\}$ with
$V_n=u_nW_n^*$ which satisfy (\ref{LTM-19}) and (\ref{LTM-19+}).

By applying ? of \cite{Lnhomp}, there is a continuous path of
unitaries $\{Z_n(t): t\in [0,1]\}$ such that
\beq\label{LTM-23}
&&Z_{n}(0)=1,\,\, Z_n(1)=V_n^*V_{n+1}\\
\andeqn &&\|[\phi_2(a), \, Z_n(t)]\|<\ep_n/2
\eneq
for all $a\in {\cal F}_n.$ Define $U(1)=V_1$ and
\beq\label{LTM-25}
U(t+n)=V_nZ_n(t)\tforal t\in [0,1),\,\,\, n=1,2,....
\eneq
We compute by (\ref{LTM-3}), (\ref{LTM-7}) and (\ref{LTM-22}) that
\beq\label{LTM-26}
\lim_{t\to\infty}{\rm ad}\,  U(t)\circ \phi_1(a)=\phi_2(a) \rforal
a\in C.
\eneq

\end{proof}

\section{Stable homotopy}

The results of this section will be applied to the case that $A=B\otimes C_k$ for a unital separable simple \CA\, $B$ and
a commutative \CA \, $C_k$ with torsion $K_0(C_k)$ and trivial $K_1(C_k).$

The following is a special case of Theorem 5.9 of \cite{Lnctaf} and it also follows from
a result of Dadarlat and Loring (\cite{DL1}).

\begin{thm}\label{KKHIT}
Let $X$ be a finite CW complex, let $k\ge 1$ be an integer and let
$P\in M_k(C(X))$ be a projection. Let $C=PM_k(C(X))P.$ Suppose
that $A$ is a unital \CA\, and suppose that  $x\in KK(C,\,A)$ and
a finite subset ${\cal P}\subset \underline{K}(C).$ Then, for any
$\ep>0$ and any finite subset ${\cal F}\subset C,$ there is an
integer $N_1>0$ and an $\ep$-${\cal F}$-multiplicative \morp\, $L:
C\to M_{N_1+1}(A)$ and a \hm\, $h_0: C\to M_{N_1}(\C)$ with finite
dimensional range such that
\beq\label{KKHIT1}
[L]|_{\cal P}=(x+[h_0])|_{\cal P}
\eneq
\end{thm}

By now several version (and more general version) of the next
theorem are known. In the following, the integer $N$ is allowed to
depend on $A$ as well as $L_1$ and $L_2$ among other things.

\begin{thm}\label{Uniq}
Let $C$ be as in {\rm \ref{KKHIT}}. For any $\ep>0$ and any finite
subset ${\cal F}\subset C.$ There exists $\dt>0,$  a finite subset
${\cal G}\subset C$ and a finite subset ${\cal P}\subset
\underline{K}(C)$ satisfying the following: Suppose that $A$ is a
unital \CA\, and $L_1, L_2: C\to A$ are two unital  $\dt$-${\cal
G}$-multiplicative \morp s such that
\beq\label{Uni1}
[L_1]|_{\cal P}=[L_2]|_{\cal P}.
\eneq
Then, there exists an integer $N\ge 1,$  and a unitary $U\in
U_0(M_{N+1}(A))$ such that
\beq\label{Uni2}
{\rm ad}\, U\circ ({\rm diag}(L_1, \Psi)\approx_{\ep} {\rm
diag}(L_2, \Psi)\,\,\,\text{on}\,\,\, {\cal F},
\eneq
where $\Psi: C\to M_N(\C)$ is a \hm.
\end{thm}

\begin{proof}
It is clear that the general case can be reduced to the case that
$C=C(X).$ Suppose that the theorem were false. We then obtain
$\ep_0>0$ and a finite subset ${\cal F}_0\subset C,$ a sequence of
finite subsets ${\cal P}_n\in \underline{K}(C)$ with ${\cal
P}_n\subset {\cal P}_{n+1}$ and $\cup_n {\cal
P}_n=\underline{K}(C),$ and a sequence of unital \CA s $\{A_n\},$
a sequence of unital \morp s $\{L_2^{(n)}\}$ and $\{L_1^{(n)}\}$
such that
\beq\label{Uni-1}
&&\hspace{-0.3in}\lim_{n\to\infty}\|L_j^{(n)}(a)L_j^{(j)}(b)-L_j^{(n)}(ab)\|=0\tforal a\, b\in C\\\label{Uni-1+}
&& [L_1^{(n)}]|_{{\cal P}_n}= [L_1^{(n)}]|_{{\cal
P}_n},\\\label{Uni-1++}
&&\hspace{-0.6in}\inf\{\sup\{\|{\rm ad}\,U\circ ({\rm diag}(L_1^{(n)}(a), \Psi_n(a))-
{\rm diag}(L_2^{(n)}(a), \Psi_n(a))\|: a\in {\cal F}_0\}\}\ge
\ep_0,
\eneq
where infimum is taken among all integer $k\ge 1,$ all possible
unital \hm s $\Psi_n: C\to M_k(\C)$ and all possible unitary $U\in
M_{k+1}(A_n).$ We may assume that $1_C\in {\cal F}_0.$

Define $B_n=A_n\otimes {\cal K},$ $B=\prod_{n=1}B_n$ and
$Q=B/\oplus_{n=1}B_n.$ Let $\pi: B\to Q$ be the quotient map.
Define $\phi_j: C\to B$ by $\phi_j(a)=\{\L_j^{(n)}(a)\}$ and
define ${\bar \phi}_j'=\pi\circ \phi_j,$ $j=1,2.$ Note that ${\bar
\phi}_j': C\to Q$ are \hm s. By 2.9 of \cite{GL1} (see also p.990
of \cite{GL2}), since $B_n$ is stable,
\beq\label{Uni-2}
&&K_i(B)=\prod_{n=1}K_i(B_n), \,\,\, K_i(B, \Z/k\Z)=\prod_{n=1}K_i(B_n, \Z/k\Z),\\
&& K_i(Q)=\prod_{n=1}K_i(B_n)/\bigoplus K_i(B_n)\andeqn\\
&&K_i(Q, \Z/k\Z)=\prod_{n=1}K_i(B_n,\Z/k\Z)/\bigoplus_{n=1}K_i(B_n, \Z/k\Z).
\eneq
Then (\ref{Uni-1+})  (and above) implies that
\beq\label{Uni-3}
[{\bar \phi_1}']=[{\bar \phi_2}']\,\,\,\text{in}\,\,\, KK(C, Q).
\eneq
Denote by $p_n$  the unit of ${\tilde B_n}$ and $q_n=p_n-1_{A_n}.$
Fix a point $\xi\in X$ define $\phi_0^{(n)}: C\to {\tilde B}_n$ by
$\phi_0^{(n)}(f)=f(\xi)q_n$ for $f\in C(X).$ Put ${\bar
\phi}_j={\bar \phi}_j'\oplus \pi\circ (\{\phi_0^{(n)}\}).$ Then
\beq\label{Uni-4}
[{\bar \phi}_1]=[{\bar \phi_2}]\,\,\,\text{in}\,\,\, KK(C, {\tilde
Q}).
\eneq

By \cite{D}, for any $\ep>0,$ there exists an integer $K>0,$ a
unitary $u\in M_{1+K}(Q)$ and a unital \hm\, $\psi: C\to
M_K(\C)\subset M_K({\tilde Q})$ such that
\beq\label{Uni-5}
{\rm ad}\, u\circ {\rm diag}({\bar \phi_1}, \psi)\approx_{\ep/4}
{\rm diag}({\bar \phi}_2,\psi)\,\,\,\text{on}\,\,\, {\cal F}_0.
\eneq
There is a unitary $V=\{V_n\}\in {\tilde B}$ such that $\pi(U)=u.$
It follows that (identifying $M_k(\C)$ with a \SCA\, of
$M_K({\tilde B}_n)$) for all sufficiently large $n,$
\beq\label{Uni-6}
{\rm ad}\, V_n\circ {\rm diag}(L_1^{(n)}\oplus \phi_0^{(n)},
\psi)\approx_{\ep/2}{\rm diag}(L_2^{(n)}\oplus \phi_0^{(n)}, \psi)
\eneq
on ${\cal F}_0.$  Put $E_n=1_{{\tilde B_n}}\oplus
1_{M_{1+K}({\tilde B}_n)}.$ Then, since $1_C\in {\cal F}_0,$
\beq\label{Uni-7}
\|[V_n, E_n]\|<\ep/2
\eneq
for all sufficiently large $n.$ Therefore, without loss of
generality, to simplify notation, we may assume that $V_n$ is a
unitary in $E_nM_{1+K}({\tilde B}_n)E_n$ and
\beq\label{Uni-8}
{\rm ad}\, V_n\circ {\rm diag}(L_1^{(n)}\oplus \phi_0^{(n)},
\psi)\approx_{\ep_0/4}{\rm diag}(L_2^{(n)}\oplus \phi_0^{(n)},
\psi)
\eneq
on ${\cal F}_0.$ Denote by $e_{n,k}'={\rm
diag}(\overbrace{1_{A_n},1_{A_n},...,1_{A_n}}^k)\in A_n\otimes
{\cal K}$ and
\beq\label{Uni-9}
e_{n,k}''={\rm
diag}(\overbrace{e_{n,k}',e_{n,k}',...,e_{n,k}'}^K)\in M_K( B_n).
\eneq
Note that $e_{n,k}''$ commutes with $\psi$ and $e_{n,k}'$ commutes
with $\phi_0^{(n)}.$ Put $e_{n,k}=e_{n,k}'\oplus e_{n,k}''$ in
$E_nM_{1+K}(B_n)E_n.$ Note that $V_n\in E_nM_{1+K}({\tilde
B}_n)E_n.$ Since $\{e_{n,k}\}$ forms an approximately central
approximate identity for $M_{1+K}(B_n),$
\beq\label{Uni-10}
\lim_{k\to\infty}\|[V_n,\, e_{n,k}]\|=0.
\eneq
Thus there is a unitary $U_{n,k}\in e_kM_{1+K}({\tilde B}_n)e_k$
for each $k$ and $n$ such that
\beq\label{Uni-11}
\lim_{k\to\infty}\|e_{n,k}V_ne_{n,k}-U_{n,k}\|=0.
\eneq
Note that there is an integer $N$ such that
$$
M_N(A_n)= ((e_{n,k}'-1_{A_n})\oplus e_{n,k}'')M_{1+K}(
B_n)((e_{n,k}'-1_{A_n})\oplus e_k'').
$$
We define
$\Psi_n(f)=(e_{n,k}'-1_{A_n})\phi_0^{(n)}(f)(e_{n,k}'-1_{A_n})\oplus
e_{n,k}\psi(f) e_{n,k}$ for $f\in C.$ Thus, for large $k,$ there
is a unitary $U_n \in M_{1+N}(A_n)$ such that
\beq\label{Uni-12}
{\rm ad}\, U_n\circ {\rm diag}(L_1^{(n)}, \Psi_n)\approx_{\ep_0/2}
 {\rm diag}(L_2^{(n)}, \Psi_n)\,\,\,\text{on}\,\,\,{\cal F}_0.
\eneq
This contradicts (\ref{Uni-1++}).

To see that one can choose $U\in U_0(M_{N+1}(A)),$ Write
\beq\label{UNi-13}
\Psi(f)=\sum_{j=1}^m f(\xi_j)p_j\tforal f\in C,
\eneq
where $\xi_1,xi_2,...,\xi_m$ are fixes points and $p_1,
p_2,...,p_m$ are mutually orthogonal projections. By replacing
$\Psi$ by $\Psi',$ where $\Psi'(f)=\Psi(f)\oplus f(\xi_1)\cdot
1_{M_{L}}$ for all $f\in C$ for some suitable integer $L\ge N+1,$
and by replacing $U$ by ${\rm diag}(U,1),$ we may assume that
$p_1=q_1\oplus 1_{M_{N+1}}.$ Define $U'=U{\rm diag}(1_{M_{N+1}},
U^*, 1_{M_{L-N-1}}).$ Then $U'\in U_0(M_{1+N'}(A)),$ where
$N'=N+L.$

\end{proof}

\begin{rem}

{\rm Note that $\Psi(f)=\sum_{j=1}^n f(\xi_j)p_j$ for all $f\in
C,$ where $\{p_1,p_2,...,p_n\}$ is a set of mutually orthogonal
projections and $\{\xi_1,\xi_2,...,\xi_n\}$ is a finite subset of
$X.$ In case that $C=M_k(C(X)),$  by choosing larger $N=Kk,$
adding another \hm\, with finite dimensional range (and replacing
$U$ by $U\oplus 1$), we can choose
$$
\Psi(f)={\rm diag}(f(\xi_1), f(\xi_2),...,f(\xi_K))
$$
for all $f\in C,$ where $\{\xi_1, \xi_2,...,\xi_K\}$ is some
$\ep$-dense subset of $X$ for some $\ep>0.$

}
\end{rem}

\begin{lem}\label{bothit}
Let $C$ be as in Theorem {\rm \ref{KKHIT}}.  For any $\ep>0,$ any
finite subset ${\cal F},$ a unital \hm\, $h: C\to A$ and
$\kappa\in Hom_{\Lambda}(\underline{K}(C), \underline{K}(SA)).$
Then there is an integer $N\ge 1,$ a unital \hm\, $h_0: C\to
M_N(\C)$ with finite dimensional range and a unitary $u\in
U(M_{1+N}(A))$ such that
\beq\label{bothit1}
\|[H(a), \, u]\|<\ep\rforal a\in {\cal F}\andeqn {\rm{Bott}}(H,u)
   =\kappa,
\eneq
where $H(c)={\rm diag}(h(c), h_0(c))$ for all $c\in C.$

\end{lem}

\begin{proof}
Let $S=\{z, 1_{C(S^1)}\}\subset C(S^1),$ where $z$ is the identity
function on the unit disk. Define $x\in
Hom_{\Lambda}(\underline{K}(C\otimes C(S^1)), A)$ as follows
$$
x|_{\underline{K}(C)}=[h]\andeqn
x|_{{\boldsymbol{\bt}}(\underline{K}(C)}=\kappa.
$$

Fix a finite subset ${\cal P}_1\subset
{\boldsymbol{\bt}}(\underline{K}(C) ).$ Choose $\ep_1>0$ and a
finite subset ${\cal F}_1\subset C$ satisfying the following:
\beq\label{bothit-1-1}
[L']|_{{\cal P}_1}=[L'']|_{{\cal P}_1}
\eneq
for any pair of $\ep_1$-${\cal  F}_1\otimes S$-multiplicative
\morp s $L', L'': C\otimes C(S^1)\to B$ (for some unital \CA\,
$B$), provided that
\beq\label{bothit-1-2}
L'\approx_{\ep_1} L''\,\,\,\text{on}\,\,\, {\cal F}_1\otimes S.
\eneq

Let $\ep>0$ and ${\cal F}\subset C$ be given finite subset. Put
$\ep_2=\min\{\ep/2, \ep_1\}$ and ${\cal F}_2={\cal F}\cup {\cal
F}_1.$

Let $\dt>0,$ ${\cal G}\subset C$ be a finite subset and ${\cal
P}\subset \underline{K}(C)$ be required by \ref{Uniq}  for
$\ep_2/2$ and ${\cal F}_2.$ Without loss of generality, we may
assume that ${\cal F}_2$ and ${\cal G}$ are in the unit ball of
$C$ and $\dt<\ep/2.$

Fix another finite subset ${\cal P}_1\subset (\underline{K}(C) )$
and define ${\cal P}_2={\cal P}\cup {\boldsymbol{\bt}}({\cal
P}_1)$ (as a subset of $\underline{K}(C\otimes C(S^1))$).

It follows from \ref{KKHIT} that there are  integers $N_1\ge 1,$ a
unital \hm\, $h_1: C\otimes C(S^1)\to M_{N_1-1}(\C)$ and a
$\dt/2$-${\cal G}\otimes S$-multiplicative \morp\, $L: C\otimes
C(S^1)\to M_{N_1}(A)$ such that
\beq\label{bothit-1}
[L]|_{{\cal P}_2}=(x+[h_1])|_{{\cal P}_2}.
\eneq
Define $H_1: C\to M_{N_1+1}(A)$ by
\beq\label{bothit-2}
H_1(f)=h(f)\oplus h_1(f\otimes 1)\tforal f\in C.
\eneq
Define $L_1: C\to M_{N_1}(A)$ by $L_1(f)=L(f\otimes 1)$ for all
$f\in C.$ Note that
\beq\label{bothit-4}
[L_1]|_{\cal P}=[H_1]|_{\cal P}
\eneq
It follows from \ref{Uniq} that there exists an integer $N_2\ge
1,$ a unitary $W\in M_{(N_2+1)N_1}(A)$ and a unital \hm\, $h_2:
C\to M_{N_2N_1}(\C)$ with finite dimensional range such that
\beq\label{botnit-5}
W^*(L_1(f)\oplus h_2(f))W\approx_{\ep/4} H_1(f)\oplus
h_2(f)\rforal f\in {\cal F}_2.
\eneq
Put $N=(N_2+1)N_1-1.$ Now define $h_0: C\to M_N(\C)$ and $H: C\to
M_{1+N}(A)$ by
\beq\label{bothit-6}
h_0(f)=h_1(f\otimes 1)\oplus h_2(f)\andeqn H(f)=h(f)\oplus
h_0(f)\rforal f\in C
\eneq
Define
\beq\label{bothit-7}
u=W^*(L(1\otimes z)\oplus 1_{M_{N_2N_1}})W.
\eneq
Then, by (\ref{botnit-5}),
\beq\label{bothit-8}
\|[H(f),\, u]\| &\le & \|(H(f)-{\rm ad}\, W\circ (L_1(f)\oplus
h_2(f))u\|+\\ &&\hspace{-0.6in}+\|[{\rm ad}\,W\circ (L_1(f)\oplus
h_2(f)),\,u]\|
+\|u(W\circ (L_1(f)\oplus h_2(f))-H(f))\|\\
&<& \ep/3+\dt/2+\ep/3<\ep
\rforal f\in {\cal F}_2.
\eneq
Define
$L_2: C\to M_{N+1}(A)$ by $L_2(f)=L_1(f) \oplus h_2(f)$ for $f\in
C.$
By (\ref{bothit-1}), we compute that
\beq\label{bothit-11}
\text{Bott}(H,\, u)|_{{\cal P}_1}
&=&\text{Bott}({\rm ad}\, W\circ L_2, \,u)|_{{\cal P}_1}\\
&=&\text{Bott}(L_2, \,L(1\otimes z)\oplus 1_{M_{N_2N_1}})|_{{\cal P}_1}\\
&=&\text{Bott}(L_1, L(1\otimes z))|_{{\cal P}_1}+\text{Bott}(h_2, 1_{M_{N_2N_1}})\\
&=&\text{Bott}(L_1, L(1\otimes z))|_{{\cal
P}_1}\\\label{bothit-11+}
&=&(x+[h_1])|_{{\boldsymbol{\bt}}({\cal P}_1)}= \kappa.
\eneq
Since $K_i(C)$ is finitely generated ($i=0,1$), by \ref{ddbot} and
by choosing a sufficiently large ${\cal P}_1,$ we may write from
(\ref{bothit-11+}) that
\beq\label{bothit-12}
\text{Bott}(H, u)=\kappa.
\eneq

\end{proof}

The following is a stable version of so-called Basic Homotopy
Lemma (see \cite{BEEK} and \cite{Lnhomp}).

\begin{thm}\label{Shomp}
Let $X$ be a finite CW complex, let $k\ge 1$ be an integer and let
$P\in M_k(C(X))$ be a projection. Let $C$ be a quotient of
$PM_k(C(X))P$ with the form $QM_k(C(F))Q,$ where $F$ is a compact
subset of $X.$ For any $\ep>0$ and any finite subset ${\cal
F}\subset C,$ there is $\dt>0$ a finite subset ${\cal G}\subset
C,$ a finite subset ${\cal P}\subset \underline{K}(C)$
 satisfying the following:

Suppose that $A$ is a unital \CA, suppose $h: C\to A$ is a unital
\hm\, and suppose that $u\in U(A)$ is a unitary such that
\beq\label{Shomp1}
\|[h(a), u]\|<\dt\tforal a\in {\cal G}\andeqn
{\rm{Bott}}(h,u)|_{\cal P}=\{0\}.
\eneq
Then there exists an integer $N\ge 1$ and a continuous path of
unitaries $\{U(t): t\in [0,1]\}$ in $M_{N+1}(A)$ such that
\beq\label{Shomp2}
U(0)=u',\,\,\, U(1)=1_{M_{N+1}(A)}\andeqn \|[h'(a),
u']\|<\ep\tforal a\in {\cal F},
\eneq
where
$$
u'={\rm diag}(u, H_0(1\otimes z))
$$
and $h'(f)=h(f)\oplus H_0(f\otimes 1)$ for $f\in C,$ where $H_0:
C\to M_N(\C)$ is a unital \hm\, (with finite dimensional range)
and $z\in C(S^1)$ is the identity function on the unit circle.

Moreover,
\beq\label{Shomp2+}
\rm{Length}(\{U(t)\})\le \pi+\ep.
\eneq

\end{thm}

\begin{proof}
Let $\ep>0$ and ${\cal F}\subset C$  be given. Without loss of
generality, we may assume that ${\cal F}$ is in the unit ball of
$C.$

Let $\dt_1>0,$  ${\cal G}_1\subset C\otimes C(S^1),$ ${\cal
P}\subset \underline{K}(C\otimes C(S^1))$  be required by
\ref{Uniq} for $\ep/4$ and ${\cal F}\otimes S.$ Without loss of
generality, we may assume that ${\cal G}_1={\cal G}_1'\otimes S,$
where ${\cal G}_1'$ is in the unit ball of $C$ and
$S=\{1_{C(S^1)}, z\}\subset C(S^1).$ Moreover, without loss of
generality, we may assume that ${\cal P}_1={\cal P}_2\cup{\cal
P},$ where ${\cal P}_2\subset \underline{K}(C)$ and ${\cal
P}\subset {\boldsymbol{\bt}}(\underline{K}(C)).$ Furthermore, we
may assume that any $\dt_1$-${\cal G}_1$-multiplicative \morp\,
$L'$ from $C$ to a unital \CA\, well defines $[L']|_{{\cal P}_1}.$

Let $\dt_2>0$ and ${\cal G}_2\subset C$ be a finite subset
required by 2.8 of \cite{Lnhomp} for $\dt_1/2$ and ${\cal G}_1'$
above.

Let $\dt=\min\{\dt_2/2, \dt_1/2, \ep/2\}$ and ${\cal G}={\cal
F}\cup {\cal G}_2.$

Suppose that $h$ and $u$ satisfy the assumption with above $\dt,$
${\cal G}$ and ${\cal P}.$ Thus, by 2.8 of \cite{Lnhomp}, there is
$\dt_1/2$-${\cal G}_1$-multiplicative \morp\, $L: C\otimes
C(S^1)\to A$ such that
\beq\label{Shomp4}
&&\|L(f\otimes 1)-h(f)\|<\dt_1/2\rforal f\in {\cal G}_1'\\
&& \|L(1\otimes z)-u\|<\dt_1/2.
\eneq

Define $y\in Hom_{\Lambda}(\underline{K}(C\otimes C(S^1)),
\underline{K}(A))$ as follows:
$$
y|_{\underline{K}(C)}=[h]|_{\underline{K}(C)}\andeqn
y|_{\boldsymbol{\bt}(\underline{K}(C))}=0.
$$
Then
\beq\label{Shomp4+1}
[L]|_{\cal P}=y|_{\cal P}.
\eneq

Define $H: C\otimes C(S^1)\to A$ by
$$
H(f\otimes g)=h(f)\cdot g(1)\cdot1_A
$$
for all $f\in C$ and $g\in C(S^1),$ where $S^1$ is identified with
the unit circle (and $1\in S^1$).

It follows that
\beq\label{Shomp3}
[H]|_{\cal P}= y|_{\cal P}=[L]|_{\cal P}.
\eneq

It follows from \ref{Uniq} that there is an integer $N\ge 1,$ a
unital \hm\, $H_0: C\otimes C(S^1)\to M_N(\C)$ with finite
dimensional range and a unitary $W\in U(M_{1+N}(A))$ such that
\beq\label{Shomp5}
W^*(H(f)\oplus H_{0}(f)W\approx_{\ep/4} L(f)\oplus H_{0}(f)
\eneq
for all $f\in {\cal F}\otimes S.$

Since $H_0$ has finite dimensional range, it is easy to construct
a continuous path $\{V': t\in [0,1]\}$ in a finite dimensional
\SCA\, of $M_N(\C)$  such that
\beq\label{Shomp6}
&&V'(0)=H_0(1\otimes z)),\,\,\,
V'(1)=1_{M_{N}(A)}\andeqn\\
&&H_0(f\otimes 1)V'(t)=V'(t)H_0(f\otimes 1)
\eneq
for all $f\in C$ and $t\in [0,1].$ Moreover,
\beq\label{Shomp6+}
\text{Length}(\{V'(t)\})\le \pi.
\eneq

Now define $U(1/4+3t/4)=W^*{\rm diag}(1, V'(t))W$ for $t\in [0,1]$
and
$$
u'=u\oplus  H_0(1_A\otimes z)\andeqn h'(f)=h(f)\oplus H_0(f\otimes
1)
$$
for $f\in C$ for $t\in [0,1].$  Then
\beq\label{Shomp7}
\|u'-U(1/4)\|<\ep/4\andeqn \|[U(t),\, h'(a)]\|<\ep/4
\eneq
for all $a\in {\cal F}$ and $t\in [1/4,1].$ The theorem follows by
connecting $U(1/4)$ with $u'$ with a short path as follows: There
is a self-adjoint element $a\in M_{1+N}(A)$ with $\|a\|\le {\ep
\pi\over{8}}$ such that
\beq\label{Shomp-8}
\exp(i a)=u'U(1/4)^*
\eneq
Define $U(t)=\exp(i (1-4t) a)U(1/4)$ for $t\in [0,1/4).$

\end{proof}

\begin{lem}\label{path}
There is $1/2>\dt>0$ satisfying the following: Let $A$ be a unital
\CA, $m\ge 1$ be an integer and
let $z\in U(M_m(A))$ be a unitary. Suppose that $V(t)\in
C([0,1-d], M_m(A))$ is a unitary such that
\beq\label{path-1}
 \|[V(1-d), \, z]\|<\dt\andeqn  V(0)=1
\eneq
for some $1/4>d>0.$  Suppose that there is a continuous path of
unitaries $W(t,s)\in C([0,1], M_m({\widetilde{SA}}))$ {\rm (}$s\in
[0,1]${\rm )} such that
\beq\label{path-2}
W(1,0)&=&1,\,\,\,\|W(t,0)-1\|<\dt\tforal t\in [1-d,1],\\
W(t,0)&=&z^*V(t)^*zV(t)\tforal t\in [0, 1-d]
\,\,\,and\,\,\,W(t,1)=1
\eneq
for all $t\in [0,1].$ Then,
\beq\label{path-3}
{\rm{bott}_1}(z, V(1-d))=0.
\eneq

\end{lem}

\begin{proof}
We  may view $W(t,0))$ as an element in $M_m({\widetilde{SA}}).$
So $W$ is a continuous path of unitaries in
$M_m({\widetilde{SA}})$ such that $W$ at $0$ is $W(t,0)$ and $W$
at $1$ is $1.$ In other words, $W(t,0)\in
U_0(M_m({\widetilde{SA}})).$

It is standard (see Lemma 9.3 of \cite{Lnhomp}, for a stronger
statement) that $W(t,s)$ can be replaced by one with additional
property: there is a positive number $M>0$ such that
\beq\label{path-4}
\|W(t,s)-W(t,s')\|\le M|s-s'|
\eneq
for all $s,s'\in [0,1]$ and $t\in [0,1].$ It follows that
$W(t,s)\in C([0,1]\times [0,1], M_m(A)).$

 Since $W(t,s)\in
M_m({\widetilde{SA}}),$ for each $s,$
$$
W(0,s)=e^{2\pi i \lambda(s)}=W(1,s)
$$
for some  continuous function $\lambda(s)\in C([0,1], M_m)_{s.a}$
with $\lambda(0)=1=\lambda(1).$

Replacing $W(t,s)$ by $e^{-2\pi i \lambda(s)} W(t,s),$ we may
assume that
\beq\label{path-4+1}
W(0,s)=W(1,s)=1.
\eneq
Thus we may assume that $W(t,s)\in C([0,1],
{\widetilde{SM_m(A)}}).$ Therefore, to simplify notation, by
denoting $M_m(A)$ by $A,$ we may assume $m=1$ for the rest of the
proof.

Put $Z(t,s)=zW(t,s).$ Then $Z(t,s)\in C([0,1]\times [0,1], A),$
\beq\label{path-5}
&&Z(1,0)=z,\,\,\,
\|Z(t,0)-z\|<\dt\rforal t\in [1-d,1],\\
&&Z(t,0)=V(t)^*zV(t)\tforal t\in [0,1-d],\,\,\,Z(t,1)=z\andeqn\\
&&Z(0,s)=Z(1,s)=z.
\eneq
Note that, for each $s\in [0,1],$ $z$ and $Z(t,s)$ are in
$$
\{f(t)\in C([0,1],A): f(0)=f(1)\}\cong C(S^1, A).
$$
Moreover,
\beq\label{path-6}
\|Z(t,s)-Z(t,s')\|\le M|s-s'|
\eneq
for all $t,s\in [0,1].$

Let $1>\ep>0.$ We assume that $\text{bott}_1(u,v)$ is well-defined
whenever $\|[u,\, v]\|<\ep.$ Moreover, we assume that if
$\|v-v'\|<\ep$ then $\text{bott}(u,v')=\text{bott}(u,v)$ (for any
such unitaries $u,v, v'$).

Let $\dt_1>0$ required by \ref{Shomp}  for $\ep/2$ above and for
$X=S^1.$

Put $\dt=\min\{\dt_1/4, \ep/4\}.$


 There is $\eta>0$ such that
\beq\label{path-6+1}
\|Z(t,s)-Z(t',s')\|<\dt/4,
\eneq
if $|t-t'|<\eta$ and $|s-s'|<\eta.$

Let $P_1:$
$$
0=s_0<s_1\cdots<s_n=1
$$ be a partition such that
\beq\label{path-7}
\|Z(t,s)-Z(t,s_j)\|&<&\dt/4\rforal t\in [0,1]
\eneq
if $s\in [s_{j-1}, s_j),$ $j=1,2,...,n.$ In particular, we assume
that $\max\{|s_j-s_{j-1}|: j=1,2,...,n\}<\eta.$

Note that, for each $s\in [0,1],$
\beq\label{path-7+}
[z]=[Z(t,s)]\,\,\,\ \text{in}\,\,\, K_1(S^1,A).
\eneq
 It follows from \ref{Uniq}  that there exists an integer $N>0$
such that for each $s\in (0,1),$ there exists $X(t,s)\in
U(M_{N+1}(C([0,1],A))$ such that
\beq\label{path-8}
X(t,s)^*z' X(t,s)&\approx_{\dt/4}& Z'(t,s)\andeqn\\
X(0,s)&=&X(1,s),
\eneq
where
\beq\label{path-9}
z'&=&{\rm diag}(z,e^{2\pi i /N},e^{4\pi i /N}\cdots e^{2N\pi i/N})\andeqn\\
Z'(t,s)&=&{\rm diag}(Z(t,s), e^{2\pi i /N},e^{4\pi i /N}\cdots
e^{2N\pi i/N}).
\eneq

Since $Z'(t,s)$ is uniformly continuous on $[0,1]\times [0,1],$
there exists $1>T_0>0$  such that
\beq\label{path-9+}
\|Z'(t,s)-z'\|<\dt/4\rforal 0<t\le T_0 \andeqn \rforal s\in [0,1].
\eneq
Note also
\beq\label{path-9++}
 \|Z'(t,s)-z'\|<\dt\rforal 1-d<t<1\andeqn
\rforal s\in [0,1].
\eneq
By choosing even smaller $T_0>0$  we may also assume that
\beq\label{path-9+2}
\|V(t)-1\|&<&\dt/4\tforal 0<t\le T_0.
\eneq

Define
\beq\label{path-9+3}
Y'(t,s)=X(T_0,s)^*X(t,s).
\eneq

Put
\beq\label{path-10}
V'(t)={\rm diag}(V(t),\overbrace{1,1,\cdots,1}^N)\rforal t\in
[0,1-d].
\eneq

Note that, by (\ref{path-8}) and (\ref{path-7}),
\beq\label{path-11}
&&\| V'(t)^*z'V'(t)X(t,s_1)^*z'^*X(t,s_1)-1\|<\dt/4+\dt/4=\dt/2
\eneq
for all  $t\in [0,1-d]$ and  by (\ref{path-7})
\beq\label{path-11+}
\|X(t,s_j)^*z'X(t,s_j)X(t,s_{j+1})^*z'^*X(t,s_{j+1})-1\|<\dt/4
\eneq
for all $t\in [0,1]$ and $j=1,2,...,n-1.$ It follows  that
\beq\label{path-12}
\|X(t,s_1)V'(t)^*z' V'(t)X(t,s_1)^*-z'\|<\dt/2
\eneq
 for all $t\in
[0,1-d]$ and
\beq\label{path-12+}
\|X(t,s_{j+1})X(t,s_j)^*z'X(t,s_j)X(t,s_{j+1})^*-z'\|<\dt/4
\eneq
for all $t\in [0,1].$
 It follows from (\ref{path-12}), (\ref{path-8}) and (\ref{path-9+}) that
\beq\label{path-13}
\| Y'(t,s_1)V'(t)^*z'V'(t)Y'(t,s_1)^*-z'\|<\dt/2+\dt/4
\eneq
for all $t\in [0,1-d]$ and
\beq\label{path-13+}
\|Y'(t,s_{j+1})Y'(t,s_{j})^*z'Y'(t,s_{j})Y'(t,s_{j+1})^*-z'\|<(\dt/2+\dt/4)+\dt/2=5\dt/4
\eneq
for all $t\in [0,1].$
 It follows from (\ref{path-9+2}) and
(\ref{path-9+3}) that
\beq\label{path-14}
\text{bott}_1(z', Y'(T_0,s_1)(V(T_0)')^*)=\text{bott}_1(z',
(V(T_0)')^*)=0.
\eneq
Since $Y'(t,s_1),V'(t,s_1)\in C([0,1-d], M_{N+1}(A))$ we conclude
(by (\ref{path-13})) and (\ref{path-14}) that
\beq\label{path-15}
\text{bott}_1(z',Y'(t,s_1)V'(t)^*)=\text{bott}_1(z',
Y'(T_0,s_1)(V(T_0)')^*)=0
\eneq
for all $t\in [0,1-d].$

We also have
\beq
&&\hspace{-0.2in}\text{bott}_1(z', Y'(T_0,s_{j+1})Y'(T_0,s_j)^*)\\\label{path-16}
&=& \text{bott}_1(z',
Y'(T_0,\,s_{j+1}))-\text{bott}_1(z',\,Y'(T_0,s_j))=0.
\eneq
Since $Y'(t,s_{j+1})Y'(t,s_j)^*$ is continuous on $[0,1],$ by
(\ref{path-13+})and (\ref{path-16}),
\beq\label{path-17}
\text{bott}_1(z', Y'(t,s_{j+1})Y'(t,s_j)^*)=0\tforal t\in [0,1].
\eneq
Thus, by \ref{Shomp}, there is an integer $K\ge 1$ and , for each
$j,$ there is a continuous path of unitaries $\{U_j(1-d,s): s\in
[0,1]\}$ in $M_{(1+N)(K+1)}(A)$ such that
\beq\label{path-18}
&&U_j(1-d,1)=Y(1-d,s_{j+1})Y(1-d,s_j)^*,\,\,\,U_j(1-d,0)=1\andeqn\\
&&\|[z'',\,U_j(1-d,s)]\|<\ep/2
\eneq
for all $s\in [0,1],$ $j=1,2,...,n,$ also,
\beq\label{path-19}
&& U_0(1-d,1)=Y(1-d,s_1)V''(1-d)^*,\,\,\, U_0(1-d,0)=1\andeqn\\
&&\|[z'',\, U_0(1-d,s)]\|<\ep/2
\eneq
for all $s\in [0,1],$ where
\beq\label{path-19+1}
z''&=&{\rm diag}(z',\Phi(z\otimes 1)),\\
V''(1-d,s_j)&=&{\rm
diag}(V'(1-d,s),\overbrace{1_{M_{1+N}},1_{M_{1+N}},...,1_{M_{1+N}}}^{K}),
\andeqn\\
 Y(1-d, s)&=&{\rm
diag}(Y'(1-d,s),\Phi(1\otimes z)),
\eneq
where  $\Phi: C(S^1\times S^1)\to M_{(1+N)K}(\C)$ is a unital \hm.

Define $U(t, s): s\in [0, n]$ as follows
\beq\label{path-20}
U(t,s)=U_0(t,s)V(t),\,\,\,\, U(t,s+j)=U_j(t,s)Y(t,s_j)\tforal s\in
[0,1],
\eneq
 $j=1,2,...,n-1.$
Note that
\beq\label{path-21}
\|[z'',\, U(1-d,s)]\|<\ep\rforal s\in [0,n].
\eneq
We have
\beq\label{path-22}
U(1-d,n)=Y(1-d,1)\andeqn U(1-d,0)=V''(1-d)
\eneq
for $t\in [0,1).$ Moreover,
\beq\label{path-23}
\|[z'',\, Y(t,1)]\|<\ep\rforal t\in [0,1]\andeqn \|[z'',\,
V''(1-d)]\|<\ep.
\eneq
Thus (by (\ref{path-21}) and (\ref{path-22})),
\beq\label{path-24}
\text{bott}_1(z'',V''(1-d)) &=& \text{bott}_1(z'', Y(1-d,1))\\
&=&\text{bott}_1(z'', Y(T_0,1))=0.
\eneq

 It follows that
\beq\label{path-25}
\text{bott}_1(z, V(1-d))=0.
\eneq

\end{proof}

\begin{lem}\label{zero}
There is $1/4>\dt>0$ satisfying the following:
 Let $A$ be a unital \CA\, and let $u\in U(A)$ be a unitary.
Suppose that $w(t)\in {\widetilde{SA}}$ is a unitary such that
\beq\label{zero1}
w(1)=w(0)=1_A,  w(t)=u^*v(t)^*uv(t) \rforal t\in [0,1-d]
\eneq
for some $0<d<1/2,$ where $v(t)\in C([0,1-d],A)$ and
\beq\label{zero2}
\|w(1)-w(t)\|<\dt \rforal t\ge 1-d.
\eneq
Suppose also  that
\beq\label{zero3}
\rm{bott}_1(u, v(1-d))=0.
\eneq
Then $[w]=0$ in $K_1(SA).$

\end{lem}

\begin{proof}
Let $1/4>\ep>0.$ Let $\dt_1>0$ and $N$ be as in \ref{Shomp}
required for $\ep.$ Let $\dt=\dt_1/2.$ Suppose that $u, w$ and $v$
are as in the lemma for $\dt.$ From inequality (\ref{zero2}) and
(\ref{zero3}), by applying   \ref{Shomp}, we obtain a continuous
path of unitaries $\{V(t): t\in [1-d, 1]\}\subset M_{1+N}(A)$ such
that
\beq\label{zero-1}
V(1-d)=\begin{pmatrix} v(1-d) & 0  &0 &\ldots & 0  \\
                    0 & \omega_N &0&\ldots &0 \\
                    0 & 0& \omega_N^2 &\ldots &0\\
                    \vdots &\vdots& &\ddots & \vdots\\
                    0 &0 &0 &\ldots & \omega_N^N
                    \end{pmatrix}\,,
                    \eneq
where $\omega_N=e^{2\pi /N},$
\beq\label{zero-2}
V(1)=1_{M_{1+N}(A)}\andeqn \|[V(t), \,u']\|<\ep
\eneq
for all $t\in [1-d, 1],$ where
$$
u'=\begin{pmatrix} u & 0  &0 &\ldots & 0  \\
                    0 & \omega_N &0&\ldots &0 \\
                    0 & 0& \omega_N^2 &\ldots &0\\
                    \vdots &\vdots& &\ddots & \vdots\\
                    0 &0 &0 &\ldots & \omega_N^N
                    \end{pmatrix}.
          $$
define
\beq\label{zero-2+}
V(t)=\begin{pmatrix} v(t) & 0  &0 &\ldots & 0  \\
                    0 & \omega_N &0&\ldots &0 \\
                    0 & 0& \omega_N^2 &\ldots &0\\
                    \vdots &\vdots& &\ddots & \vdots\\
                    0 &0 &0 &\ldots & \omega_N^N
                    \end{pmatrix}
                    \eneq
for $t\in [0,1-d].$ Define $ W(t)=(u')^*V(t)^*u'V(t) $ and
$W'(t)=w(t)\oplus 1_{M_N(A)}.$

Put  $T_s\in M_2(M_{1+N}(A))$ as follows:
\beq\label{zero-3}
T_s=\begin{pmatrix}\cos ({\pi s\over{2}}) & -\sin ({\pi s\over{2}})\\
                     \sin({\pi t\over{2}}) & \cos ({\pi t\over{2}})
                     \end{pmatrix}.
                     \eneq
Since $V(0)=V(1)=1,$  for each $s\in [0,1],$
\beq\label{zero-4}
\hspace{-0.2in}W_s(t)=\begin{pmatrix} u' & 0\\
                                      0 &1\end{pmatrix}  T_s\begin{pmatrix}V(t)^* & 0\\
                                                                                                   0 & 1\end{pmatrix}
                                                     T_s^{-1} \begin{pmatrix} u' &0\\
                                                                                              0 &
                                                                                              1\end{pmatrix}^*
                                                              T_s\begin{pmatrix}V(t)& 0\\
                                                                                              0 & 1\end{pmatrix} T_s^{-1}\in M_{2(1+N)}({\widetilde{SA}}).
                                                                                              \eneq
  Moreover,
  \beq\label{zero-5}
  W_0(t)=W(t)\oplus 1_{M_N(A)}\andeqn W_1(t)= 1_{M_{2N}(A)}.
    \eneq
It follows that $[W(t)]=0$ in $K_1(SA).$

On the other hand, for $t\in [0, 1-d),$
\beq\label{zero-6-}
W(t)=W'(t)
\eneq
and for $t\in [1-d,1],$
\beq\label{zero-6}
\|W'(t)-W(t)\|\le
\|W'(t)-1\|+\|1-(u')^*V(t)^*u'V(t)\|<\dt+\ep<1/2.
\eneq
Therefore
$$
[W'(t)]=0\,\,\,\text{in} \,\,\,K_1(SA).
$$
It follows that $[w(t)]=0$ in $K_1(SA).$

\end{proof}

\section{The Main Theorem}

In \ref{LTM}, $K_i(C)$ is assumed to be finitely generated
($i=0,1$). Therefore $KK(C,A)=KL(C,A).$ In particular, one may
write that $K_1(M_{\phi_1,\phi_2})=K_0(A)\oplus K_1(C).$ Moreover,
one may write $K_1(C)=Tor(K_1(C))\oplus G_1$ for some finitely
generated abelian free group $G_1.$ Note that $Tor(K_1(C))\subset
{\rm ker}R_{\phi_1, \phi_2}.$It is then easy to check that
$R_{\phi_1, \phi_2}(K_1(M_{\phi_1,\phi_2}))=\rho_A(K_0(A))$
implies that ${\tilde \eta}_{\phi_1,\phi_2}=0.$ Theorem \ref{LTM}
is a special case of \ref{TM} below and the proof of \ref{TM} does
not depend on that of \ref{LTM}. However, \ref{LTM} is handy as
soon as \ref{LTRBOT} and \ref{LKT} are available. It perhaps also
helps to understand the additional problems involved in the proof of
\ref{TM}.  Note that, the approximation argument in the proof of
\ref{LTM} considers a finite subset of $C$ at a time. So as in
\ref{CDT}, it suffices to consider $KL(C,A).$ It is therefore
interesting to understand that this could not be done for the
general situation. However, using $KK(C,A)$ (instead of $KL(C,A)$)
enables us to have a consistent information on $KK(C_{n+1}, A).$
It should also be noted that the full power of ${\tilde
\eta}_{\phi_1,\phi_2}=0$ will be used. In other words, we have to
use more information about  rotation maps than what is used in the
proof of \ref{LTRBOT}. It probably worth to point out that the
proof would be simpler if  \hm s $\psi_n$ from $C_n$ to $C_{n+1}$ is injective.
Among other things, one could apply \ref{LTM}.
Unfortunately, in general, AH-algebra can not be written such
a way so that each \hm\, from $C_n$ to $C_{n+1}$ is injective.

\begin{lem}\label{PPVuV}
Let $C$and $A$  be unital \CA s and let $\phi_1,\, \phi_2: C\to A$
be two unital monomorphism. Let $(a_{ij})\in
M_n(M_{\phi_1,\phi_2}).$ Suppose that $\ep>0$ is a positive number
and $c_{ij}\in C$ such that
\beq\label{PPVuV1}
\|\pi_0((a_{ij}))-(\phi_1(c_{ij}))\|<\ep.
\eneq
Then
\beq\label{PPVuV2}
\|\pi_1((a_{ij}))-(\phi_2(c_{ij}))\|<\ep.
\eneq
\end{lem}

\begin{proof}
This is almost evident. In fact there are $(b_{ij})\in M_n(C)$
such that
\beq\label{PPVuV-1}
\pi_0((a_{ij}))=\phi_1((b_{ij})).
\eneq
Here we continue to use $\phi_1$ for $\phi_1\otimes {\rm
id}_{M_n}.$ Then this and (\ref{PPVuV1}) imply that
\beq\label{PPVuV-2}
\|(b_{ij})-(c_{ij})\|<\ep,
\eneq
since $\phi_1$ is injective.  But by (\ref{PPVuV-1})
$$
\pi_1((a_{ij}))=(\phi_2(b_{ij})).
$$
Using $\phi_2$ for $\phi_2\otimes {\rm id}_{M_n},$ it then follows
from (\ref{PPVuV-2}) that
$$
\|\pi_1((a_{ij}))-\phi_2((c_{ij}))\|<\ep.
$$

\end{proof}

\begin{NN}\label{PreVuV}
{\rm Now let $A$ be a unital AH-algebra and $B$ be a unital \CA.
Suppose that $\phi_1, \phi_2: A\to B$ are two unital
monomorphisms.

Suppose that  $[\phi_1]=[\phi_2]$ in $KK(A,B).$ Then the following
 exact sequence splits:
$$
\begin{array}{ccccccc}
0 \to  & \underline{K}(B)  &\to &
\underline{K}(M_{\phi_1, \phi_2})
&{\stackrel{\pi_0}{\rightleftarrows}}_{\theta}& \underline{K}(A) &\to 0
\end{array}
$$

}
\end{NN}

\begin{lem}\label{VuV}
Let $A$ and $B$ and $\phi_1,$ and $\phi_2$ be as in \ref{PreVuV}.
Moreover, we assume that ${\tilde \eta}(\eta)(\phi_1,\phi_2)=0.$

Let $C=PM_k(C(X))P,$ where $X$ is a finite CW complex, $k\ge 1$ is
an integer, $P\in M_k(C(X))$ is a projection,  and let $\psi: C\to
A$ be a unital \hm.

For any $1/2>\ep>0,$ any finite subset ${\cal F}\subset C$ and any
finite subset ${\cal P}\subset \underline{K}(C),$ there are
integers $N_1\ge 1,$  an $\ep/2$-${\cal F}$-multiplicative \morp\,
$L: \psi(C)\to M_{1+N_1}(M_{\phi_1,\phi_2}),$  a unital \hm\,
$h_0: \psi(C)\to M_{N_1}(\C),$
and a continuous path of unitaries $\{V(t): t\in [0,1-d]\}$ of
$M_{1+N_1}(B)$ for some $1/2>d>0,$ such that $[L]|_{\cal P}$ is
well defined, $V(0)=1_{M_{1+N_1}(B)},$
\beq\label{VuV1}
[L\circ \psi]|_{\cal P}=(\theta\circ [\psi]+[h_0\circ
\psi])|_{\cal P},
\eneq
\beq\label{VuV2}
\pi_t\circ L\circ \psi\approx_{\ep} {\rm ad}\, V(t)\circ
((\phi_1\circ \psi)\oplus (h_0\circ
\psi))\,\,\,\text{on}\,\,\,{\cal F}
\eneq
for all $t\in (0,1-d],$
\beq\label{VuV3}
\pi_t\circ L\circ \psi\approx_{\ep} {\rm ad}\, V(1-d)\circ
((\phi_1\circ \psi)\oplus (h_0\circ
\psi))\,\,\,\text{on}\,\,\,{\cal F}
\eneq
for all $t\in (1-d,1),$  and
\beq\label{VuV33}
\pi_1\circ L\circ \psi\approx_{\ep}\phi_2\circ \psi\oplus h_0\circ
\psi\,\,\,\text{on}\,\,\,{\cal F},
\eneq
where $\pi_t: M_{\phi_1,\phi_2}\to B$ is the point-evaluation at
$t\in (0,1).$

\end{lem}

\begin{proof}
Let $\ep>0$ and let ${\cal F}\subset C$ be a finite subset.

Let $\dt_1>0,$ ${\cal G}_1\subset C$ be a finite subset and ${\cal
P}\subset \underline{K}(C)$ be a finite subset required by
\ref{Shomp} for $\ep/4$ and ${\cal F}$ above.

Let $\ep_1=\min\{\dt_1/2, \ep/4\}$ and ${\cal F}_1={\cal F}\cup
{\cal G}_1.$ We may assume that ${\cal F}_1$ is in the unit ball
of $C.$ We may also assume that $[L']|_{\cal P}$ is well defined
for any $\ep_1$-${\cal F}_1$-multiplicative \morp\, from $C$ to
(any unital \CA).

Let $\dt_2>0$ and ${\cal G}\subset C$ be a finite subset and
${\cal P}_1\subset \underline{K}(C)$ be finite subset required by
\ref{Uniq}  for $\ep_1/2$ and ${\cal F}_1.$ We may assume that
$\dt_2<\dt_1/2,$ ${\cal G}\supset {\cal F}_1$ and ${\cal
P}_1\supset {\cal P}.$ We also assume that ${\cal G}$ is in the
unit ball of $C.$

It follows from \ref{KKHIT} that there exists an integer $K_1\ge
1,$ a unital \hm\, $h_0': \psi(C)\to M_{K_1}(\C)$ and a
$\dt_2/2$-${\cal G}$-multiplicative \morp\, $L_1:  \psi(C)\to
M_{K_1+1}(M_{\phi_1, \phi_2})$ such that
\beq\label{VuV-1}
[L_1\circ \psi]|_{{\cal P}_1}=(\theta\circ [\psi]+[h_0'\circ
\psi])|_{{\cal P}_1}.
\eneq

Note that $\pi_0: M_{\phi_1, \phi_2}\to B$ has an image in
$\phi_0(A).$ We also have, by viewing $\pi_0$ as a \hm\, from
$M_{\phi_1,\phi_2}$ to $A,$
\beq\label{VuV-2}
[\pi_0]\circ \theta\circ [\psi]=[\psi]\,\,\,\text{in}\,\,\, KK(C,
A).
\eneq

Moreover, by viewing $\pi_0$ as a map from $M_{\phi_1,\phi_2}$ to
$B,$
\beq\label{VuV-3}
[\pi_0]\circ \theta\circ
[\psi]=[\phi_1\circ\psi]\,\,\,\text{in}\,\,\, KK(C, B).
\eneq
Furthermore, for each $t\in (0,1],$
\beq\label{VuV-4}
[\pi_t]\circ \theta\circ [\psi]=[\phi_1\circ \psi].
\eneq


By \ref{Uniq}, we obtain an integer $K_0,$ a unitary $V_{00}\in
U(M_{1+K_1+K_0}(\phi_1(A)))$ and a unital \hm\, $h: C\to
M_{K_0}(\C)$ such that
\beq\label{VuV-5--1}
{\rm ad}\, V_{00}\circ (\pi_0\circ L_1\oplus h)\approx_{\ep_1/2}
(\phi_1\oplus h_0'\oplus h)\,\,\,\text{on}\,\,\,{\cal F}_1.
\eneq
Write $V_{00}=(\phi_1(a_{ij}))$ for $a_{ij}\in A.$  Since $\phi_1$
is injective, it follows that $(a_{ij})\in U(M_{K_1+K_0}(A)).$ Let
$V_{00}'=\phi_2((a_{ij})).$ The assumption that
$[\phi_1]=[\phi_2]$ implies that $[V_{00}]=[V_{00}']$ in $K_1(B).$
By adding another $h$ in (\ref{VuV-5--1}) and replacing $V_{00}$
by $V_{00}\oplus 1_{M_{K_0}},$ if necessary, we may assume that
$V_{00}$ and $V_{00}'$ are in the same component of
$U(M_{1+K_1+K_0}(B)).$
 One obtains a continuous path of unitaries $\{Z(t): t\in [0,1]\}$
in $M_{1+K_1+K_0}(B)$ such that
\beq\label{VuV-5--2}
Z(0)=V_{00}\andeqn Z(1)=V_{00}'.
\eneq
It follows that $Z\in M_{1+K_1+K_0}(M_{\phi_1,\phi_2}).$ By
replacing $L_1$ by ${\rm ad}\, Z\circ (L_1\oplus h)$ and using a
new $h_0',$ we may assume that
\beq\label{VuV-5--3}
\pi_0\circ L_1\approx_{\ep_1/2}\phi_1\oplus
h_0'\,\,\,\text{on}\,\,\,{\cal F}_1.
\eneq
It follows from \ref{PPVuV} that we may also assume that
\beq\label{VuV-5--4}
\pi_1\circ L_1\approx_{\ep_1/2} \phi_2\oplus
h_0'\,\,\,\text{on}\,\,\, {\cal F}_1.
\eneq

There is a partition:
\beq\label{VuV-5}
0=t_0<t_1<\cdots <t_n=1
\eneq
such that
\beq\label{VuV-6}
\pi_{t_{i}}\circ L_1\approx_{\dt_2/2}\pi_t\circ
L_1\,\,\,\text{on}\,\,\,{\cal G}
\eneq
for all $t_i\le t\le t_{i+1},$ $i=1,2,...,n-1.$

 By applying \ref{Uniq}  again, we
obtain an integer $K_2\ge 1, $
a unital \hm\, $h_{00}: C\to M_{K_2}(\C),$
 and  a unitary $V_i\in M_{1+K_1+K_2}(B)$ such that
\beq\label{VuV-7}
{\rm ad}\, V_i\circ (\phi_1\oplus h_0'\oplus
h_{00})\approx_{\ep_1/2} (\pi_{t_i}\circ L_1 \oplus
h_{00})\,\,\,\text{on}\,\,\, {\cal F}_1.
\eneq
Note that, by (\ref{VuV-6}), (\ref{VuV-7}) and (\ref{VuV-5--3})
$$
\|[\phi_1\oplus h_0'\oplus h_{00}(a), V_1]\|<\dt_2/2+2\ep_1\rforal
a\in {\cal F}_1.
$$
Put $K_3=K_1+K_2$ and $H_0=\phi_1\oplus h_0'\oplus h_{00}.$

It follows from \ref{bothit} that there is an integer $J_1\ge 1,$
a unital \hm\, $F_1: C\to M_{J_1}(\C)$ and a unitary $W_1\in
U(M_{K_3+1+J_1}(B))$ such that
\beq\label{VuV-8}
\|[H_1(a), \, W_1]\|<\dt_2/2\rforal a\in {\cal F}_1\andeqn
\text{Bott}(H_1, W_1)=\text{Bott}(H_0, V_1),
\eneq
where $H_1=\phi_1\oplus h_0'\oplus h_{00}\oplus F_1.$ Define
$V_1'=W_1^*{\rm diag}(V_1, 1_{M_{J_1}}).$ Note that
\beq\label{VuV-9}
\|[H_1(a), \,V_1']\|<\dt_2+2\ep_1\rforal a\in {\cal F}_1\andeqn
\text{Bott}(H_1, V_1')=0.
\eneq
Therefore,
\beq\label{VuV-11-}
\|[H_1(a),\, W_1^*{\rm diag}(V_1, 1_{M_{J_1}}){\rm diag}(V_2^*,
1_{M_{j_1}})]\|<\dt_2+2\ep_1\rforal a\in {\cal F}_1.
\eneq
It follows from \ref{bothit} that there is an integer $J_2\ge 1,$
a unital \hm\, $F_2: C\to M_{J_2}(\C)$ and a unitary $W_2\in
U(M_{J_2+J_1+K_3+1}(B))$ such that
\beq\label{VuV-11}
\|[H_2(a), \, W_2]\|<\dt_2/2\rforal a\in {\cal F}_1\andeqn
\text{Bott}(H_2, W_2)=\text{Bott}(H_1(a),{\bar  V_1'}{\bar
V_2}^*),
\eneq
where $H_2=H_1\oplus F_2,$ ${\bar V_1'}={\rm diag}(V_1',
1_{M_{J_2}(B)})$ and ${\bar V_2}={\rm diag}(V_2,
1_{M_{J_1+J_2}(B)}).$ Thus
\beq\label{VuV-12}
\|[H_2(a),\,{\bar V_1'}{\bar V_2}^*W_2^*]\|<\dt_2+2\ep_1\rforal
a\in {\cal F}_1\andeqn
\eneq
\beq\label{VuV-13}
\text{Bott}(H_2,{\bar V_1'}{\bar V_2}^*W_2 ^*)=0.
\eneq
Put $V_2'=W_2{\bar V_2}.$ By continuing this process, we obtain an
integer $N_0\ge 1,$ a unital \hm\, $F_0: C\to M_{N_0}(\C),$
and unitaries $v_i\in M_{1+K_2+N_0}(B)$ such that
\beq\label{VuV-14}
{\rm ad}\, v_i\circ (\phi_1\oplus h_0'\oplus
F_{00})\approx_{\dt_2+2\ep_1} \pi_{t_i}\circ (L_1\oplus
F_{00})\,\,\,\text{on}\,\,\,{\cal F}_1,
\eneq
\beq\label{VuV-15}
&&\|[\phi_1\oplus  h_0'\oplus F_{00}(a),\,v_iv_{i+1}^*]\|<\dt_2+2\ep_1\rforal a\in {\cal F}_1\\
&&\text{Bott}(\phi_1\oplus h_0'\oplus F_{00}, v_1)=0\andeqn \text{Bott}(\phi_1\oplus h_0'\oplus F_{00}, v_iv_{i+1}^*)=0.
\eneq
It follows from \ref{Shomp} that there is an integer $N_1'\ge 1,$
a unital \hm\, $F_0': C\to M_{N_1'}(\C)$ and a continuous path of
unitaries $\{w_i(t): t\in [t_{i-1}, t_i]\}$ such that
\beq\label{VuV-17}
&&\hspace{-0.2in}w_1(0)=1, \,\,\,w_1(t_1)=v_1',\,\,\, w_i(t_{i-1})=v_{i-1}'(v_{i}')^*, w_i(t_i)=1, \,\,\, i=2,3,...,n\andeqn\\
&&\|[\phi_1\oplus h_0'\oplus F_{00}\oplus F_{0}'(a), \,w_i(t)]\|<\ep/2\rforal a\in {\cal F}, \,\,\,i=1,2,...,n,
\eneq
where $v_i'={\rm diag}(v_i, 1_{M_{N_1'}}(B)),$ $i=1,2,...,n.$
Define $V(t)=w_i(t)v_i'$ for $t\in [t_{i-1}, t_i],$ $i=1,2,...,n.$
Then $V(t)\in C([0,1-d], M_{N_1'}(B)).$ Moreover,
\beq\label{VuV-18}
{\rm ad}\, V(t)\circ (\phi_1\oplus h_0'\oplus F_{00}\oplus
F_0')\approx_{\ep/2} \pi_t\circ L_1\oplus F_{00}\oplus F_0'
\,\,\,\text{on}\,\,\,{\cal F}.
\eneq

Define  $h_0=h_0'\oplus F_{00}\oplus F_0',$ $L=L_1\oplus
F_{00}+F_0'$ and $d=1-t_{n-1}.$ Then, by (\ref{VuV-18}),
(\ref{VuV2}) and (\ref{VuV3}) hold. From (\ref{VuV-5--4}),
(\ref{VuV33}) also holds.

\end{proof}

Let $C$ be a unital separable \CA\, and $A$ be a unital \CA. Suppose that $\phi_1, \phi_2: C\to A$
are two unital monomorphisms.

Define
$$
M_{\phi_1, \phi_2}=\{f\in C([0,1],A): f(0)=\phi_1(a)\andeqn
f(1)=\phi_2(a)\,\,\,{\rm for\,\,\,some}\,\,\,a \in C\}.
$$
Denote by $\pi_t: M_{\phi_1,\phi_2}\to A$ the point-evaluation at
$t$ ($t\in [0,1]$). We will also use $\pi_t$ for the \hm\, from
$M_n(M_{\phi_1, \phi_2})\to M_n(A).$ Note also $\pi_0$ gives a
surjective \hm\, from $M_{\phi_1, \phi_2}$ onto $C$ and ${\rm
ker}\pi_0=SA.$

\begin{NN}\label{0AKK}

{\rm
 Let $C$ be a unital separable amenable \CA\, in ${\cal N}$ with finitely generated $K_i(C)$ ($i=0,1$), let $A$ be a unital separable \CA\, and let $\phi_1, \phi_2: C\to A$ be two unital \hm s.
 In what follows, we will continue to use $\phi_1$ and $\phi_2$ for the induced \hm s from $M_k(C)$ to $M_k(A).$
Suppose that $v\in U(A)$ and
\beq\label{0AKK1}
\|v^*\phi_1(z_j)v-\phi_2(z_j)\|<\ep<1/2,\,\,\,j=1,2,...,n
\eneq
for some $z_1, z_2,...,z_n\in U(M_k(C)).$ Define $W_j(t)\in
U(M_2(C([0, 1], M_k(A))$ as follows
\beq\label{0AKK2}
W_j(t)=(T_tVT_t^{-1})^*{\rm diag}(\phi_1(z_j),
1_{M_k})T_tVT_t^{-1},
\eneq
where
$$
V={\rm diag}(v, 1_{M_k})\andeqn
T_t=\begin{pmatrix}\cos ({\pi t}) & -\sin ({\pi t})\\
                     \sin({\pi t}) & \cos ({\pi t})
                     \end{pmatrix}
                     $$

Note that $W_j(0)={\rm diag}(v^*\phi_1(z_j)v, 1)$ and $W_j(1)={\rm
diag}(\phi_1(z_j),1).$ Connecting $W_j(0)$ with ${\rm
diag}(\phi_2(z),1)$ by a short path, we obtain a continuous path
of unitary $Z_j(t)$ such that $Z_j(0)={\rm diag}(\phi_2(z_j),1),$
$Z(1/4)=W(0)$ and $Z_j(1)={\rm diag}(\phi_1(z_j),1)$ and
$\|Z_j(t)-Z_j(1/4)\|<1/2$ for $t\in [0,1/4).$ Thus $Z\in M_{2k}(M_{\phi_1,\phi_2}).$
Note that, if $Z_j\in M_{2k}(M_{\phi_1,\phi_2})$ such that
$Z_j'(0)=Z_j(0),$ $Z_j'(t)=Z_j(t)$ for all $t\in [1/4, 1]$ and
$$
\|Z_j'(t)-Z_j'(1/4)\|<1/2
$$
for all $t\in [0,1/4].$ Then $(Z_j')^*Z_j\in U_0(M_2(M_{\phi_1,
\phi_2})).$ Thus the map $\gamma: K_1(C)\to K_1(M_{\phi_1,
\phi_2})$ defined by $\gamma([z_j])=[Z_j]$ is well defined. One
easily verifies that, with sufficiently small $\ep,$ since
$K_1(C)$ is finitely generated, $\gamma$ defines a \hm. Let
\beq\label{0AKK3}
h_j={\rm diag}(\log (\phi_2(z_j)^*V^*\phi_1(z_j)V,
1),\,\,\,j=1,2,...,n
\eneq
We may specifically use
\beq\label{0AKK4}
Z_j(t)={\rm diag}(\phi_2(z_j), 1)\exp(i 4t h_j)\rforal t\in
[0,1/4].
\eneq

If we use $\phi_1$ and $\phi_2$ for the induced \hm s from
$M_k(C\otimes C')$ to $A\otimes C',$ for some commutative \CA\,
$C'$ with finitely generated $K_i(C')$ ($i=0,1$), we obtain a
\hm\, $\gamma: K_1(C\otimes C')\to K_1(M_{\phi_1,\phi_2}\otimes
C')$ provided that $\ep$ is small (and $z_j$ are in $M_k(C\otimes
C')$).

Let ${\cal F}\subset C$ be a finite subset and $\ep>0.$ Suppose
that there is a unitary $v\in U(A)$ such that
\beq\label{0AKK5}
{\rm ad}\, v\circ \phi_1\approx_{\ep}
\phi_2\,\,\,\text{on}\,\,\,{\cal F}.
\eneq
Let $U(t)'=T_tVT_t^{-1}.$ Define
\beq\label{0AKK6}
L(c)(t)=(U({4t-1\over{3}})')^*{\rm diag}(\phi_1(c),0)U(({4t-1\over{3}})'\tforal
t\in [1/4, 1]
\eneq
and $L(c)(t)= 4tL(c)(1/4) + (1-4t){\rm diag}(\phi_2(c),0).$ Note $L$
maps $C$ into $M_2(M_{\phi_1, \phi_2}).$ Thus, since $K_i(C)$ is
finitely generated ($i=0,1$), as in Proposition 2.4 of
\cite{Lnhomp}, if $\ep$ is sufficiently small and ${\cal F}$ is
sufficiently large, there is $\gamma\in
Hom_{\Lambda}(\underline{K}(C), \underline{K}(M_{\phi_1,\phi_2}))$
such that
\beq\label{0AKK7}
[L]|_{\cal P}=\gamma|_{\cal P}
\eneq
for a given finite subset ${\cal P}\subset \underline{K}(C).$ Since $K_i(C)$ is finitely generated,
we may assume that $[L]=\gamma.$ It should be noted that,
with sufficiently small $\ep$ and sufficiently large ${\cal F},$ 
\beq\label{0AKK7+}
[\pi_0]\circ \gamma=[{\rm id}_C].
\eneq

Using section 2 and section 3, one also computes that
\beq\label{0AKK8}
\int_0^1\tau({dZ_j(t)\over{dt}}Z_j(t))dt=\tau(h_j)\tforal \tau\in
T(A).
\eneq
Therefore, if $R_{\phi_1,\phi_2}\circ \gamma(K_1(C))=0,$ then
\beq\label{0AKK8+}
\tau(h_j)=0\tforal \tau\in T(A).
\eneq

 On the other hand, for any given
$\eta>0$ and a finite set $\{z_1,z_2,...,z_n\}$ of generators of
$K_1(C),$ by (\ref{0AKK3}),
\beq\label{0AKK9}
\tau(h_j)<\eta
\eneq
for all $\tau\in T(A),$ provided that $\ep$ is sufficiently small
and ${\cal F}$ is sufficiently large.

Thus we have the following:

}

\end{NN}

\begin{lem}\label{AKK}
Let $C$ be a unital separable amenable \CA\, in ${\cal N}$ with
finitely generated $K_i(C)$ ($i=0,1$). Let $A$ be a unital
separable \CA\, and let $\phi_1, \phi_2: C\to A$ be two unital \hm
s. Suppose $\{z_1,z_2,...,z_m\}\subset  U(M_k(C))$ which forms  a
set of generators of $K_1(C)$ and $\eta>0$ is a positive number.

Then, there is $\ep>0$ and a finite subset ${\cal F}\subset C$
satisfying the following: Suppose that $A$ is  a unital separable
\CA\, and suppose that $\phi_1, \phi_2: C\to A$ are  two unital
\hm s such that
\beq\label{AKK1}
{\rm ad}\, v\circ \phi_1\approx_{\ep}\phi_2\,\,\,\text{on}\,\,\,
{\cal F}.
\eneq
Then it defines a \hm\, $\gamma\in Hom_{\Lambda}(\underline{K}(C),
\underline{K}(M_{\phi_1, \phi_2}))$ as
 in
\rm{(\ref{0AKK})}. Moreover,
\beq\label{AKK2}
\tau(\log(\phi_2(z_j)^*V^*\phi_1(z_j)V))<\eta,\,\,\,j=1,2,...,m
\eneq
for all $\tau\in T(A),$ where $V={\rm diag}(\overbrace{v,v,...,v}^k).$
\end{lem}

\begin{NN}\label{theta}

{\rm Let $C$ and $A$ be as in \ref{0AKK}, and $\phi_1: C\to A$ be
a unital \hm. Suppose that $v\in U(A)$ and
\beq\label{theta1}
\|v^*\phi_1(z_j)v-\phi_1(z_j)\|<\ep<1/2
\eneq
for some $z_1, z_2,...,z_n\in U(M_k(C)).$ Define $W_j(t)\in
U(M_2(C([0, 1], M_k(A))$ as follows
\beq\label{theta2}
W_j(t)={\rm diag}(\phi_1(z_j),1)(T_tVT_t^{-1})^*{\rm
diag}(\phi_1^*(z_j), 1)T_tVT_t^{-1},
\eneq
where
$$
V={\rm diag}(\overbrace{v,v,...,v}^n, 1_{M_k})\andeqn T_t=\begin{pmatrix}\cos ({\pi t}) & -\sin ({\pi t})\\
                     \sin({\pi t}) & \cos ({\pi t})
                     \end{pmatrix}
                     $$
 Put $b_j= {\rm diag}({1\over{2\pi i}}\log (\phi_1(z_j)^*v^*\phi_1(z_j)v, 1),\,\,\,j=1,2,...,n.$

Note that $W_j(0)={\rm diag}(\phi_1(z_j)v^*\phi_1(z_j)v, 1)$ and
$W_j(1)=1_{M_{2k}}.$ We obtain a continuous path of unitary
$Z_j(t)$ such that $Z_j(t)=\exp(i2\pi 4tb_j)$ for $t\in [0,1/4]$
and $Z(t)=W({4t-1\over{3}}).$  Thus $Z\in M_{2k}({\tilde SA}).$
Suppose that
$$
\|[\phi_1(a), \, v]\|<\dt\rforal {\cal F}
$$
for a sufficiently large finite subset ${\cal F}$ of $C$ and a
sufficient small $\dt.$ We may assume that  $\text{Bott}(\phi_1,
v)$ is well defined. Exactly as in \ref{0AKK} and \ref{AKK}, we
obtain a \hm\, in $Hom_{\Lambda}(\underline{K}(C),
\underline{K}(SA)).$ It will be denoted by
$\Gamma(\text{Bott}(\phi_1,v)).$

Suppose that $\Theta: \underline{K}(C)\to \underline{K}(SA)$ is
another \hm. Suppose there is $V\in C([0, 1-d], A)$ (for some
$0<d<1/4$) has the following property: Suppose that $z\in
U(M_k({\tilde C\otimes C_n}))$ (in $U(M_k(C))$, or in
$U(M_k({\tilde S(C\otimes C_n)})))$), where $C_n$ is a separable
commutative \CA\, with $K_0(C_n)=\Z/k\Z$ and $K_1(C_n)=\{0\},$
such that $\Theta([z])$ can be represented by a unitary
$\theta(z)\in M_k({\tilde S(A\otimes C_n)})$ (in $M_k(C),$ or  in
$M_k({\tilde S(C\otimes C_n)}))$) such that
\beq\label{theta2+1}
\theta(z)(t)=\phi_1(z)V^*(t)\phi_1(z)^*V(t)
\eneq
for $t\in [0, 1-d]$ and, for $t\in [1-d,1],$
\beq\label{theta2+2}
\|\theta(z)(t)-\phi_1(z)V(1-d)^*\phi_1(z)^*V(1-d)\|<\dt<1/2,
\eneq
where $\phi_1=\phi_1\otimes {\rm id}.$  By \ref{zeroV}, if
$\text{Bott}(\phi_1, u)([z])=\text{bott}_1(V,\phi_1(z)),$ then
$\Theta([z])=0$ (with sufficiently small $\dt$). Conversely, if
$\Theta([z])=0,$ by \ref{path}, $\text{Bott}(\phi_1, u)([z])=0.$
It follows that $\Gamma(\text{Bott}(\phi_1, u))$ is  defined
independent of the choice of the path $V(t)$ and
$\Gamma(\text{Bott}(\phi_1,u))=0$ if and only if
$\text{Bott}(\phi_1, u)=0.$

Using section 3 and section 2, one computes that
\beq\label{theta3}
\int_0^1 \tau({dZ_j(t)\over{dt}}Z_j(t))=\tau(b_j)\tforal \tau\in
T(A).
\eneq
If $\phi_1$ and $\phi_2$ are as in \ref{0AKK}, then
\beq\label{theta4}
R_{\phi_1,\phi_2}(Z_j(t))(\tau)=\tau(b_j)\rforal \tau\in T(A).
\eneq
In particular, if $\rho_A(\text{Bott}(\phi, v))=0,$ then
\beq\label{theta5}
\Gamma(\text{Bott}(\phi,v))(K_1(C))\subset {\rm ker}\,R_{\phi_1,
\phi_2}.
\eneq

 }
\end{NN}

\begin{thm}\label{TM}
Let $C$ be a unital AH-algebra and let $A$ be a unital separable
simple \CA\, of tracial rank zero. Suppose that $\phi_1, \phi_2:
C\to A$ are two unital monomorphisms. Then there exists a
continuous path of unitaries $\{u(t):t\in [0,\infty)\}\subset A$
such that
\beq\label{TM-2}
\lim_{t\to\infty}{\rm ad}\, u(t)\circ \phi_1(c)=\phi_2(c)\tforal
c\in C
\eneq
if and only if
\beq\label{TM-1}
[\phi_1]&=&[\phi_2]\,\,\,\text{in}\,\,\,KK(C,A),\\
\tau\circ \phi_1&=&\tau\circ \phi_2\tforal \tau\in T(A)\andeqn
{\tilde\eta}_{\phi_1,\phi_2}=\{0\}.
\eneq
\end{thm}

\begin{proof}
The ``only if " part follows from \ref{NecC}. We need to show that
``if" part of the theorem.

 Let $C=\lim_{n\to\infty}(C_n, \psi_n)$
be as in \ref{IND}. Let $\{{\cal F}_n\}$ be an increasing sequence
of finite subsets of $C$ such that $\cup_{n=1}^{\infty}{\cal F}_n$
is dense in $C.$

Define
$$
M_{\phi_1, \phi_2}=\{f\in C([0,1],A): f(0)=\phi_1(a)\andeqn
f(1)=\phi_2(a)\,\,\,{\rm for\,\,\,some}\,\,\,a \in C\}.
$$
Denote by $\pi_t: M_{\phi_1,\phi_2}\to A$ the point-evaluation at
$t$ ($t\in [0,1]$).
 Note also $\pi_0$ gives a
surjective \hm\, from $M_{\phi_1, \phi_2}$ onto $C$ and ${\rm
ker}\, \pi_0=SA.$

Since $C$ satisfies the Universal Coefficient Theorem, the
assumption that   $[\phi_1]=[\phi_2]$ in $KK(C,A)$ implies the
following
 exact sequence splits:
\beq\label{nFM}
\begin{array}{ccccccc}
0 \to  & \underline{K}(SA)  &\to &
\underline{K}(M_{\phi_1, \phi_2})
&{\stackrel{\pi_0}{\rightleftarrows}}_{\theta}& \underline{K}(C) &\to 0.
\end{array}
\eneq
Furthermore, since $\tau\circ \phi_1=\tau\circ \phi_2$ for all
$\tau\in T(A)$ and ${\tilde \eta}_{\phi_1,\phi_2}=0,$ we may also
assume that
\beq\label{FM-1-1}
{\rm R}_{\phi_1,\phi_2}(\theta(x))=0\rforal x\in K_1(C).
\eneq

By \cite{DL}, one has
\beq\label{FM-1}
\lim_{n\to\infty}(\underline{K}(C_n), [\psi_n])=\underline{K}(C).
\eneq

Since each $K_i(C_n)$ is finitely generated, there is an integer
$K(n)\ge 1$ such that
\beq\label{NNFM-00}
Hom_{\Lambda}(F_{K(n)}\underline{K}(C_n),
F_{K(n)}\underline{K}(A))=Hom_{\Lambda}(\underline{K}(C_n),
\underline{K}(A)).
\eneq

Let $\dt_n'>0,$ ${\cal G}_n'\subset C$  and ${\cal P}_n'\subset
\underline{K}(C)$  be finite subsets corresponding to $1/2^{n+2}$
and ${\cal F}_n$ required by \ref{LNHOMP}. Without loss of
generality, we may assume that ${\cal G}_n'\subset \psi_{n,
\infty}({\cal G}_{n})$ and ${\cal P}_n'= [\psi_{n, \infty}]({\cal
P}_n)$ for some finite subset ${\cal G}_n$ of $C_n$ and for some
finite subset ${\cal P}_n\subset \underline{K}(C_{n}).$  We may
assume that ${\cal P}_n$ contains a set of generators of
$F_{K(n)}\underline{K}(C_{n}),$ $\dt_n'<1/2^{n+3}$ and ${\cal
F}_n\subset {\cal G}_n'.$ We also assume that $\text{Bott}(h',
u')|_{{\cal P}_n}$ is well defined whenever $\|h'(a),\,
u']\|<\dt_n'$ for all $a\in {\cal G}_n'$ and for any unital \hm\,
$h'$ and unitary $u'.$ Note that $\text{Bott}(h',u')|_{{\cal
P}_n}$ defines $\text{Bott}(h',u').$

We further assume that
\beq\label{NNF0}
\text{Bott}(h,u)|_{{\cal P}_n}=\text{Bott}(h',u)|_{{\cal P}_n}
\eneq
provided that $h\approx_{\dt_n'}h'$ on ${\cal G}_n'.$


We may also assume that $\dt_n'$ is smaller than $\dt/3$ for that
$\dt$ in \ref{ddbot} for $C_n$ and ${\cal P}_n.$

Let $k(n)\ge n$ and $\eta_n>0$ ( in place of $\dt$)  be required by \ref{CMhit} for
$\psi_{n,\infty}({\cal G}_n)$( in place of ${\cal F}$), ${\cal P}_n$ (in place of ${\cal P}$) and $\dt_n'/4$ (in place of $\ep$).  For
$C_n,$ since $K_i(C_n)$ ($i=0,1$) is finitely generated, by
choosing larger $k(n),$ we may assume that
$(\phi_{k(n),\infty})_{*i}$ is injective on $(\phi_{n,
k(n)})_{*i}(K_i(C_n)),$ $i=0,1.$
Since $K_i(C_n)$ is finitely generated, by (\ref{NNFM-00}), we may
further assume that $[\phi_{k(n), \infty}]$ is injective on
$[\phi_{n, k(n)}](\underline{K}(C_n)),$ $n=1,2,....$

By passing to a subsequence, to simplify notation, we may assume that $k(n)=n+1.$

Let $\dt_n=\min\{\eta_n, \dt_n/2'\}.$

By \ref{LTRBOT},  there are unitaries $v_n\in U(A)$ such that
\beq\label{NFM2}
&&{\rm ad}\, v_n\circ \phi_1\approx_{\dt_{n+1}/4}
\phi_2\,\,\,\text{on}\,\,\,\psi_{n+1,\infty}({\cal
G}_{n+1}),\\\label{NFM3}
&&\rho_A(\text{bott}_1(\phi_2,\, v_n^*v_{n+1}))(x)=0
\rforal x\in \psi_{n+1,\infty}(K_1(C_{n+1}))\andeqn\\
&&\|[\phi_2(a),\,v_n^*v_{n+1}]\|<\dt_{n+1}/2\rforal a\in
\psi_{n+1,\infty}({\cal G}_{n+1})
\eneq
($K_1(C_{n+1})$ is finitely generated).

Note that, by (\ref{NNF0}), we may also assume that
\beq\label{NFM3+}
\text{Bott}(\phi_1, v_{n+1}v_n^*)|_{{\cal
P}_n}=\text{Bott}(v_n^*\phi_1v_n,v_n^*v_{n+1})|_{{\cal
P}_n}=\text{Bott}(\phi_2, v_n^*v_{n+1})|_{{\cal P}_n}.
\eneq
In particular,
\beq\label{NFM3+1}
\text{bott}_1(v_n^*\phi_1v_n,v_n^*v_{n+1})(x)=\text{bott}_1(\phi_2,
v_n^*v_{n+1})(x)
\eneq
for all $x\in \psi_{n+1,\infty}(K_1(C_{n+1})).$

 By applying \ref{AKK}, without loss of generality, we may assume that
$\phi_1$ and $v_n$ define $\gamma_n\in
Hom_{\Lambda}(\underline{K}(C_{n+1}), \underline{K}(M_{\phi_1,
\phi_2})).$ Note that $[\pi_0]\circ \gamma_n=[{\rm id}_{C_{n+1}}].$  Furthermore, by \ref{AKK}, without loss of
generality, we may assume that
\beq\label{0TM-10-1}
\tau(\log(\phi_2\circ\psi_{n+1,\infty}(z_j)^*V_n^*\phi_1\circ \psi_{n+1, \infty}(z_j)V_n))<\dt_{n+1},\,\,\,j=1,2,...,r(n),
\eneq
where $\{z_1,z_2,...,z_{r(n)}\}\subset U(M_k(C_{n+1}))$ which forms a
set of generators of $K_1(C_{n+1})$ and where $V_n={\rm
diag}(\overbrace{v_n,v_n,...,v_n}^k).$

Let $H_n=[\phi_{n+1}](\underline{K}(C_{n+1})).$ Since
$\cup_{n=1}[\phi_{n+1,
\infty}](\underline{K}(C_n))=\underline{K}(C)$ and $[\pi_0]\circ \gamma_n=[{\rm id}_{C_{n+1}}],$ we conclude that
\beq\label{0TM-10-}
\underline{K}(M_{\phi_1,\phi_2})=\underline{K}(SA)+\cup_{n=1}^{\infty}\gamma_n(H_n).
\eneq
Thus, by passing to a subsequence, we may further assume that
\beq\label{0TM-10--}
\gamma_{n+1}(H_n)\subset
\underline{K}(SA)+\gamma_{n+2}(H_{n+1}),\,\,\,n=1,2,....
\eneq

By identifying $H_n$ with $\gamma_{n+1}(H_n),$ we may write $j_n:
\underline{K}(SA)\oplus H_n\to \underline{K}(SA)\oplus H_{n+1}.$
By (\ref{0TM-10-}), the inductive limit is $\underline{K}(M_{\phi_1,\phi_2}).$

From the definition of $\gamma_n,$ we note that,
$\gamma_{n}-\gamma_{n+1}\circ [\psi_{n+1}]$ maps
$\underline{K}(C_{n+1})$ into $\underline{K}(SA).$

By \ref{theta}
$$
\Gamma(\text{Bott}(\phi_1\circ \psi_{n+2,\infty},
v_{n}v_{n+1}^*))|_{H_{n}}=(\gamma_{n+1}-\gamma_{n+2}\circ
[\psi_{n+2}])|_{H_n}
$$
gives a \hm\, $\xi_n: H_n\to \underline{K}(SA).$
Put $\zeta_n=\gamma_{n+1}|_{H_n}.$ Then
\beq\label{0TM-19}
j_n(x,y)=(x+\xi_n(y),[\psi_{n+2}](y))\rforal (x, y)\in
\underline{K}(SA)\oplus H_n.
\eneq

Thus, we obtain the following  diagram:
\beq\label{LDG1}
\begin{array}{ccccccc}
0 \to  & \underline{K}(SA)  &\to & \underline{K}(SA)\oplus H_n
&\to
& H_n &\to 0\\\nonumber
 &\| & &\hspace{0.4in}\| \hspace{0.15in}\swarrow_{\xi_n} \hspace{0.05in}\downarrow_{[\psi_{n+2}]} &&
 \hspace{0.2in}\downarrow_{[\psi_{n+2}]} &\\\label{LDG2}
 0 \to  & \underline{K}(SA)  &\to & \underline{K}(SA)\oplus H_{n+1} &\to & H_{n+1} &\to 0\\\nonumber
  &\| & &\hspace{0.4in}\| \hspace{0.1in}\swarrow_{\xi_{n+1}}\downarrow_{[\psi_{n+3}]} &&
\hspace{0.2in} \downarrow_{[\psi_{n+3}]} &\\
 0 \to  & \underline{K}(SA)  &\to & \underline{K}(SA)\oplus H_{n+2} &\to & H_{n+2} &\to 0\\
 \end{array}
\eneq
By the assumption that ${\tilde \eta}_{\phi_1,\phi_2}=0,$ $\theta$
also gives the following
$$
{\rm ker}R_{\phi_1, \phi_2}={\rm ker}\rho_A\oplus K_1(C).
$$

Define $\theta_n=\theta\circ [\psi_{n+2, \infty}]$ and
$\kappa_n=\zeta_n-\theta_n.$ Note that
\beq\label{NFM11-1}
\theta_n=\theta_{n+1}\circ [\psi_{n+2}].
\eneq
We also have that
\beq\label{NFM-11-2}
\zeta_n-\zeta_{n+1}\circ [\psi_{n+2}]=\xi_n.
\eneq

 Since $[\pi_0]\circ (\zeta_n-\theta_n)=0,$ $\kappa_n$ maps $H_n$
into $\underline{K}(SA).$ It follows that
\beq\nonumber
\kappa_n-\kappa_{n+1}\circ [\psi_{n+2}] &=&
\zeta_n-\theta_n-\zeta_{n+1}\circ [ \psi_{n+2}]+\theta_{n+1}\circ
[\psi_{n+2}]\\\label{NFM11}
&=&\zeta_n-\zeta_{n+1}\circ[\psi_{n+2}]=\xi_n
\eneq

It follows from \ref{VuV} that  there are integers $N_1\ge 1,$  a
${\dt_{n+1}\over{4}}$-$\psi_{n+1}({\cal G}_{n+1})$-multiplicative
\morp\, $L_n: \psi_{n+1, \infty}(C_{n+1})\to
M_{1+N_1}(M_{\phi_1,\phi_2}),$ a unital \hm\, $h_0: \psi_{n+1,
\infty}(C_{n+1})\to M_{N_1}(\C),$
and a continuous path of unitaries $\{V_n(t): t\in [0,3/4]\}$ of
$M_{1+N_1}(A)$ such that $[L_n]|_{{\cal
P}_{n+1}'}$ is well defined, $V_n(0)=1_{M_{1+N_1}(A)},$
\beq\label{NVuV1}
[L_n\circ \psi_{n+1,\infty}]|_{{\cal P}_n}=(\theta\circ
[\psi_{n+1,\infty}]+[h_0\circ \psi_{n+1,\infty}])|_{{\cal P}_n},
\eneq
\beq\label{NVuV2}
\hspace{-0.4in}\pi_t\circ L_n\circ
\psi_{n+1,\infty}\approx_{\dt_{n+1}/4} {\rm ad}\, V_n(t)\circ
((\phi_1\circ \psi_{n+1,\infty})\oplus (h_0\circ
\psi_{n+1,\infty}))
\eneq
on $\psi_{n+1,\infty}({\cal G}_{n+1})$ for all $t\in (0,3/4],$
\beq\label{NVuV3}
\hspace{-0.4in}\pi_t\circ L_n\circ
\psi_{n+1,\infty}\approx_{\dt_{n+1}/4} {\rm ad}\, V_n(3/4)\circ
((\phi_1\circ \psi_{n+1,\infty})\oplus (h_0\circ
\psi_{n+1,\infty}))
\eneq
on $\psi_{n+1,\infty}({\cal G}_{n+1})$ for all $t\in (3/4,1),$ and
\beq\label{NVuV33}
\pi_1\circ L_n\circ
\psi_{n+1,\infty}\approx_{\dt_{n+1}/4}\phi_2\circ
\psi_{n+1,\infty}\oplus h_0\circ \psi_{n+1,\infty}
\eneq
on $\psi_{n+1,\infty}({\cal G}_{n+1}),$ where $\pi_t:
M_{\phi_1,\phi_2}\to A$ is the point-evaluation at $t\in (0,1).$

Note that $R_{\phi_1,\phi_2}(\theta(z))=0$ for all $x\in
\phi_{n+1,\infty}(K_1(C_{n+1})).$  As computed in \ref{0AKK},
\beq\label{NVuV4}
\tau(\log((\phi_2(z)\oplus h_0(z)^*V_n(3/4)^*(\phi_1 (z)\oplus
h_0(z))V_n(3/4)))=0
\eneq
for $z=\psi_{n+1, \infty}(y),$ where $y$ is in a set of generators
of $K_1(C_{n+1})$ and for all $\tau\in T(A).$

Define $W_n'={\rm daig}(v_n,1)\in M_{1+N_1}(A).$ Then
$\text{Bott}((\phi_1\oplus h_0)\circ \psi_{n+1, \infty},\,
W_n'(V_n(3/4)^*)$ defines a \hm\, ${\tilde \kappa}_n\in
Hom_{\Lambda}(\underline{K}(C_{n+1}),\underline{K}(A)).$ By
(\ref{0TM-10-1})
\beq\label{NNFM1-}
\tau(\log((\phi_2\oplus h_0)\circ \psi_{n+1,\infty}(z_j)^*{\tilde V}_n^*(\phi_1\oplus h_0)\circ \psi_{n+1,\infty}(z_j){\tilde V}_n))<\dt_{n+1},
\eneq
$j=1,2,...,r(n),$ where ${\tilde V}_n={\rm diag}(V_n,1).$
Then, by  (\ref{NVuV4}), exactly as in the argument from
(\ref{ltrbot-17})--(\ref{ltrbot-19}), we compute that
\beq\label{NNFM1}
\rho_A({\tilde \kappa}_n(z_j))(\tau) <\dt_{n+1},\,\,\,j=1,2,....
\eneq
It follows from \ref{CMhit} that there is a unitary $w_n'\in U(A)$
such that
\beq\label{NNFM2}
\|[\phi_1(a), w_n']\|<\dt_{n+1}'/4\rforal a\in \psi_{n,\infty}({\cal G}_n)\andeqn\\
\text{Bott}(\phi_1\circ \psi_{n,\infty},\,
w_n')|_{K_1(C_n)}=-{\tilde \kappa}_n|_{[\psi_{n}](K_1(C_{n}))}.
\eneq
By (\ref{NNF0}),
\beq\label{NNFM3}
\text{Bott}(\phi_2\circ \psi_{n,\infty},\,v_n^*w_n'v_n)|_{{\cal
P}_n}=-{\tilde \kappa}_n|_{[\psi_n](\underline{K}(C_n))}.
\eneq
Put $w_n=v_n^*w_n'v_n.$

It follows from \ref{theta} that
\beq\label{NNFM4}
\Gamma(\text{Bott}(\phi_1\circ \psi_{n, \infty}, w_n'))=-\kappa_n
\andeqn \Gamma(\text{Bott}(\phi_1\circ \psi_{n+1, \infty},
w_{n+1}'))=-\kappa_{n+1}.
\eneq
We also have
\beq\label{NNFM5}
\Gamma(\text{Bott}(\phi_1\circ \psi_{n,\infty},
v_nv_{n+1}^*))|_{H_n}=\zeta_n-\zeta_{n+1}\circ [\psi_{n+2}]=\xi_n.
\eneq
But, by (\ref{NFM11}),
\beq\label{NNFM6}
-\kappa_n +\xi_n+\kappa_{n+1}\circ [\psi_{n+2}]=0.
\eneq
It follows from (\ref{NNFM4}), (\ref{NNFM5}), (\ref{NNFM6}) and
\ref{theta} that
\beq\label{NNFM7}
\hspace{-0.2in} -\text{Bott}(\phi_1\circ \psi_{n, \infty},\,w_n')
+\text{Bott}(\phi_1\circ \psi_{n, \infty}, \,v_nv_{n+1}^*)
+\text{Bott}(\phi_1\circ\psi_{n,\infty}, w_{n+1}') =0.
\eneq
Define $u_n=v_nw_n^*,$ $n=1,2,....$ Then, by (\ref{NFM2}) and
(\ref{NNFM2}),
\beq\label{NFM16}
{\rm ad}\, u_n\circ \phi_1\approx_{\dt_n'/2} \phi_2\rforal a\in
\psi_{n, \infty}({\cal G}_n).
\eneq

From (\ref{NFM3+}), (\ref{NNF0}) and (\ref{NNFM7}), we compute that
\beq\label{NNFM8}
&&\hspace{-0.6in}\text{Bott}(\phi_2\circ \psi_{n, \infty},u_n^*u_{n+1})\\
&=& \text{Bott}(\phi_2\circ \psi_{n, \infty}, w_nv_n^*v_{n+1}w_{n+1}^*)\\
&=& \text{Bott}(\phi_2\circ \psi_{n, \infty}, w_n)+\text{Bott}(\phi_2\circ \psi_{n, \infty}, v_n^*v_{n+1})\\
&&\hspace{1.8in}+\text{Bott}(\phi_2\circ \psi_{n, \infty}, w_{n+1}^*)\\
&=&\text{Bott}(\phi_1\circ \psi_{n, \infty},w_n')+\text{Bott}(\phi_1\circ \psi_{n, \infty}, v_{n+1}v_n^*)\\
&&\hspace{1.8in}+\text{Bott}(\phi_1\circ \psi_{n, \infty}, (w_{n+1}')^*)\\
&=&-[-\text{Bott}(\phi_1\circ \psi_{n, \infty}, w_n')+\text{Bott}(\phi_1\circ \psi_{n, \infty}, v_nv_{n+1}^*)\\
&&\hspace{1.8in}+\text{Bott}(\phi_1\circ \psi_{n, \infty}, w_{n+1}')]\\
&=&0
\eneq

Therefore, by \ref{LNHOMP}, there exists a continuous path of
unitaries $\{z_n(t): t\in [0,1]\}$ of $A$ such that
\beq\label{FM-16}
&&z_n(0)=1,\,\,\, z_n(1)=u_n^*u_{n+1}\andeqn\\\label{FM-16+}
&&\|[\phi_2(a),\, z_n(t)]\|<1/2^{n+2}\rforal a\in {\cal F}_n\andeqn t\in [0,1].
\eneq
Define
$$
u(t+n-1)=u_nz_{n+1}(t)\,\,\,t\in (0,1].
$$
Note that $u(n)=u_{n+1}$ for all integer $n$ and $\{u(t):t\in [0,
\infty)\}$ is a continuous path of unitaries in $A.$ One estimates
that, by (\ref{NFM16}) and (\ref{FM-16+}),
\beq\label{FM-17}
{\rm ad}\, u(t+n-1)\circ \phi_1
&\approx_{\dt_n'}& {\rm
ad}\,z_{n+1}(t)\circ \phi_2
\\
 &\approx_{1/2^{n+2}}& \phi_2
\,\,\,\,\,\,\,\,\,\,\,\,\,\,\,\,\,\,\hspace{0.5in}\text{on}\,\,\,\,{\cal
F}_n
\eneq
 for all $t\in (0,1).$

It then follows that
\beq\label{FM-18}
\lim_{t\to\infty}u^*(t)\phi_1(a) u(t)=\phi_2(a)\rforal a\in C.
\eneq

\end{proof}

\section{Applications}

\begin{thm}\label{AT1}
Let $C$ be a unital AH-algebra with real rank zero and let $B$ be
a unital separable  \CA\, with tracial rank zero. Suppose that
$\phi_1,\phi_2: C\to B$ are two unital monomorphisms. Then there
exists a continuous path of unitaries $\{u(t): t\in [0, \infty)\}$
such that
$$
\lim_{t\to\infty} {\rm ad}\, u(t)\circ \phi_1(a)=\phi_2(a)\tforal
a\in C
$$
if and only if
\beq\label{AT1-1}
[\phi_1]=[\phi_2]\,\,\,\text{in}\,\,\, KK(C,B)\andeqn
{\tilde{\eta}}_{\phi_1,\phi_2}=0.
\eneq

\end{thm}

\begin{proof}
To prove the ``if" part of the theorem, it suffices to show that
\beq\label{AT1-1-1}
\tau\circ \phi_1(a)=\tau\circ \phi_2(a)\rforal a\in C\andeqn
\rforal \tau\in T(A).
\eneq
It then suffices to show that (\ref{AT1-1-1}) holds for all $a\in
C_{s.a.}.$ The assumption that $[\phi_1]=[\phi_2]$ in $KK(C,A)$
implies that (\ref{AT1-1-1}) holds for the case that $a$ is a
projection, whence for any self-adjoint elements with finite
spectrum. Since $C$ is assumed to have real rank zero, one
concludes that (\ref{AT1-1-1}) holds for all self-adjoint
elements.

\end{proof}

\begin{cor}\label{AC1}
Let $A$ be a unital separable simple \CA\, with tracial rank zero
and satisfying the UCT. Suppose that $\af: A\to A$ is a unital
endomorphisms. Then $\af$ is asymptotically inner if and only if
\beq\label{AC1-1}
[\af]=[{\rm id}]\,\,\, \text{in}\,\,\, KK(A,A)\andeqn
{\tilde\eta}_{{\rm id}_A,\af}=0.
\eneq
\end{cor}

\vspace{0.2in}

Now we consider the problem when $C\rtimes_{\af}\Z$ can be
embedded into a unital simple AF-algebra, where $C$ is a unital
AH-algebra and $\af$ is an automorphism (see \cite{V1}, \cite{V2}, \cite{V3},  \cite{Bn1}, \cite{Bn2},
\cite{M}, \cite{Pi} and \cite{Lnemb2} for the  background).

\vspace{0.1in}
\begin{lem}\label{Lat2}
Let $A$ and $B$ be two unital separable \CA s and let $\phi_1,
\phi_2: A\to B$ be two unital monomorphisms. Suppose that  $A$
satisfies the Universal Coefficient Theorem and
$$
[\phi_1]=[\phi_2]\,\,\,\text{in}\,\,\, KK(A,B)\andeqn \tau\circ
\phi_1=\tau\circ \phi_2
$$
for all $\tau\in T(A).$ Suppose that $K_1(B)=0$ and
$K_0(B)=\rho_B(K_0(B))$ is torsion free and divisible. Suppose
also that
$$
\rho_A(K_0(A))=R_{\phi_1,\phi_2}(K_1(M_{\phi_1, \phi_2})).
$$
Then, there is $\theta\in Hom_{\Lambda}(\underline{K}(A),
\underline{K}(B))$ such that the following exact sequence splits:
$$
0\to \underline{K}(B)\stackrel{[\imath]}{\to}
\underline{K}(M_{\phi,
\psi})\stackrel{[\pi_0]}{\rightleftarrows}_{\theta}
\underline{K}(A)\to 0
$$
and
$$
R_{\phi_1,\phi_2}\circ \theta=0.
$$
Consequently ${\tilde\eta}_{\phi_1,\phi_2}=0.$

\end{lem}

\begin{proof}
From the assumption there is $\theta\in \in
Hom_{\Lambda}(\underline{K}(A), \underline{K}(B))$ such that the
following exact sequence splits:
\beq\label{Lat-1}
0\to \underline{K}(SB)\stackrel{[\imath]}{\to}
\underline{K}(M_{\phi_1,
\phi_2})\stackrel{[\pi_0]}{\rightleftarrows}_{\theta}
\underline{K}(A)\to 0.
\eneq
Homomorphism $\theta|_{K_1(A)}$ gives the following splitting
exact sequence:
$$
0\to K_0(B)\stackrel{\imath_*}{\to} K_1(M_{\phi_1,
\phi_2})\stackrel{(\pi_0)_{*1}}{\rightleftarrows}_{\theta|_{K_1(A)}}
K_1(A)\to 0.
$$
Using $\theta,$ we may write
\beq\label{Lat-2}
K_1(M_{\phi_1,\phi_2})=K_0(B)\oplus K_1(A).
\eneq
Since
$K_0(B)=\rho_B(K_0(B))=R_{\phi_1,\phi_2}(K_1(M_{\phi_1,\phi_2})),$
there is an isomorphism $\gamma: K_1(A)\to {\rm
ker}R_{\phi_1,\phi_2}$ such that $(\pi_0)_{*1}\circ \gamma={\rm
id}_{K_1(A)}.$ Therefore
\beq\label{Lat-3}
(\theta|_{K_1(A)}-\gamma)(x)\in K_0(B)\tforal x\in K_1(A).
\eneq
 However, since $K_0(B)$ is torsion free, if $x\in Tor(K_1(A))$
 then
 \beq\label{Lat-4}
(\theta|_{K_1(A)}-\gamma)(x)=0.
\eneq
Let $\rho_k^{1}$ be the map from $K_0(B)\to K_0(B, \Z/k\Z).$ Then
since $K_0(B)$ is divisible,
\beq\label{Lat-5-}
\rho_k^{1}=0.
\eneq
In particular,
\beq\label{Lat-5}
\rho_k^{1}\circ (\theta|_{K_1(A)}-\gamma)=0.
\eneq
Define $\theta'$ as follows $\theta'|_{K_1(A)}=\gamma,$
$\theta'|_{K_0(A)}=\theta|_{K_0(A)}$ and
$\theta'|_{K_i(A,\Z/k\Z)}= \theta|_{K_i(A,\Z/k\Z)}$ for $i=0,1$
and $k=2,3,....$

Now put $C=M_{\phi_1, \phi_2}.$
 Since $\theta|_{K_i(A,\Z/k\Z)}=\theta'|_{K_i(A,\Z/k\Z)},$
 to show that $\theta'\in
 Hom_{\Lambda}(\underline{A},\underline{C}),$ it suffices to show
 the following diagram commutes.
$$
{\small \put(-160,0){$K_0(A)$} \put(0,0){$K_0(A,{\Z }/k\Z )$}
\put(180,0){$K_1(A)$} \put(-85,-40){$K_0(C)$}
\put(0,-40){$K_0(C,{\Z}/k{\Z})$} \put(105,-40){$K_1(C)$} \put(-85,
-70){$K_0(C)$} \put(0,-70){$K_1(C, {\Z}/k{\Z})$}
\put(105,-70){$K_1(C)$} \put(-160,-110){$K_0(A)$}
\put(0,-110){$K_1(A,{\Z}/k{\Z})$} \put(180,-110){$K_1(A)$}
\put(-120, 2){\vector(1,0){95}} \put(70,1){\vector(1,0){95}}
\put(-123,-3){\vector(1,-1){30}} \put(30,-3){\vector(0,-1){25}}
\put(180,-2){\vector(-1,-1){30}} \put(-45,-38){\vector(1,0){35}}
\put(70,-38){\vector(1,0){25}} \put(-147, -90){\vector(0,1){85}}
\put(-75,-60){\vector(0,1){15}} \put(115, -45){\vector(0,-1){15}}
\put(190,-7){\vector(0,-1){85}} \put(-7,-68){\vector(-1,0){35}}
\put(95,-68){\vector(-1,0){25}} \put(-123,-102){\vector(1,1){30}}
\put(175, -105){\vector(-1,1){30}} \put(30,-104){\vector(0,1){30}}
\put(-5, -108){\vector(-1,0){100}}
\put(170,-108){\vector(-1,0){95}} \put(-112,-12){$\theta$}
\put(15, -15){$\theta $} \put(145, -12){$\gamma$}
\put(-122,-92){$\theta$} \put(15, -88){$\theta$} \put(160,
-88){$\gamma$} }
$$
 From (\ref{Lat-1}), since
$K_1(B)=0$ and $K_0(B)$ is torsion free and divisible, we may
write  that
\beq\label{Lat-6}
K_0(C, \Z/k\Z)&=&K_1(B,\Z/k\Z)\oplus
K_0(A,\Z/k\Z)=K_0(A,\Z/k\Z)\andeqn\\\label{Lat-7}
K_1(C,\Z/k\Z)&=&K_0(B,\Z/k\Z)\oplus K_1(A,\Z/k\Z)=K_1(A,\Z/k\Z).
\eneq
Thus, by (\ref{Lat-4}), (\ref{Lat-5}), (\ref{Lat-2}),
(\ref{Lat-6}) and (\ref{Lat-7}), one checks the above 12-term
diagram is commutative.

Thus, by replacing $\theta$ by $\theta',$ the lemma follows.
\end{proof}

We will need the following:

\begin{lem}{\rm (Corollary 5.5 of \cite{Lnemb2})}\label{ATL}
Let $A$ be a unital AH-algebra, let $\af\in Aut(A)$ and let $B\in
{\cal N}$ be a unital separable simple \CA\, with tracial rank
zero. Suppose that there is a unital monomorphism $h: A\to B$ such
that there exists a continuous path of unitaries $\{U(t): t\in [0,
\infty)\}$ of $B$ such that
\beq
\lim_{t\to\infty}{\rm ad}\, U(t)\circ h\circ\af(a)=h(a)\tforal
a\in A.
\eneq
Then $A\rtimes_{\af}\Z$ can be embedded into a unital simple
AF-algebra.
\end{lem}

\begin{thm}\label{AT2}
Let $A$ be a unital AH-algebra and let $\af\in Aut(A)$ be an
automorphism. Then $A\rtimes_{\af}\Z$ can be embedded into a
unital simple AF-algebra if and only if $A$ admits a faithful
$\af$-invariant tracial state $\tau.$
\end{thm}

\begin{proof}
Only ``if "  part needs a proof. Now let $t$ be an $\af$-invariant
faithful tracial state of $A.$ It follows from Theorem 4.1  of
\cite{Lnemb2} that there is a unital separable simple AF-algebra
$B_0$ with a unique tracial state $\tau$ with
$$
(K_0(B_0), K_0(B_0)_+, [1_{B_0}])=(\D, \D_+, 1),
$$
where $\D=\tau(K_0(B_0))$ a countable dense subgroup of $\R,$ and
a unital monomorphism $h: A\to B_0$ such that
$$
\tau\circ h=t.
$$

By embedding $B_0$ to $B_0\otimes Q$ for the UHF-algebra with
$K_0(Q)=\Q,$ we may further assume that  $\D$ is divisible.

Since $t(\af(a))=t(a)$ for all $a\in A,$ we have
$$
\tau\circ h(a)=\tau\circ h\circ\af(a)\tforal a\in A.
$$
This implies that
$$
h_{*0}=(h\circ \af)_{*0}.
$$
Since $K_1(B_0)=\{0\}$ and $K_0(B_0)$ is torsion free, we conclude
that
$$
[h]=[h\circ \af]\,\,\,\text{in}\,\,\, KL(A,B_0).
$$
However, since $K_0(B_0)$ is divisible (and $K_1(B_0)=\{0\}$), we
actually have
$$
[h]=[h\circ \af]\,\,\,\text{in}\,\,\, KK(A,B_0).
$$
Let
$$
M_{h,h\circ \af}=\{f\in C([0,1], B_0):  f(0)=h(a)\andeqn
f(1)=h\circ \af(a)\,\,\,\text{for\,\,\,some}\,\,\, a\in A\}.
$$
We have
\beq\label{AT-0}
K_1(M_{h,h\circ\af})=K_1(SB_0)\oplus K_1(A)=K_0(B_0)\oplus K_1(A).
\eneq

 Let $\D_1$ be a divisible countable dense subgroup of
$\R$ which contains $\D$ and $R_{h, h\circ \af}(K_1(M_{h,h\circ
\af})).$ There exists a unital separable simple AF-algebra $B$
with
$$
(K_0(B), K_0(B)_+, [1_B])=(\D_1, (\D_1)_+, 1).
$$
There is a unital  embedding $j: B_0\to B$ which induces the
natural embedding $j_*: \D\to \D_1.$ We continue to use $\tau$ for
the unique tracial state of $B.$

Let
$$
{\tilde M}=\{f\in C([0,1], B):f(0)=h(a)\andeqn f(1)=h\circ
\af(a)\,\,\,\text{for\,\,\,some}\,\,\, a\in A\}.
$$
The embedding $j$ induces a unital embedding ${\tilde j}:
M_{h,h\circ \af}\to {\tilde M}.$ We also have
\beq\label{AT-1}
[j\circ h]=[j\circ h\circ \af]\,\,\,\text{in}\,\,\, KK(A, B)
\andeqn \tau\circ j\circ h=\tau\circ j\circ\circ h\circ\af.
\eneq
Note that, we may write that
\beq\label{AT-2}
K_1({\tilde M})=K_0(B)\oplus K_1(A).
\eneq
Moreover, it is easy to see that ${\tilde j}_{*1}|_{K_0(B_0)}=
j_{*0}.$
 Let $z\in M_l({\tilde M}).$ Then there is $l'\ge l$ and
unitary $z'\in M_{l'}(M_{h, h\circ \af})$ such that
$$
\begin{pmatrix} z &0\\
                0 & 1\end{pmatrix} (z')^*
\in K_1(SB)=K_0(B)
$$
for some unitary $z'\in M_{l' }(M_{j\circ h,j\circ h\circ \af}).$
In particular,
\beq\label{AT-3}
{\rm im}R_{j\circ h, j\circ h\circ \af}\subset {\rm im}R_{h,
h\circ \af}+\D_1 \subset \D_1.
\eneq
It follows from \ref{Lat2} that
\beq\label{AT-4}
{\tilde \eta}_{j\circ h,j\circ h\circ \af}=0.
\eneq
By \ref{TM}, there exists a continuous path of unitaries $\{U(t):
t\in [0, \infty)\}$ of $B$ such that
\beq\label{AT-5}
\lim_{t\to\infty}{\rm ad}\, U(t)\circ j\circ h(a)=j\circ h\circ
\af(a)\tforal a\in A.
\eneq
It follows from \ref{ATL} that $A\times_\af\Z$ can be embedded
into a unital simple AF-algebra.

\end{proof}

\begin{cor}
Let $A$ be a unital simple AH-algebra and let $\af\in Aut(A)$ be
an automorphism. Then $A\rtimes_{\af}\Z$ can be embedded into a
unital simple AF-algebra.
\end{cor}

\begin{proof}
Since $A$ is simple, it is well-known that there exists at least
one $\af$-invariant tracial state which is faithful.

\end{proof}

\section{Strong asymptotic unitary equivalence}

Let $\phi_1, \phi_2: C\to A$ be two unital monomorphisms which are
asymptotically unitarily equivalent. A natural question is: can
one find a continuous path of unitaries $\{u_t: t\in [0,
\infty)\}$ of $B$ such that
$$
u_0=1_B\andeqn \lim_{t\to\infty}{\rm ad}\, u_t\circ
\phi_1(a)=\phi_2(a)\tforal a\in A?
$$
It certainly desirable to have an affirmative answer. Corollary
\ref{SUC} could be easily proved from the main theorem \ref{TM}
(without introducing $H_1(K_0(C), K_1(B)))$). However, in general,
unfortunately, the answer is negative. In this section we give a
detailed discussion of this phenomenon.

\begin{df}
{\rm Let $A$ and $B$ be two unital \CA.  Suppose that $\phi_1,
\phi_2: A\to B$ are two unital \hm s. We say that $\phi_1, \phi_2$
are strongly asymptotically unitarily equivalent if there exists a
continuous path of unitaries $\{u_t: t\in [0, \infty)\}$ of $B$
such that
$$
u_0=1_B\andeqn \lim_{t\to\infty}{\rm ad}\, u_t\circ
\phi_1(a)=\phi_2(a)\tforal a.
$$

If the answer to the above question is negative, then when they
$\phi_1$ and $\phi_2$ are strong asymptotically unitarily
equivalent?}

\end{df}

\begin{df}\label{Dsu}
Let $A$ be a unital \CA\, and $B$ be another \CA. Recall
(\cite{Lnequ}) that
$$
H_1(K_0(A), K_1(B))=\{x\in K_1(B): \phi([1_A])=x,\,\phi\in
Hom(K_0(A), K_1(B))\}.
$$
\end{df}

Exactly as in 3.3 of \cite{Lnequ}, we have the following:

\begin{prop}\label{Sup}
Let $A$  be a unital separable \CA\, and let $B$ be a unital \CA.
Suppose that $\phi: A\to B$ is a unital \hm\, and $u\in U(B)$ is a
unitary. Suppose that there is a continuous path of unitaries
$\{u(t): t\in [0,\infty)\}\subset B$ such that
\beq\label{sup-1}
u(0)=1_B\andeqn \lim_{t\to\infty}{\rm ad}\, u(t)\circ \phi(a)={\rm
ad}\, u\circ \phi(a)
\eneq
for all $a\in A.$  Then
$$
[u]\in H_1(K_0(A), K_1(B)).
$$

\end{prop}

\begin{proof}
In fact, one has
\beq\label{sup-2}
\lim_{t\to\infty}uu(t)^*\phi(a)u(t)u^*=\phi(a)\rforal a\in C
\eneq
for all $t\in [0, \infty).$ Define $\psi_t: A\otimes C(S^1)\to B$
by
$$
\psi_t(a\otimes f)=\phi(a)f(u(t)u^*)\rforal a\in A\andeqn f\in
C(S^1).
$$
Note that
$$
\lim_{t\to\infty}\|[\phi(a), \,u(t)u^*]\|=0\rforal a\in A.
$$
Therefore
$$
\lim_{t\to\infty}\|\psi_t(bc)-\psi_t(b)\psi_t(c)\|=0
$$
for all $b,c\in A\otimes C(S^1).$ Thus $[\{\psi_t\}]\in
KK(C\otimes C(S^1),A).$ We may view that $[\{\psi_t\}]\in
Hom_{\Lambda}(\underline{K}(A\otimes C(S^1)), \underline{K}(B)).$
Let $\kappa: K_0(A)\to K_1(B)$ be defined by
$$
\kappa=[\{\psi_t\}]\circ \boldsymbol{\bt}|_{K_0(A)}.
$$
Then $\kappa([1_A])=[u^*],$ since $u(t)\in U_0(B)$ for all $t\in
[0, \infty).$ This implies that $-[u]\in H_1(K_0(A), K_1(B)).$
Hence $[u]\in H_1(K_0(A), K_1(B)).$

\end{proof}

\begin{lem}\label{SL1}
Let $C$ be a unital AH-algebra and let $A$ be a unital separable
simple \CA\, with tracial rank zero. Suppose that $\phi_1, \phi_2:
C\to A$ be two monomorphisms such that there is a sequence of
unitaries $\{u_n\}\subset A$ such that the conclusions
\rm{((\ref{ltrbot-2})} and \rm{(\ref{ltrbot-3}))} of
\rm{\ref{LTRBOT}} hold. Then, in \ref{LTRBOT}, we may further
require that $u_n\in U_0(A),$ if $H_1(K_0(C), K_1(A))=K_1(A).$

\end{lem}

\begin{proof}
Let $x_n=[u_n]$ in $K_1(A).$ Then, since $K_1(A)=H_1(K_0(C),
K_0(A)),$ there is a \hm\, $\kappa_{n,0}: K_0(C)\to K_1(A)$ such
that $\kappa_{n,0}([1_C])=-x_n.$ Let $\kappa_{n,1}: K_1(C)\to
K_0(A)$ be zero map. By the Universal Coefficient Theorem, there
is $\kappa_n\in KK(C,A)$ such that
\beq\label{sl1-1}
(\kappa_n)|_{K_i(C)}=\kappa_{n,i},\,\,\,i=0,1.
\eneq
There is, for each $n,$  a positive number $\eta_n<\dt_n$ such
that
\beq\label{sl1-2}
{\rm ad}\, u_n\circ \phi_1\approx_{\eta_n}
\phi_2\,\,\,\text{on}\,\,\, {\cal F}_n.
\eneq
It follows from \ref{CCM} that there is a unitary $w_n\in U(A)$
such that
\beq\label{sl1-3}
\|[\phi_2(a), w_n]\|&<&(\dt_n-\eta_n)/2\rforal a\in {\cal
F}_n\andeqn\\
\text{Bott}(\phi_2, w_n)&=&\kappa_n.
\eneq
Put $v_n=u_nw_n,$ $n=1,2,....$ Then, we have
$$
{\rm ad}\,v_n\circ
\phi_1\approx_{\dt_n}\phi_2\,\,\,\text{on}\,\,\,{\cal F}_n\andeqn
$$
$$
\rho_A(\text{bott}_1(\phi_2, v_n^*v_{n+1}))=0\andeqn
[v_n]=[u_n]-x_n=0.
$$

\end{proof}

\begin{thm}\label{ST1}
Let $C$ be a unital AH-algebra and let $A$ be a unital separable
simple \CA. Suppose that $K_1(K_0(C),K_1(A))=K_1(A)$ and suppose
that $\phi_1, \phi_2: C\to A$ are two unital monomorphisms which
are asymptotically unitarily equivalent. Then there exists a
continuous path of unitaries $\{u(t): t\in [0, \infty)\}$ such
that
$$
u(0)=1\andeqn \lim_{t\to\infty}{\rm ad}u(t)\circ
\phi_1(a)=\phi_2(a)\rforal a\in C.
$$
\end{thm}

\begin{proof}
By \ref{NecT}, we must have
$$
[\phi_1]=[\phi_2]\,\,\,\text{in}\,\,\, KK(C,A),\,\,\, {\tilde
\eta}_{\phi_1,\phi_2}=0\andeqn
$$
$$
\tau\circ \phi_1=\tau\circ \phi_2.
$$

By \ref{SL1}, in the proof of \ref{TM}, we may assume that $v_n\in
U_0(A),$ $n=1,2,....$ It follows that $\xi_n([1_C])=0,$
$n=1,2,....$ Therefore $\kappa_n([1_C])=0.$ This implies that
$\gamma_n\circ \boldsymbol{\bt}([1_C])=0.$ Hence $w_n\in U_0(A).$
It follows that $u_n\in U_0(A).$ Therefore there is a continuous
path of unitary $\{U(t): t\in [-1,0]\}$ in $A$ such that
$$
U(-1)=1_A\andeqn U(0)=u(0).
$$
The theorem then follows.

\end{proof}

\begin{cor}\label{Scc2}
Let $C$ be a unital AH-algebra and let $A$ be a unital separable
simple \CA\, with tracial rank zero. Let $\phi: C\to B$ be a
unital monomorphism and let $u\in U(A).$ Then there exists a
continuous path of unitaries $\{U(t): t\in [0,\infty)\}\subset A$
such that
$$
U(0)=1_B\andeqn \lim_{t\to\infty}{\rm ad}\, U(t)\circ \phi(a)={\rm
ad}\, u\circ \phi(a)\rforal a\in C
$$
if and only if $[u]\in H_1(K_0(C), K_1(A)).$
\end{cor}

\begin{proof}
This follows from \ref{Sup} and the proof of \ref{ST1}.
\end{proof}

\begin{cor}\label{SCZ}
Let $C$ be a unital AH-algebra with $K_0(C)=\Z\cdot [1_C]\oplus G$
and let $B$ be a unital separable simple \CA\, with tracial rank
zero. Suppose that $\phi_1, \phi_2: C\to B$ are two unital
monomorphisms such that
\beq\label{scz1}
[\phi_1]&=&[\phi_2]\,\,\,\text{in}\,\,\,KK(C,B),\,\,\,{\tilde
\eta}_{\phi_1,\phi_2}=0\andeqn\\
\tau\circ \phi_1&=&\tau\circ \phi_2\rforal \tau\in T(B).
\eneq
Then there exists a continuous path of unitaries $\{u_t:t\in [0,
\infty)\}$ of $B$ such that
$$
u_0=1_B\andeqn \lim_{t\to\infty}{\rm ad}\, u_t\circ
\phi_1(a)=\phi_2(a)\rforal a\in A.
$$

\end{cor}

\begin{proof}
 Let $x\in K_1(B).$ Then one defines $\gamma: K_0(C)\to K_1(B)$
by $\gamma([1_C])=x$ and $\gamma|_G=0.$ This implies that
$$
H_1(K_0(C), K_1(B))=K_1(B).
$$
Thus the corollary follows from \ref{ST1}.
\end{proof}

\begin{cor}\label{SUC}
Let $X$ be a compact metric space and let $B$ be a unital
separable simple \CA\, with tracial rank zero. Suppose that
$\phi_1, \phi_2: C(X)\to B$ are two unital monomorphisms. Then
there exists a continuous path of unitaries $\{u(t): t\in
[0,\infty)\}\subset B$ such that
$$
u(0)=1_B \andeqn \lim_{t\to\infty}{\rm ad}\, u(t)\circ
\phi_1(a)=\phi_2(a)
$$
for all $a\in C(X)$ if and only if
\beq\label{suc1}
[\phi_1]&=&[\phi_2]\,\,\,\text{in}\,\,\,KK(C,B),\,\,\,{\tilde
\eta}_{\phi_1,\phi_2}=0\andeqn\\
\tau\circ \phi_1&=&\tau\circ \phi_2\rforal \tau\in T(B).
\eneq

\end{cor}

\begin{proof}
Note that $K_0(C(X))=\Z\cdot [1_C]\oplus G$ for some  subgroup
$G.$

\end{proof}

\begin{cor}\label{Sdiv}
Let $A$ be a unital AH-algebra and let $B$ be a unital separable
simple \CA\, with tracial rank zero such that $K_1(B)$ is
divisible. Suppose that $\phi_1, \phi_2: C\to B$ are two unital
monomorphisms such that
\beq\label{Sdvi1}
[\phi_1]&=&[\phi_2]\,\,\,\text{in}\,\,\,KK(C,B),\,\,\,{\tilde
\eta}_{\phi_1,\phi_2}=0\andeqn\\
\tau\circ \phi_1&=&\tau\circ \phi_2\tforal \tau\in T(B).
\eneq
Then there exists a continuous path of unitaries $\{u(t): t\in
[0,\infty)\}\subset B$ such that
$$u(0)=1_B \andeqn
\lim_{t\to\infty}{\rm ad}\, u(t)\circ \phi_1(a)=\phi_2(a)\tforal
a\in C.
$$

\end{cor}

\begin{proof}
Let $x\in K_1(B).$ Since $K_1(B)$ is divisible, a map from $[1_C]$
to $x$ can be extended to an element in $Hom(K_0(C), K_1(B)).$ So
$H_1(K_0(C), K_1(B))=K_1(B).$

\end{proof}

\begin{rem}\label{SR1}
{\rm Let $C$ be a unital AH-algebra and $A$ be a unital separable
simple \CA\, of tracial rank zero. Suppose that $\phi: C\to B$ is
a unital monomorphism. Denote by $ASU(\phi)$ the class of all
unital monomorphisms which are asymptotically unitarily equivalent
to $\phi.$ Let $\psi\in ASU(\phi).$ Then
$$
\psi(a)=\lim_{t\to\infty}{\rm ad}\, u_t\circ \phi(a)\rforal a\in A
$$
 for some continuous path of unitaries
$\{u_t: t\in [0, \infty)\}$ of $A.$ If there is another continuous
path of unitaries $\{v_t: t\in [0, \infty)\}$ of $A$ such that
$$
\psi(a)=\lim_{t\to\infty}{\rm ad}\, v_t\circ \phi(a)\rforal a\in
A,
$$
then,  as in the proof of \ref{Sup},
$$
[u_0v_0^*]\in H_1(K_0(C), K_1(A)).
$$
Define $u(\psi)$ to be the element in $K_1(A)/H_1(K_0(C), K_1(A))$
represented by $[u_0].$ Let $\psi_1, \psi_2$ be in $ASU(\phi).$
Then, by \ref{ST1}, $\psi_1$ and $\psi_2$ are strongly
asymptotically unitarily equivalent if and only if
$$
u(\psi_1)=u(\psi_2).
$$

Moreover, for any $x\in K_1(B),$ there is a unitary $u\in B$ such
that $x=[u].$ Define $\psi= {\rm ad}\, u\circ \phi.$ then
$u(\psi)=x.$

}

\end{rem}

\begin{df}\label{ann}
{\rm Let $A$ be a unital separable \CA. Denote by $\text{Ainn}(A)$
the group of all asymptotically inner automorphisms. A
automorphism $\af$ is said to be strong asymptotically inner if
there is a continuous path of unitaries $\{u(t): t\in [0,
\infty)\}$ of $A$ such that
$$
u(0)=1_A\andeqn \lim_{t\to\infty}u(t)^*au(t)=\af(a)\rforal a\in A.
$$
Denote by $\text{Ainn}_0(A)$ the subgroup of all strong
asymptotically inner automorphisms. }
\end{df}

By \ref{SR1}, we have the following:

\begin{thm}\label{ainn}
Let $A$ be a unital separable amenable simple \CA\, with tracial
rank zero satisfying the Universal Coefficient Theorem. Then
$$
{\text{Ainn}(A)\over{\text{Ainn}_0(A)}}\cong K_1(A)/H_1(K_0(A),
K_1(A)).
$$
\end{thm}

\end{document}